\documentstyle[11pt,amssymb,oldlfont]{article}

\font\tencmmib=cmmib10 \font\sevencmmib=cmmib7 \font\fivecmmib=cmmib5
\newfam\cmmibfam \textfont\cmmibfam=\tencmmib
\scriptfont\cmmibfam=\sevencmmib \scriptscriptfont\cmmibfam=\fivecmmib

 \font\teneufm=eufm10 at 12 pt
\font\seveneufm=eufm9 \font\fiveeufm=eufm7
\newfam\eufmfam \textfont\eufmfam=\teneufm
\scriptfont\eufmfam=\seveneufm \scriptscriptfont\eufmfam=\fiveeufm
\def\frak#1{{\fam\eufmfam\relax#1}}

\font\script=eusm10 at 12 pt
\font\sscrpt=eusm8

\font\helv=cmssbx10 at 12 pt
\def\nz{\hfill\break}
\def\isom{\ifmmode\ \cong\ \else$\isom$\fi}
\def\scsi{\scriptsize}
\def\scsis{\tiny}

\def\aa{\frak a}  \def\pp{\frak
    p}  \def\Gg{\frak g}
   \def\nn{\frak n} 
\def\tt{\frak t}  \def\ll{\frak
    l} \def\Ss{\frak s} \def\cc{\frak c}
\def\uu{\frak u}  \def\ee{\frak
    e}

\def\AA{\frak A} \def\BB{\frak B} \def\HH{\frak
    H} \def\DD{\frak D} 
\def\LL{\frak L} \def\JJ{\frak J}

   \def\CC{\frak C}

\def\fO{{\Bbb O}} 

 \def\fZ{{\Bbb Z}}
 \def\fC{{\Bbb C}}
\def\fP{{\Bbb P}} \def\fB{{\Bbb B}}
\def\fS{{\Bbb S}} \def\fQ{{\Bbb Q}}
\def\fR{{\Bbb R}} \def\fH{{\Bbb H}}

\def\bfb{\hbox{\scsi \bf b}}

\def\hbf#1{\hbox{\scsi\bf #1}}

\def\bs{\backslash}

\def\gD{\Delta} \def\gd{\delta} \def\gG{\ifmmode {\Gamma} \else$\gG$\fi}
\def\gg{\gamma} \def\gs{\sigma} \def\gS{\Sigma} 
 \def\gt{\theta}  \def\ga{\alpha}
\def\gb{\beta}  \def\gl{\lambda} \def\gL{\Lambda}
\def\gO{\Omega} \def\go{\omega}   \def\gr{\rho}
\def\grr{\varrho} \def\ge{\varepsilon}  

 \def\-#1{\overline{#1}}
\def\~#1{\tilde{#1}}    \def\cO{\ifmmode {\cal O} \else$\cO$\fi}

\def\cA{{\cal A} }

\def\cC{{\cal C} }
\def\cD{{\cal D} }

\def\cE{{\cal E} }

\def\cH{{\cal H} }

\def\cJ{{\cal J} }

\def\cL{{\cal L} }
\def\cM{{\cal M} }
\def\cN{{\cal N} }
\def\cO{{\cal O} }
\def\cP{{\cal P} }
\def\cR{{\cal R} }

\def\cS{{\cal S} }
\def\cT{{\cal T} }
\def\cU{{\cal U} }

\def\cV{{\cal V} }

\def\cZ{{\cal Z} }

\def\scC{\hbox{{\script C}} }

\def\scH{\hbox{{\script H}} }

\def\scL{\hbox{{\script L}} }
\def\scM{\hbox{{\script M}} }

\def\scQ{\hbox{{\script Q}} }
\def\scS{\hbox{{\script S}} }
\def\scT{\hbox{{\script T}} }

 \def\ad{\hbox{ad}} \def\dim{\hbox{dim}}
  
\def\Aut{\hbox{Aut}}  
  \def\Ker{\hbox{Ker}}
 \def\rank{\hbox{rank}}

\def\isom{\ifmmode\ \cong\ \else$\isom$\fi}

\def\nni{\supset}  \def\ovv{\ifmmode
  \overline{\gG \bs \cal D}\ \else$\ovv$\fi}

\def\xg{\ifmmode {X_{\gG}} \else$\xg$\fi} \def\xgeq{\ifmmode {\xg =
    \gG\bs\cD} \else$\xgeq$\fi}
\def\xgs{\ifmmode {X_{\gG}^*} \else$\xgs$\fi}
\def\xgc{\ifmmode {\overline{X}_{\gG}} \else$\xgc$\fi}

\def\cIf{\ifmmode {\cal I}_5%
  \else$\cIf$\/\fi} \def\cIt{\ifmmode {\cal I}_{10} \else$\cIt$\fi}
\def\JG{\ifmmode {G_{25,920}} \else$\JG$\fi}

\def\pts{\ifmmode \hbox{\helv P}^2_6 \else $\pts$\fi} \def\pthrees{\ifmmode
  \hbox{\helv P}^3_6 \else $\pthrees$\fi}

\def\pos{\ifmmode \hbox{\helv P}^1_6 \else $\pos$\fi}

\def\hra{\hookrightarrow} 
\def\lra{\longrightarrow} 
\def\lla{\longleftarrow} 
\def\llra{\longleftrightarrow} \def\ra{\rightarrow} \def\Ra{\Rightarrow}
 \def\bs{\ifmmode {\setminus} \else$\bs$\fi}

\def\inn{\subset} \def\nni{\supset} 
\def\ove{\overline}

\def\ende{\hfill $\Box$ \vskip0.25cm }

\def\p{\prime} \def\sdprod{\rtimes}
\def\boldone{\hbox{\Large\boldmath $1_{\hbf{n}}$}}
\def\boldOne{\hbox{\Large\bf 1}}

\def\hbf#1{\hbox{\scsi\bf #1}}
\def\hbfs#1{\hbox{\scsis\bf #1}}

\newtheorem{theorem}{Theorem}[section]

\newtheorem{lemma}[theorem]{Lemma}
\newtheorem{proposition}[theorem]{Proposition}
\newtheorem{corollary}[theorem]{Corollary}

\newenvironment{definition}{\refstepcounter{theorem}
\vskip.2cm\noindent{\bf Definition \thesection.\arabic{theorem}
  }}{\vskip.2cm\noindent}

\newenvironment{example}{\refstepcounter{theorem}
\vskip.2cm\noindent{\bf Example \thesection.\arabic{theorem}
  }}{\vskip.2cm\noindent}

\newenvironment{remark}{\refstepcounter{theorem}
\vskip.2cm\noindent{\bf Remark \thesection.\arabic{theorem}
  }}{\vskip.2cm\noindent}

\begin{document}

\title{Modular subvarieties of arithmetic quotients of bounded symmetric
  domains}
\author{Bruce Hunt \\ FB Mathematik, Universit\"at \\ Postfach 3049 \\
  67653 Kaiserslautern}
\maketitle

\begin{center}
{\Large\bf Introduction}
\end{center}
A reductive $\fQ$-simple algebraic group $G$ is of {\it hermitian type}, if
the symmetric space $\cD$ defined by $G(\fR)$ is a hermitian symmetric
space. A discrete subgroup $\gG\inn G(\fQ)$ is {\it arithmetic}, if for
some faithful rational representation $\gr:G\lra GL(V)$, and for some
lattice $V_{\fZ}\inn V_{\fQ}$, $\gG$ is commensurable with
$\gr^{-1}(GL(V_{\fZ}))$. A non-compact hermitian symmetric space is
holomorphically equivalent to a bounded symmetric domain. If this is the
case, $\gG$ acts on $\cD$ preserving the natural Bergmann metric, and
$X_{\gG}=\gG\bs \cD$ is, if $\gG$ is torsion-free, a complex manifold, in
general not compact (it is compact exactly when $G$ is anisotropic). We
call spaces $X_{\gG}$ arithmetic quotients of bounded symmetric domains
even when $\gG$ has torsion; it is known that in this case $X_{\gG}$ is a
$V$-manifold (in the sense of Satake), locally the quotient of a smooth
space by a finite group action. There is a natural compactification of
$X_{\gG}$, the Satake compactification $\xgs$, which has the property: the
complement $\xgs-X_{\gG}$ is a finite disjoint union of arithmetic quotients of
bounded symmetric domains of lower dimension. A natural problem in this
respect is to consider, in addition to the data above, a symmetric
subdomain $\cD'\inn \cD$ which has the property that the restriction of the
action of $\gG$ to $\cD'$ is discrete, say by a discrete subgroup $\gG'$,
resulting in a commutative square
$$\begin{array}{ccc}\cD' & \hra & \cD \\ \downarrow & & \downarrow \\
  X_{\gG'}=\gG'\bs \cD' & \hra & \gG\bs \cD = X_{\gG}.
\end{array}$$
It is spaces arising as $X_{\gG'}$ that we will refer to as {\it modular
  subvarieties}. This situation ensuses in particular if $N\inn G$ is a
reductive subgroup of hermitian type; then the symmetric space $\cD_N$ of
$N$ has a holomorphic symmetric embedding $\cD_N\inn \cD$. As a matter of
notation, we refer to such subgroups as {\it symmetric} subgroups.
This explains the title and describes the topic of this
paper.

We are in fact concerned with a very special set of modular subvarieties
which have very special behavior at the cusps. Upstairs in the universal
covers (in the domains $\cD'$ and $\cD$), this behavior defines a notion of
incidence, implying a relation between (real) parabolic and symmetric
subgroups. Conversely, one can define a notion of incidence, group
theoretically, between parabolic and symmetric subgroups $P$ and $N$, which
implies the desired behavior of the subdomain $\cD_N$ (the $\cD'$ in the
notations above) near the cusp $F$ which corresponds to $P$. It is the
geometric point of view we will adapt in this paper. Starting with a
rational boundary component $F$ (respectively the corresponding parabolic
$P$), we define what it means for $F$ and a symmetric subdomain $\cD_N$
(respectively the parabolic $P$ and a symmetric subgroup $N$) to be
incident. For $F$ and $\cD_N$ this is easy to formulate: $\cD_N$ and $F$
are incident if $F$ is a rational boundary component of $\cD_N$ and maximal
with this property. For the
subgroups $N$ and $P$ the condition is more complicated to formulate, but
consists essentially in a maximality condition plus the condition
$N=N_1\times N_2$, and $N_1$ is a hermitain Levi factor of $P$ (if
$\dim(F)>0$) or just that $N$ is irreducible such that $F$ is a boundary
component of $\cD_N$ (if $\dim(F)=0$), see below for details.

This gives rise to a notion of incidence for $\fQ$-subgroups $P$ and $N$,
by insisting that incidence holds as above for the groups of real points
$P(\fR)$ and $N(\fR)$. Everything being defined over $\fQ$, one can proceed
to form an arithmetic quotient by an arithmetic subgroup $\gG\inn G(\fQ)$,
and $\gG$ will act naturally both on $F$ and on $\cD_N$. The quotients will
be a boundary variety $W$ and a modular subvariety $X_{\gG'}$ as above, and
the notion of incidence now becomes $W\inn X_{\gG'}^*$, where $X_{\gG'}^*\inn
\xgs$ is the embedding of Baily-Borel embeddings described in
\cite{S2}. Thus the problem has two aspects
\begin{itemize}\item[a)] The question of the {\it existence} of
  $\fQ$-subgroups $N$ which are incident with a given $P$, and
\item[b)] the action of $\gG$ on the set of such subgroups.
\end{itemize}
The first point a) was treated in detail in \cite{hyp} and \cite{sym};
this paper contains a preliminary study of b).

Similar questions have been asked and answered for parabolic
subgroups, but the situation for symmetric subgroups is totally different. The
basic cause for this is the following elementary fact: if two
$\fQ$-parabolics are conjugate (over the algebraic closure), then they are
in fact $\fQ$-conjugate, i.e., conjugate by an element of $G(\fQ)$, whereas the
corresponding statement for symmetric subgroups is totally wrong. The
reason behind this is the very fundamental
fact that the homogenous spaces $G/P$ are {\it projective} for parabolic
subgroups $P$, while $G/N$ is {\it affine}
for a reductive subgroup $N$, a property
which in fact characterizes reductive subgroups $N$, as was proved
in \cite{BHC}. A consequence of this is that it holds
for the parabolic subgroup $P$ that $G(\fQ)/P(\fQ)=(G/P)(\fQ)$,
while for symmetric subgroups $N$ only the inclusion $G(\fQ)/N(\fQ)\inn
(G/N)(\fQ)$ holds. The elements of the left-hand side are those subgroups
$N'(\fQ)\inn G(\fQ)$ which are $G(\fQ)$-conjugate to $N(\fQ)$, while on the
right-hand side we have those $G$-conjugates of $N$ which are {\it defined
  over $\fQ$}. Considering only the former leads to the definition of {\it
  rational} symmetric subgroups (Definition \ref{D4.1}).
{}From this point on we restrict our attention
to symmetric subgroups $N$ which are incident with a given parabolic $P$
(or the other way around, parabolic subgroups $P$ which are incident with a
given symmetric subgroup $N$). This means we consider for each $b\in
\{1,\ldots, s\},\ s=\rank_{\fQ}G$, a fixed isomorphism class (with the
exception $b=s$ for the two exceptional domains, for which there are three,
resp. two such isomorphism classes) of the groups of real points $N(\fR)$
for subgroups $N$ (conjugate to a fixed symmetric subgroup $N_{\bfb}$
incident with $P_{\bfb}$, a standard parabolic),
and consider pairs $(N,P)$ of incident
symmetric and parabolic subgroups. We consider rational symmetric subgroups
$N$ representing points of $G(\fQ)/N_{\bfb}(\fQ)$, i.e., $N$ is
$G(\fQ)$-conjugate to
$N_{\bfb}$. However, these turn out to still be too many such subgroups in
the sense that there are in general infinitely many $\gG$-equivalence
classes for any $\gG$. We are led to introduce the notion of $\gG$-{\it
  integral} symmetric subgroups (Definition \ref{D5.1}),
of which there {\it are} finitely many
$\gG$-orbits. To formulate this we assume $\dim(F)>0$ and let
$N_{\bfb}=N_{\hbf{b,1}}\times N_{\hbf{b,2}}$ (resp. $N=N_1\times N_2$) be
the decomposition of $N_{\bfb}$ (resp. of $N$) mentioned above. Then the
definition in this case is:
\[ \hbox{ $N=gN_{\bfb}g^{-1}$ is integral }
\iff g N_1 g^{-1}\cap \gG = g(N_{\hbf{b,1}}\cap \gG) g^{-1},\]
and it depends in fact on the choice of $N_{\bfb}$ (that is, on the choice
of maximal $\fQ$-split torus and order on it) as well as on $\gG$. (For
$\dim(F)=0$, the condition is simply $gNg^{-1}\cap \gG = g(N_{\bfb}\cap
\gG)g^{-1}$.) Fixing
this data $N_{\bfb}$ and $\gG$ leads to a finite number of
$\gG'$-equivalence classes of
$\gG$-integral symmetric subgroups conjugate to $N_{\bfb}$ for any
arithmetic group $\gG'\inn G(\fQ)$,
which we show by utilizing the basic
finiteness result of \cite{BHC}.

We now consider $\gG$-integral symmetric subgroups $N$ and arbitrary
arithmetic subgroups $\gG'\inn G(\fQ)$, let $\gG'_N=N\cap
\gG'$ and consider the corresponding modular subvarieties they define,
$X_{\gG'_N}\inn  X_{\gG'}$; we call these {\it integral} modular
subvarieties (Definition \ref{D9.2a}). As mentioned above, the inclusion
extends to the Baily-Borel embeddings $X_{\gG'_N}^*\inn X_{\gG'}^*$. We now
take $\gG$ to be $G_{\fZ}$ for some rational representation $\gr:G\lra
GL(V)$, that is $\gG=\gr^{-1}(GL(V_{\fZ}))$ for some $\fZ$-structure
$V_{\fZ}$ on $V$. Putting all the pieces together yields our main result,
which we now formulate. For this we refer to the notations $\nu_b(\gG'),\
b=1,\ldots, s$ and $\mu_b(\gG,\gG'),\ b=1,\ldots, s$ of Definition \ref{D6.1}
and \ref{D7.1}, respectively, for the number of $b^{th}$ boundary varieties
and the number of $b^{th}$ integral modular subvarieties, respectively. We
let $W_{b,i},\ b=1,\ldots, s,\ i=1,\ldots,\nu_b(\gG')$ be the corresponding
boundary varieties on the Satake compactification, $Y_{b,j},\ b=1,\ldots,
s\ j=1,\ldots, \mu_b(\gG,\gG')$ the corresponding $\gG$-integral modular
varieties, everything on the arithmetic quotient $X_{\gG'}$. Then
the main result of the paper is the following.
\begin{theorem} Let $\gG$ be as above, $\gG'\inn G(\fQ)$ arithmetic, and
  $X_{\gG'}\inn X_{\gG'}^*$ the Satake compactification,
  $X_{\gG'}^*-X_{\gG'}=\sum_{b,i}W_{b,i}$. Then $\Xi:=\sum_{b,j}Y_{b,j}$ is
  a complete (finite, non-empty) set of $\gG'$-equivalence classes of
  $\gG$-integral modular subvarieties, such that for each $W_{b,i}$, there
  is at least one $Y_{b,j}$ incident to $W_{b,i}$.
\end{theorem}

This gives us a {\it well-defined, non-empty, finite} set of subvarieties
of the Baily-Borel embedding $X_{\gG'}^*\inn \fP^N$ for any subgroup
$\gG'\inn \gG$ of finite index. Furthermore these have a prescribed
behavior near the cusps. For example, if $f:\cD\lra \fC$ is a modular form
whose zero divisor $D_f$ on $X_{\gG'}^*$ contains the union of the integral
modular subvarieties, then $f$ is a cusp form for $\gG'$.

In the case that the arithmetic quotient $\xg$ is a moduli space of abelian
varieties with some {\sc Pel} structure, it is not difficult to see the
moduli interpretation of the modular subvarieties $X_{\gG_N}$ determined by
{\it rational} symmetric subgroups. See
Example \ref{example}, where an interesting case is discussed in more
detail. That is the case $G=Sp(4,\fQ)$, and it turns out that in this case
many modular subvarieties which come from $\fQ$-groups $N$ conjugate
to the standard one have a nice moduli interpretation: they parameterize
abelian varieties (surfaces in this case) with real multiplication by a
real quadratic field $k$. The group $N$ is $G(\fQ)$-conjugate by the
standard one $N_{\hbf{1}}$ precisely when the field $k$ splits into a
product $\fQ\times \fQ$, and in this case the real multiplication
``degenerates'' into two copies of multiplication by $\fQ$, in other words
the abelian surface is no longer simple but {\it isogenous} to a {\it
  product}. In this particular case the set of $\gG$-integral subgroups
($\gG=Sp(4,\fZ)$) corresponds to the abelian surfaces which are actually
isomorphic to a product. We prove more generally the following:
\begin{theorem} Let $G,\ S,\ P_{\bfb},\ N_{\bfb}$ and $\gG$
  be as above ($b<t$), $\gG'\in G(\fQ)$ arithmetic, and let $X_{\gG'_N}$ be a
  modular subvariety of $X_{\gG'}$ for $N$ rational symmetric, conjugate to
  $N_{\bfb}$. Then $X_{\gG'_N}$ is a finite quotient of a product, and the
  set of $\gG'$-equivalence classes of such modular subvarieties forms a
  locus in $X_{\gG'}$ where the corresponding abelian varieties are
  isogenous to products, i.e., are not simple.
  If $N$ is $\gG$-integral, then $X_{\gG'_N}$ is a product, and
  the set of $\gG'$-equivalence classes of such modular subvarieties forms
  a locus in $X_{\gG'}$ where the corresponding abelian varieties split
  while preserving the endomorphisms (but not necessarily the polarizations).
\end{theorem}

An important application of the main result is to define a simplicial
complex which is an analog for reductive groups of what the Tits building
is for parabolic groups. Recall this complex $\cT(G)$
is constructed by forming the
simplicial complex of the set of all rational parabolics, partially ordered
by the inverse of the inclusion. Taking $\gG$-equivalence classes of the
parabolics gives rise to a {\it finite} quotient complex $\cT(G)/\gG$,
whose vertices are in one to one correspondence with the boundary varieties
of $\xgs$. However, the Tits building $\cT(G)$ is a {\it rational}
invariant of $G$, not depending in any way on the arithmetic group
$\gG$. As an analog we can define a complex $\cS(\gG)$ by replacing the
maximal $\fQ$-parabolics by the $\gG$-integral symmetric subgroups incident
with the maximal $\fQ$-parabolics. The arithmetic group $\gG$ acts on
$\cS(\gG)$, and the quotient complex $\cS(\gG)/\gG$ is again finite, by our
main result above, and the vertices are in one to one correspondence with
the $\gG$-integral modular subvarieties which are incident with a rational
boundary variety. In this case the complex $\cS(\gG)$ itself
is an {\it integral}, not a
rational, object. This will be dealt with elsewhere.

We now sketch the contents of the paper. The first paragraph is preliminary
and recalls some mostly known facts on the classification of the rational
groups of hermitian type for which we could find no reference.
In the second we recall some notions and results
from \cite{sym} on which the theory is based. The third papagraph
essentially describes all arithmetic subgroups of the rational groups of
hermitain type, in terms of maximal orders and ideals therein in division
algebras $D$ for the classical cases and in terms of maximal orders of
exceptional Jordan algebras for the exceptional cases. Theoretically this
section could have been dispensed with, but it does help one get the
feeling for the modular subvarieties later on. In the fourth paragraph we
introduce the notion of $\gG$-integral symmetric subgroups and derive the
basic finiteness result. In the fifth paragraph we discuss the
compactifications of the arithmetic quotients, and finally in the sixth
paragraph we define precisely modular subvarieties and derive the main
results above.

\noindent{\bf Thanks} I would like to acknowledge helpful discussions with
Steven Weintraub which led to the definition of $\gG$-integral, which is
more or less the central contribution of the paper.

\noindent {\bf Notations:} For an algebraic $k$-group $H$, the group of
$K$-valued points for a field extension $K|k$ will be denoted $H_K$ or
$H(K)$; similarly, for a vector space $V$ defined over $k$, $V_k$ will
denote the set of $k$-points, and for the ring of integers $\cO_k$,
$V_{\cO_k}$ will denote a $\cO_k$-lattice in $V_k$.
Throughout, $s$ will denote the $\fQ$-rank of a $\fQ$-group $G$, $t$ will
denote the $\fR$-rank of $G(\fR)$, and $f$ will denote the degree of $k$
over $\fQ$, when $k$ is a totally real number field fixed in a
discussion. Usually $d$ will denote the degree of a division algebra $D$,
and $n$ will denote the dimension of a $D$-vector space $V$. For a group
$G$ and a subset $\Xi\inn G$, the normalizer (resp. centralizer) of $\Xi$
in $G$ will be denoted $\cN_G(\Xi)$ (resp. $\cZ_G(\Xi)$).
\tableofcontents

\section{Rational groups of hermitian type}
\subsection{Notations}
We now fix some notations to be in effect for the rest of the paper. We
will be dealing with algebraic groups defined over $\fQ$, which give rise
to hermitian symmetric spaces, groups of {\it hermitian type}, as we will
say. As we are interested in the automorphism groups of domains, we may,
without restricting generality, assume the group is {\it centerless}, and
{\it simple} over $\fQ$. We will also assume $G$ is Zariski connected.
Henceforth, if not indicated otherwise (occasionally $G$ will denote a
reductive group; in sections 2.1 and 2.2 $G$ will be a real Lie group) $G$
will denote such an algebraic group. To avoid complications, we exclude in
this paper the following case:

\vspace*{.2cm}{\bf Exclude:}\hspace*{2cm} All non-compact real factors of
$G(\fR)$ are of type $SL_2(\fR)$.

\vspace*{.2cm}
\noindent Finally, we shall
only consider {\it isotropic} groups. This implies the hermitian symmetric
space $\cD$ has no compact factors. By our assumptions, then, we have
\begin{itemize}\item[(i)] $G=Res_{k|\fQ}G'$, $k$ a totally real number
  field, $G'$ absolutely simple over $k$.
\item[(ii)] $\cD=\cD_1\times \cdots \times \cD_f$, each $\cD_i$ a non-compact
  irreducible hermitian symmetric space, $f=[k:\fQ]$.
\end{itemize}
\subsubsection{Real parabolics}
This material is presented in detail in \cite{BB} and \cite{sym}, 1.2, so
we just mention enough to fix notations. We work in this section in the
category of real Lie groups. $G$ will denote a
connected reductive real Lie group of hermitian type, such that the
symmetric space $\cD=G/K$ is irreducible. In a well-known manner one fixes
a maximal set of strongly orthogonal (absolute) roots, defining a
subalgebra $\aa\inn \Gg$, such that $A=\exp(\aa)$ is a maximal $\fR$-split
torus which will be fixed throughout this discussion. The set of strongly
orthogonal roots is ordered, defining an order on $A$, which determines a
set of simple $\fR$-roots $\gD_{\fR}=\{\eta_1,\ldots, \eta_t\},\
t=\rank_{\fR}G=\dim(A)$, in the $\fR$-root system
$\Phi_{\fR}:=\Phi(\aa,\Gg)$. For each $b\in \{1,\ldots,t\}$, the
one-dimensional subtorus $A_b$ is defined: $\aa_b=\bigcap\limits_{i\neq
  b}\Ker(\eta_i),\ A_b=\exp(\aa_b)$. We also set $\nn=\sum\limits_{\eta\in
  \Phi_{\fR}^+} \Gg^{\eta},\ N=\exp(\nn)$.
The {\it standard maximal
  $\fR$-parabolic}, $P_b,\ b=1,\ldots, t$, is the group generated by
$\cZ_G(A_b)$ and $N$; equivalently it is the semidirect product (Levi
decomposition)
\begin{equation}\label{e2.2} P_b=\cZ_G(A_b)\rtimes U_b,
\end{equation}
where $U_b$ denotes the unipotent radical.
For real parabolics of hermitian type one has a very useful refinement of
(\ref{e2.2}). This is explained in detail in \cite{SC} and especially in
\cite{S}, \S III.3-4. First we have the decomposition of $\cZ_G(A_b)$
\begin{equation}\label{e3.4} \cZ_G(A_b)=M_b\cdot L_b \cdot \cR_b,
\end{equation}
where $M_b$ is compact, $L_b$ is the {\it hermitian Levi factor}, $\cR_b$ is
reductive (of type $\bf A_{\hbox{\scsi \bf b-1}}$), and the product is almost
direct (i.e., the factors have finite intersection). Secondly, the
unipotent radical decomposes,
\begin{equation}\label{e3.5} U_b=\cZ_b\cdot V_b,
\end{equation}
which is a direct product, $\cZ_b$ being the center of $U_b$. For this
decomposition the groups are both Zariski connected and connected in the
real Lie groups. The action of
$\cZ_G(A_b)$ on $U_b$ can be explicitly described, and is the basis for the
compactification theory of \cite{SC}.
The main results can be found in \cite{S}, III \S3-4, and can be summed up
as follows.
\begin{theorem}\label{t4.1} In the decomposition of the standard parabolic
  $P_b$ (see (\ref{e3.4}) and (\ref{e3.5}))
$$P_b=(M_b\cdot L_b \cdot \cR_b)\rtimes \cZ_b\cdot V_b,$$
the following statements hold.
\begin{itemize}\item[(i)] The action of $M_b\cdot L_b$ is trivial on
  $\cZ_b$, while on $V_b$ it is by means of a symplectic representation
  $\gr:M_b\cdot L_b \lra Sp(V_b,J_b)$, for a symplectic form $J_b$ on
  $V_b$.
\item[(ii)] $\cR_b$ acts transitively on $\cZ_b$ and defines a homogenous
  self-dual (with respect to a bilinear form) cone $C_b\inn \cZ_b$, while
  on $V_b$ it acts by means of a representation $\gs:\cR_b\lra GL(V_b,I_b)$
  for some complex structure $I_b$ on $V_b$.
\end{itemize}
Furthermore the representations $\gr$ and $\gs$ are compatible in a natural
sense. The decomposition and the representations in fact are valid for the
corresponding real algebraic group $G$ and its algebraic subgroups.
\end{theorem}
Finally, there is a one to one correspondence between the maximal real
parabolics $P$ (each of which is conjugate to a unique $P_b$) and the
boundary components $F$ (each of which is the image of a unique standard
boundary component $F_b$), given by $P\llra F$, where $P=\cN_G(F)$. In
particular, $P_b=\cN_G(F_b)$.

\subsubsection{Roots}
We now return to the notation used above, $G=Res_{k|\fQ}G'$ the
$\fQ$-simple group of hermitian type.
and introduce a few notations concerning the root systems involved. Let
$\gS_{\infty}$ denote the set of embeddings $\gs:k\lra \fR$; this set is in
bijective correspondence with the set of infinite places of $k$. We
denote these places by $\nu$, and if necessary we denote the
corresponding embedding by
$\gs_{\nu}$. For each $\gs\in \gS_{\infty}$, the group $^{\gs}G'$ is the
algebraic group defined over $\gs(k)$ by taking the
set of elements $g^{\gs},\ g\in G'$. For each infinite
prime $\nu$ we have
$G_{k_{\nu}}\cong (^{\gs_{\nu}}G')_{\fR}$, and the decomposition of $\cD$
above can be written
$$\cD=\prod_{\gs\in \gS_{\infty}}\cD_{\gs},\quad
\cD_{\gs}:=(^{\gs}G')_{\fR}/K_{(\gs)}=(^{\gs}G')_{\fR}^0/K_{(\gs)}^0.$$
We set $G_{\gs}=({^{\gs}G}')^0_{\fR}$ and note that the discussion of the
last section applies to $G_{\gs}$ for each $\gs$. For convenience we now
index the components $\cD_{\gs}$ by $i\in \{1,\ldots, f\}$.
For each $\cD_i$ we have $\fR$-roots $\Phi_{i,\fR}$,
of $\fR$-ranks $t_i$ and simple $\fR$-roots
$\{\eta_{i,1},\ldots,\eta_{i,t_i}\},\ i=1,\ldots, d$. For each factor we
have standard parabolics $P_{i,b_i}$ $(1\leq b_i\leq t_i$) and standard
boundary components $F_{i,b_i}$. The standard parabolics of $G_{\fR}^0$ and
boundary components of $\cD$ are then
products
\begin{equation}\label{e3.2} P_{\bfb}(\fR)^0=P_{1,b_1}\times \cdots \times
  P_{d,b_d},\quad F_{\bfb}=F_{1,b_1}\times \cdots \times
  F_{d,b_d},\quad ({\bf b}=(b_1,\ldots, b_d)),
\end{equation}
where $P_{i,b_i}\inn G_{\gs_i}$, $P_{\bfb}\inn G$ is a maximal
$\fQ$-parabolic, and as above $P_{\bfb}(\fR)^0=\cN_{G_{\fR}}(F_{\bfb})^0$.
Furthermore, there is
a $\fQ$-subgroup $L_{\bfb}\inn G$ such that
\begin{equation}\label{e3.3} \Aut(F_{\bfb})^0=L_{\bfb}(\fR)^0,\
  L_{\bfb}(\fR)^0=L_{1,b_1}\times
  \cdots\times L_{d,b_d},
\end{equation}
where $L_{i,b_i}\inn P_{i,b_i}$ is the hermitian Levi component as above.
As far as the domains are concerned, any of the boundary components
$F_{i,b_i}$ may be the {\it improper} boundary component $\cD_i$,
which is indicated by setting $b_i=0$. Consequently, $P_{i,0}=L_{i,0}=G_i$
and in (\ref{e3.2}) and (\ref{e3.3})
any ${\bf b}=(b_1,\ldots,b_d),\ 0\leq b_i\leq t_i$ is admissible.

Since $G'$ is isotropic, there is a positive-dimensional $k$-split torus
$S'\inn G'$, which we fix. Then  ${^{\gs}S}'$ is a maximal $\gs(k)$-split
torus of $^{\gs}G'$ and there is a canonical isomorphism $S'\ra {^{\gs}S}'$
inducing an isomorphism $\Phi_k=\Phi(S',G')\lra
\Phi_{\gs(k)}({^{\gs}S}',{^{\gs}G}')=:\Phi_{k,\gs}$.
The torus $Res_{k|\fQ}S'$ is
defined over $\fQ$ and contains $S$ as maximal $\fQ$-split torus; in fact
$S\cong S'$, diagonally embedded in $Res_{k|\fQ}S'$. This yields an
isomorphism $\Phi(S,G)\cong \Phi_k$, and the root systems
$\Phi_{\fQ}=\Phi(S,G)$, $\Phi_k$ and $\Phi_{k,\gs}$ (for all $\gs\in
\gS_{\infty}$) are identified
by means of the isomorphisms.

In each group $^{\gs}G'$ one chooses a maximal $\fR$-split torus
$A_{\gs}\nni {^{\gs}S}'$, contained in a maximal torus defined over
$\gs(k)$. Fixing an order on $X(S')$ induces one also on $X({^{\gs}S}')$
and $X(S)$. Then, for each $\gs$, one chooses an order on $X(A_{\gs})$
which is compatible with that on $X({^{\gs}S}')$, and $r:X(A_{\gs})\lra
X({^{\gs}S}')\cong X(S)$ denotes the restriction homomorphism. The canonical
numbering on $\gD_{\fR,\gs}$ of simple $\fR$-roots of $G$ with respect to
$A_{\gs}$ is compatible by restriction with the canonical numbering of
$\gD_{\fQ}$ (\cite{BB}, 2.8). Recall also that each $k$-root in $\Phi_k$ is
the restriction of at most one simple $\fR$-root of $G'(\fR)$ (which is a
simple Lie group). Let $\gD_k=\{\gb_1,\ldots,\gb_s\}$; for $1\leq i\leq s$ set
$c(i,\gs)$:= index of the simple $\fR$-root of $^{\gs}G'$ restricting on
$\gb_i$. Then $i<j$ implies $c(i,\gs)<c(j,\gs)$ for all $\gs\in \gS$.

Each simple $k$-root defines a unique standard boundary component: for
$b\in \{1,\ldots,s\}$,
\begin{equation}\label{e9.1} P_{\bfb}:=\prod_{\gs\in \gS_{\infty}}
  P_{c(b,\gs)} \quad (\hbox{ resp. }F_{\bfb}:=\prod_{\gs\in
    \gS_{\infty}}F_{c(b,\gs)}),
\end{equation}
which is the product of standard (with respect to $A_{\gs}$ and
$\gD_{\fR,\gs}$) parabolics $P_{c(b,\gs)}\inn G_{\gs}$ (resp.
boundary components $F_{c(b,\gs)}$ of $\cD_{\gs}$). It
follows that $\overline{F}_{\hbox{\scsi\bf j}}\inn
\overline{F}_{\hbox{\scsi\bf i}}$ for $1\leq i\leq j\leq
s$. Furthermore, setting $o_{\bfb}:=\prod o_{c(b,\gs)}$, then (\cite{BB},
p.~472)
\begin{equation}\label{e9.2} F_{\bfb}=o_{\bfb}\cdot L_{\bfb}(\fR)^0,
\end{equation}
where $L_{\bfb}(\fR)^0$ denotes the hermitian Levi component (\ref{e3.3}).
As these are the only boundary components of
interest to us, we will henceforth refer to any conjugates of the
$F_{\bfb}$ of (\ref{e9.1}) by elements of $G$ as {\it rational boundary
  components} (these should more precisely be called rational with respect
to $G$), and to the conjugates of the parabolics $P_{\bfb}$ as the {\it
  rational parabolics}.

\subsection{Classification}\label{classification}
For the convenience of the reader we sketch the classification of rational
groups of hermitian type. As $G=Res_{k|\fQ}G'$ for an absolutely simple
$G'$ over $k$ we need only classify these.
\subsubsection{Classical cases}
By means of the correspondence
given by Weil in \cite{W} if $G$ is of classical type, classifying the
(semi)simple $k$-groups of interest to us is equivalent to
classifying the central (semi)simple $k$-algebras with involution such that
$\Aut(A,*)$ is of hermitian type. We now just list the cases, the
possible $k$-groups $G'$, the set of $\fR$-points of $G'$ as well as of the
$\fQ$-group $G$, and the corresponding domains. Let $k$ be a totally real
number field of degree $f$ over $\fQ$. For the bounded symmetric domains we
shall use the notations $\bf I_{\hbf{p,q}},\ II_{\hbf{n}},\ III_{\hbf{n}},\
IV_{\hbf{n}},\ V,\ VI$. In what follows $G'$ will be simple but not
necessarily centerless.
\begin{itemize}\item[{\bf O}] Orthogonal type
\begin{itemize}\item[{\bf O.1}] split case: $G'=SO(V,h)$, $V$ a $k$-vector
  space of dimension $n+2$, $h$ a symmetric bilinear form such that,
  at all real primes $\nu$, $h_{\nu}$ has signature $(n,2)$.
  $$G'(\fR)\cong SO(n,2),\quad
  G(\fR)\cong \prod^f_{i=1}SO(n,2)_i,\quad \cD\cong {\bf IV_{\hbf{n}}}\times
  \cdots \times {\bf IV_{\hbf{n}}},\ f\hbox{ factors}.$$
\item[{\bf O.2}] non-split case: $G'=SU(V,h)$, $V$ a right $D$-vector space
  of dimension $n$, $h$ is a skew-hermitian form; here $D$ is a
  quaternion division algebra, central simple over $k$, and for
  all real primes, either
  \begin{itemize}\item $D_{\nu}\cong \fH,\ \ G_{\nu}'\cong SU(\fH^n,h)$,
    $h$ a skew-hermitian form on $\fH^n$,
                 \item $D_{\nu}\cong M_2(\fR),\ \ G_{\nu}'\cong SO(2n-2,2).$
  \end{itemize}
and in the first case $h$ has Witt index $[{n\over 2}]$.
Number the primes such that for $\nu_1,\ldots, \nu_{f_1}$ the first case
occurs and for $\nu_{f_1+1},\ldots, \nu_f$ the second occurs. Then
$G'(\fR)\cong SU(\fH^n,h)$, and
$$G(\fR)\cong \prod^{f_1}_{i=1}SU(\fH^n,h)_i\times \prod_{i=f_1+1}^f
SO(2n-2,2)_i,\quad \cD\cong \underbrace{\bf II_{\hbf{n}}\times \cdots
    \times {\bf
    II_{\hbf{n}}}}_{f_1\ \hbox{\small factors}}\times
    \underbrace{\bf IV_{\hbf{2n-2}}\times \cdots \times {\bf
    IV_{\hbf{2n-2}}}}_{f-f_1\ \hbox{\small factors}}.$$
\end{itemize}
\item[{\bf S}] Symplectic type
\begin{itemize}\item[{\bf S.1}] split case: $G'=Sp(2n,k),\ G'(\fR)\cong
  Sp(2n,\fR),\ G(\fR)\cong
  (Sp(2n,\fR))^f,\ \cD\cong ({\bf III_{\hbf{n}}})^f.$
\item[{\bf S.2}] non-split case: $G'=SU(V,h)$, where $V$ is an
  $n$-dimensional right vector space over a quaternion division algebra
  $D$, central over $k$, which is however now required to be totally
  indefinite, and $h$ is a hermitian form on $V$. Then
  $G'(\fR)=Sp(2n,\fR)$, and
$$G(\fR)\cong (Sp(2n,\fR))^f,\ \cD=({\bf III_{\hbf{n}}})^f.$$
\end{itemize}
\item[{\bf U}] Unitary type
\begin{itemize}\item[{\bf U.1}] split case: $G'=SU(V,h)$, where $V$ is an
  $n$-dimensional $K$-vector space, $K|k$ an imaginary quadratic extension,
  and $h$ is a hermitian form. Let for each real prime $\nu$
  $(p_{\nu},q_{\nu})$ denote the signature of $h_{\nu}$. Then
$$G(\fR)\cong \prod_{\nu}SU(p_{\nu},q_{\nu}),\ \ \cD={\bf
  I_{\hbf{p$_{\nu_1}$,q$_{\nu_1}$}}}\times \cdots \times {\bf
  I_{\hbf{p$_{\nu_f}$,q$_{\nu_f}$}}}.$$
\item[{\bf U.2}] non-split case: $G'=SU(V,h)$, where $D$
  is a division algebra of degree $d$, central simple over $K$ ($K$ as in
  ${\bf U.1}$) with a $K|k$-involution and $V$ is an
  $n$-dimensional right $D$-vector space with hermitian form $h$.
  If $d=1$ this reduces to ${\bf
    U.1}$, so we may assume $d\geq 2$. Again
  letting $(p_{\nu},q_{\nu})$ denote the local signatures, we have
$$G(\fR)\cong \prod_{\nu}SU(p_{\nu},q_{\nu}),\ \ \cD\cong {\bf
  I_{\hbf{p$_{\nu_1}$,q$_{\nu_1}$}}}\times \cdots \times {\bf
  I_{\hbf{p$_{\nu_f}$,q$_{\nu_f}$}}}.$$
\end{itemize}
\end{itemize}

\subsubsection{Exceptional groups}
The exceptional groups can be classified by results of Ferrar as we now
describe. A general reference to non-associative algebra used here is
\cite{shaffer}. See also \cite{faulk} for an excellent survey and further
references.
\begin{definition}\label{d46b.1} For an alternative algebra with involution
$(\AA,-)$ let $\gg=(\gg_1,\ldots,\gg_n)$ be
a diagonal matrix with coefficients in $k$, and set
$ \LL(\AA^n,\gg):=\{g\in M_n(\AA) \Big| \gg g^*
\gg^{-1}=g\}.$
One defines the Jordan algebra $J(\AA,\gg)$ by taking $n=3$,
$$J(\AA,\gg):=\LL(\AA^3,\gg).$$
If $\AA$ is an octonion algebra we call
$J(\AA,\gg)$ the {\em exceptional simple Jordan algebra} defined by $\AA$ and
$\gg$.
\end{definition}
In particular, the following cases for exceptional simple Jordan algebras
can occur over $\fR$:
\begin{equation}\label{e46c.0}\begin{minipage}{12cm}
\begin{itemize}\item[(i)] $J^c=J(\CC,(1,1,1))$
(the compact form)
\item[(ii)] $J^b=\eta J(\CC,(1,-1,1))\eta^{-1},\ \eta=diag(1,i,1)$,
\item[(iii)] $J^s=J(\fO,(1,1,1))$ (the split form).
\end{itemize}
\end{minipage}
\end{equation}
There is only one $\fR$-form for the split octonion algebra $\fO$.
Furthermore, for an algebraic number
field $k$, there are $3^t$ isomorphism classes of
Jordan algebras, where $t$ denotes the number of real primes of $k$.
The Jordan algebras $J^c, J^b,$ and $J^s$ have the following explicit matrix
realisations (see \cite{druck}, p. 33)
\begin{equation}\label{e46c.1} J^c \hbox{ (respectively $J^s$) } \isom \left\{
g=
\left(\begin{array}{ccc} \xi_1 & x_3 & \-x_2 \\ \-x_3 & \xi_2 & x_1 \\ x_2 &
\-x_1 & \xi_3 \end{array}\right)\Big| \xi_i\in \fR,\ x_i\in \CC \hbox{
(respectively $x_i\in \fO$)} \right\}
\end{equation}
\begin{equation}\label{e46c.2} J^b\isom \left\{g=\left(\begin{array}{ccc}
\xi_1 &
ix_3 & \-x_2 \\ i\-x_3 & \xi_2 & ix_1 \\ x_2 & i\-x_1 & \xi_3
\end{array}\right) \Big| \xi_i\in \fR, x_i \in \CC \right\},
\end{equation}
and the algebra of Definition \ref{d46b.1} is given explicitly as a
matrix algebra as follows:
\begin{equation}\label{e46c.3} J(\AA,(\gg_1,\gg_2,\gg_3))=\left\{x=\left(
\begin{array}{ccc} \xi_1 & \gg_2x_3 & \gg_3\-x_2 \\ \gg_1\-x_3 & \xi_2 &
\gg_3x_1 \\ \gg_1x_2 & \gg_2\-x_1 & \xi_3 \end{array}\right) \Big| \xi_i\in
k,\ x_i\in \AA\right\}.
\end{equation}
Utilizing composition algebras and Jordan algebras one can construct, with
the following construction of {\it Tits algebras}, exceptional Lie
algebras.
\begin{definition}\label{d46f.1} Let $\AA$ be a composition algebra over $k$,
and $\JJ=J(\BB,\gg)$ a Jordan algebra over another composition algebra $\BB$
as in Definition \ref{d46b.1}. Set:
$$\scL(\AA,\JJ)=Der(\AA)\oplus (\AA_0\otimes\JJ_0)\oplus Der(\JJ).$$
One defines a multiplication $[\cdot,\cdot]$ on $\scL(\AA,\JJ)$, which
extends the $[,]$ products on $Der(\AA)$ and $Der(\JJ)$, by the rules:
\begin{itemize}\item[(a)] $[\cdot,\cdot]$ is bilinear and $[x,x]=0$ for all
$x\in \scL(\AA,\JJ)$;
\item[(b)] $[\cdot,\cdot]$ restricts to the usual commutator on $Der(\AA)$
and $Der(\JJ)$, and these are orthogonal with respect to $[\cdot,\cdot]$,
i.e., $[D,E]=0$ for all $D\in Der(\AA),\ E\in Der(\JJ)$;
\item[(c)] For $D\in Der(\AA), E\in Der(\JJ), D+E$ acts on $\AA_0\otimes
\JJ_0$ by:
$$[D+E,a\otimes x]=D(a)\otimes x + a\otimes E(x);$$
\item[(d)] $[\cdot,\cdot]$ is defined on $\AA_0\otimes \JJ_0$ by the formula:
$$[a\otimes x,b\otimes y]={1\over 3}T(x\circ y)<a,b>+(a*b)\otimes (x*y) +
{1\over 2} T(a\cdot b)<x,y>,$$
where the $*$ and $<,>$ products are defined as in \cite{druck} (in the
cases which we require the definition simplifies somewhat and will be
described below).
\end{itemize}
This makes $\scL(\AA,\JJ)$ a Lie algebra.
\end{definition}
For later use we mention that (\ref{e46c.3}) allows us to write elements in
$J(\AA,\gg)$ in the following way:
\begin{equation}\label{e46d.1} x=\sum_{i=1}^3 \xi_ie_{ii} +
\sum_{i=1}^3x_i[j,k],\ x_i[j,k]:=\gg_kxe_{jk} + \gg_j\-xe_{kj},
\end{equation}
and the second sum is over cyclic permutations $(i,j,k)$ of $(1,2,3)$. In
these terms the norm and trace forms are given by (see \cite{faulk}, 4.11)
\begin{equation}\label{e46d.2} N(x)=\xi_1\xi_2\xi_3 +
\xi_1\gg_2\gg_3n(x_1)+\gg_1\xi_2\gg_3n(x_2) + \gg_1\gg_2\xi_3n(x_3),
\end{equation}
\begin{equation}\label{e46d.3} T(x)=\xi_1+\xi_2+\xi_3.
\end{equation}
In the first formula $n(a)=a\cdot\-a$ is the norm in $\AA$. The norm above is
of course analogous to the determinant in a usual matrix algebra. In
particular, $N(x)\neq 0$ is a neccessary and sufficient condition for $x$ to
be {\em invertible} in $\JJ$, i.e., $N(x)\neq0 \iff \exists_{y\in \JJ}$ with
$x\cdot y=1, x^2\cdot y=x$, and the inverse of $x$ is given by:
\begin{equation}\label{e46d.4} x^{-1}={x^{\#}\over N(x)},
\end{equation}
where $x^{\#}$ satisfies $x\cdot x^{\#}=N(x)\cdot 1$, or explicitly
\begin{equation}\label{e46d.5} x^{\#}=\sum(\xi_j\xi_k-\gg_j\gg_kn(x_i))e_{ii}
+\sum(\gg_i(\-{x_jx_k})-\xi_ix_i)[j,k].
\end{equation}

\paragraph{$\bf E_6$}
There are two constructions leading to
the real Lie algebra of hermitian type $\ee_{6(-14)}$. On the one hand
there is the algebra
$\scL(\fC,J^b)$ (see Definition \ref{d46f.1}), where $J^b$
is isomorphic to the Jordan algebra $\JJ(\CC,(1,-1,1))$ of Definition
\ref{d46b.1} and is given explicitly as a matrix algebra in (\ref{e46c.2}).
Note that in this case the general definition of the algebra
$\scL(\fC,J^b)$ simplifies to
\begin{equation}\label{e54.1} \scL(\fC,J^b)\isom i\cdot J^b_0\oplus Der(J^b),
\end{equation}
which, identifying $J^b_0$ with the right translations by traceless
algebra elements $\cR_{J^b_0}$, is nothing but Albert's twisted
$\scL(\JJ)_{\gl}=\sqrt{\gl}\cR_{\JJ_0}\oplus
Der(\JJ)$, $\scL(\JJ)
= \cR_{\JJ_0}\oplus Der(\JJ)$, as mentioned in \cite{ferrar1}, p.~62.
In our case $\gl=-1$ and $\JJ=J^b$, and this implies the $\bf \ee_6$-form is
of {\em outer} type (see \cite{ferrar1}, \S4 and Theorem 5 b), p.~70). The
Lie multiplication with respect to the decomposition in (\ref{e54.1}) is
given as follows. Writing an element of $\scL(\fC,J^b)$ as $x=i\otimes
A+D,\ A\in J^b_0$ and $D\in Der(J^b)$ and identifying $i\otimes J^b_0$ and
$J^b_0$ so that $x=A+D$, the Lie multiplication is given by
\begin{equation}\label{e54a.2} [A+D,A'+D']= \left(D(A')-D'(A)\right)+
\left([D,D']-[L(A),L(A')]\right),
\end{equation}
where $L(A)$ is left multiplication in $J^b$ by $A$ (cf.~\cite{druck}, 3.2,
p.~46). It turns out that this construction is insufficient to describe all
$k$-forms for number fields $k$.

The other description of $\ee_{6(-14)}$ is as $\scL(\CC,
J^b_1)$, where $J^b_1$ is isomorphic to the Jordan algebra
$J(\AA^c_1,(1,-1,1))=J(\fC,(1,-1,1))$ which can be explicitly described in
matrix terms as
\begin{eqnarray}\label{e54.3} J^b_1 & \cong &
             \hbox{{\script H}}_3(\fC,(1,-1,1)) \\
   & = & \left\{ \left(
\begin{array}{ccc} r_1 & \ga_3 & \-{\ga}_2 \\
-\-{\ga}_3 & r_2 & \ga_1 \\
\ga_2 & -\-{\ga}_1 & r_3 \end{array}\right) \Big| r_i\in \fR, \ga_i\in \fC
\right\}. \nonumber
\end{eqnarray}
It is then clear that $J_1^b\isom \DD^+$ for an associative algebra $\DD$
whose traceless elements form a Lie algebra of type $\Ss\uu(2,1)$. With this
information we can exhibit an explicit isomorphism:
$$\scL(\fC,J^b) \stackrel{\sim}{\lra}\scL(\CC, J^b_1).$$
By means of the isomorphism (\ref{e54.1}) we may represent an element as a
$k$-linear transformation of $J^b$, i.e., as an element of $\CC\otimes
M_3(k)$. Write an element in $\scL(\CC,J^b_1)$ as follows: $D+c\otimes a +
\ad y$, where $D\in Der(\CC), c\in \CC_0, a\in (J^b_1)_0, y\in M_3(k)$. The
isomorphism is given by (\cite{ferrar2}, 2.1)
\begin{eqnarray}\label{e54a.1}
\psi:\hbox{{\script L}}(\CC, J^b_1) & \lra & \hbox{{\script L}}(\fC,J^b) \\
D+c\otimes a +\ad y & \mapsto & D\otimes 1 + (c\otimes a)_r +
(\-c\otimes {^ta})_l + I \otimes (y_r+{^ty}_l). \nonumber
\end{eqnarray}

Now we have the following result of Ferrar concerning $k$-forms of $\ee_6$:
\begin{theorem}[\cite{ferrar2}, p.~201]\label{t54.1} If $L$ is a Lie algebra
of type $\bf \ee_{\hbf{6}}$ over an algebraic number field $k$, then
$$L\isom \scL(\AA_k,J(\BB,\gg))$$
as in Definition \ref{d46f.1} for some octonion algebra $\AA$ and Jordan
algebra $J(\BB,\gg)$ as in Definition \ref{d46b.1}, with $\BB$ an
alternative $k$-algebra of dimension two,
and $\gg=diag(\gg_1,\gg_2,\gg_3)$ is a diagonal $k$-matrix.
\end{theorem}
For our situation of $k$-forms of the $\fR$-algebra $\ee_{6(-14)}$ this
means:
\begin{corollary}\label{c55.1} Any $k$-form of $\ee_{6(-14)}$ (with $k$
totally real) is of the form
$$\scL(\CC_k,(J^b_1)_k),$$
where $\CC_k$ is an anisotropic octonion algebra over $k$ and $(J^b_1)_k$ is
a $k$-form of the algebra (\ref{e54.3}).
\end{corollary}
As a corollary of this we get a classification of $k$-groups of hermitian
$\bf E_6$ type:
\begin{corollary}\label{c55.2} Let $G'$ be an absolutely almost simple
$k$-group of hermitian type, in the class of structures of type $\bf E_6$.
Then $(G')^0\sim \hbox{\em Aut}(\scL(\CC_k,(J^b_1)_k))^0$
with the notations of the
preceeding corollary, where ``$\sim$'' means isogenous.
\end{corollary}
Since an octonion algebra $\AA$ over $k$ is uniquely determined up to
isomorphism by the
set of real primes at which it ramifies, the totally definite (Cayley)
algebra $\CC_k$ is unique, and
we need only apply the classification of $k$-forms of the Lie algebra
$\Ss\uu(2,1)$ to get a complete classification of $k$-forms of $J^b_1$, and
hence a classification of the $k$-forms of $\ee_{6(-14)}$.
There are essentially three cases which can
occur (let $\DD$ denote the associative algebra with involution
and $\DD^-$ the $k$-form of the Lie algebra $\uu(2,1)$):
\begin{equation}\label{e55.1}
\begin{minipage}{12cm}
\begin{itemize}\item[(i)] $(V,h)$ is a $k$-vector space with
hermitian form $h$ of Witt index 1, represented by a matrix $H$,
and $\DD^-=\{g\in End(V)\big|gH-Hg=0\}$.
\item[(ii)] $(V,h)$ is a $k$-vector space with {\em anisotropic} hermitian
form $h$, represented by a matrix $H$, and
$\DD^-=\{g\in End(V)\big|gH-Hg=0\}$.
\item[(iii)] $D$ is a central simple division algebra of degree three over
  an imaginary quadratic extension $K$ of $k$ with a $K|k$-involution,
  and $\DD=D$.
\end{itemize}
\end{minipage}
\end{equation}
Considering the Tits index of these $\fQ$-groups,
note that since a $\fQ$-split torus is all the more $\fR$-split, it follows
that the set of split roots of the index of $G$ (usually drawn white in the
Tits index) are a subset of the split roots of $G(\fR)$. This gives a simple
criterion for deciding which indices may give rise to the given $\fR$-form.
Looking now at the list of $\bf E_6$ indices (of outer type)
in \cite{tits}, the following four possibilities arise for
$\fQ$-forms of $E_{6(-14)}$: $^2E^{78}_{6,0},\ ^2E^{35}_{6,1},\
^2E^{29}_{6,1},\ ^2E^{16'}_{6,2}$. However, as shown in \cite{kazdan}, the
index $^2E^{29}_{6,1}$ does not give rise to a bounded symmetric domain, but
rather has symmetric space $\bf E IV$ in the notation of \cite{Helg}. The
argument is roughly as follows. If $H\inn G$ is the anisotropic kernel, of
type $\bf D_4$, then, since dim$[U,U]=8$ for a maximal unipotent subgroup (in
the maximal $\fQ$-parabolic $P_{\ga_1}\cap P_{\ga_6}$), it follows that
$H\inn End_{\fQ}([U,U])$, a relation preserved upon tensoring with $\fR$, so
that $P_{\ga_2}$ is still not defined over $\fR$; thus the index of $G(\fR)$
is $^1E^{28}_{6,2}$, giving rise to the symmetric space denoted $\bf E IV$
in \cite{Helg}.
Hence there are only three possible Tits indices, namely
$^2E^{16'}_{6,2}$, $^2E^{35}_{6,1}$ and $^2E^{78}_{6,0}$ for $k$-forms of
$\ee_{6(-14)}$, and it may hold that the three possibilities in
(\ref{e55.1}) coincide with the three possible indices.

\paragraph{$\bf E_7$} There are two
constructions utilizing the Tits algebra
leading to the real form of type $\ee_{7(-25)}$. On the one
hand there is the algebra
$\scL(\AA,J^b)\isom \scL(\AA,J^c)$ (see Definition
\ref{d46f.1}), where $\AA\isom M_2(\fR)$ and $J^b$ (respectively
$J^c$) is isomorphic to the
Jordan algebra $J(\CC,(1,-1,1))$ (respectively is the Jordan algebra
 $J(\CC,(1,1,1))$) in
Definition \ref{d46b.1} and is given explicitly as a matrix algebra in
(\ref{e46c.2}) (respectively in (\ref{e46c.1})). In this case the direct sum
decomposition analogous to (\ref{e54.1}) is (\cite{druck}, 4.6, p. 50)
\begin{equation}\label{e55a.1} \scL(\AA, J^b)\isom (\AA_0\otimes J^b)\oplus
Der(J^b).\end{equation}
The multiplication is given by the rules
\begin{equation}
\begin{minipage}{14cm}
\begin{itemize}\item[(i)] $[a\otimes A, b\otimes B]={1\over 2}[a,b]\otimes
A\circ B + {1\over 2}Tr(ab)[L(A),L(B)],$ for $a,b\in \AA_0,\ A,B\in J^b$;
\item[(ii)] $[D,b\otimes B]=b\otimes D(B),\ D\in Der(J^b), b\in \AA_0, B\in
J^b$;
\item[(iii)] $[D,D']=$ usual commutator of $D,D'\in Der(J^b)$.
\end{itemize}\end{minipage}\end{equation}
The other description of $\ee_{7(-25)}$ is as the algebra
$\scL(\CC,\cJ\cO_6(\fR))$, where $\cJ\cO_6(\fR)$ is the Jordan algebra
$\scH_3(M_2(\fR),(1,1,1))$, and is given as a
matrix algebra by (\ref{e46c.3}). Of course we could derive an explicit
isomorphism as in (\ref{e54a.1}) between the two. But in this case it turns
out that the first description is sufficient to get all $k$-forms. Namely, we
have the following result of Ferrar:
\begin{theorem}[\cite{ferrar3}, Theorem 4.3]\label{t55a.1}
Let $k$ be an algebraic
number field $k$ and let $L$ be a $k$-form of the Lie algebra $\ee_7$. Then
$$L\isom \scL(\AA,\JJ)$$
as in Definition \ref{d46f.1} for some quaternion algebra $\AA$ over $k$, and
exceptional simple Jordan algebra $\JJ$ over $k$.
\end{theorem}
For the case at hand here, namely $k$-forms of the $\fR$-algebra
$\ee_{7(-25)}$, this implies
\begin{corollary}\label{t55a.2} Let $G'$ be an almost absolutely simple
$k$-group of hermitian type, of type $\bf E_7$, and let $\Gg'$ be the Lie
algebra. Then
$$\Gg'\cong\scL(\AA_k,\JJ_k),$$
where $\JJ_k$ is exceptional simple such that for each real prime of $k$,
$(\JJ_k)_{\nu}\isom J^b$ or $J^c$, and $\AA_k$ is a
quaternion algebra over $k$ which splits at all infinite primes $\nu$.
\end{corollary}
There are the following possibilities over $\fQ$:
\begin{equation}\label{e55b.1}
\begin{minipage}{14cm}\begin{itemize}\item[(i)] $\AA$ is
split, $\JJ_{\fQ}$ is a $\fQ$-form of $J^b$;
\item[(ii)] $\AA$ is split, $\JJ_{\fQ}$ is a $\fQ$-form of $J^c$;
\item[(iii)] $\AA$ is division, $\JJ_{\fQ}$ is a $\fQ$-form of $J^b$;
\item[(iv)] $\AA$ is division, $\JJ_{\fQ}$ is a $\fQ$-form of $J^c$.
\end{itemize}\end{minipage}\end{equation}
There are three possible Tits indices, namely $E^{28}_{7,3}, E^{31}_{7,2}$
and $E^{133}_{7,0}$. It is rather clear that the first (respectively the last)
case above gives rise to $E^{28}_{7,3}$ (respectively to $E^{133}_{7,0}$),
and it seems natural to expect the other two cases to give rise to
$E^{31}_{7,2}$.

\subsection{Boundary components}
We briefly discuss the rational boundary components occuring in each of the
cases. Again we tabulate this, giving the Tits index in each case and
describing the boundary components. We also describe, in the classical
cases, the corresponding isotropic subspaces of the vector space
$V$. Throughout, $\cD^*$ denotes the union of $\cD$ and the rational
boundary components.

\begin{itemize}\item[{\bf O.1}] The Tits index is $D_{n,s}$ (for
  $n\equiv2(4)$), ${^2D}_{n,s}$ (for $n\equiv0(4)$) or $B_{n,s}$
  (for $n$ odd), where $s$ is the Witt
  index of $h$. The corresponding diagrams are (the top diagrams are for
  the case $s=2$, the lower ones giving the left ends for $s=1$):

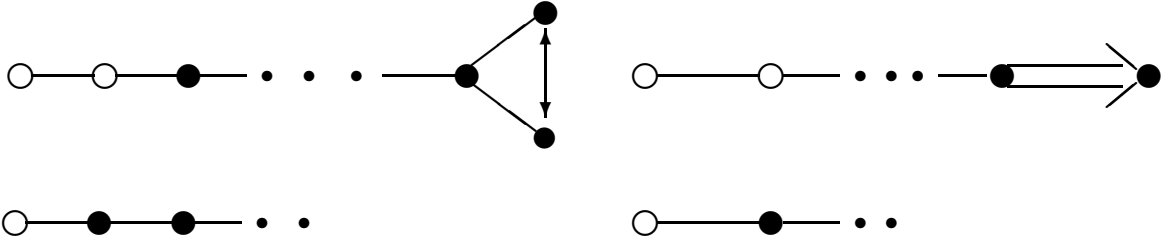
\setlength{\unitlength}{0.005500in}%
\begin{picture}(1102,222)(64,609)
\thicklines
\put(160,760){\circle{22}}
\put(315,760){\circle*{10}}
\put(355,760){\circle*{10}}
\put( 80,760){\circle{22}}
\put( 90,760){\line( 1, 0){ 60}}
\put(170,760){\line( 1, 0){ 60}}
\put(250,760){\line( 1, 0){ 45}}
\put(880,760){\circle*{10}}
\put(910,760){\circle*{10}}
\put(675,760){\circle{22}}
\put(795,760){\circle{22}}
\put(807,760){\line( 1, 0){ 53}}
\put(687,760){\line( 1, 0){ 97}}
\put(400,760){\circle*{10}}
\put(580,820){\circle*{22}}
\put(579,701){\circle*{20}}
\put(235,620){\circle*{22}}
\put(310,620){\circle*{10}}
\put(350,620){\circle*{10}}
\put( 75,620){\circle{22}}
\put(240,760){\circle*{22}}
\put(155,620){\circle*{22}}
\put(687,620){\line( 1, 0){ 97}}
\put(505,760){\circle*{22}}
\put(935,760){\circle*{10}}
\put(1015,760){\circle*{22}}
\put(1155,760){\circle*{22}}
\put(880,620){\circle*{10}}
\put(910,620){\circle*{10}}
\put(675,620){\circle{22}}
\put(795,620){\circle*{22}}
\put(425,760){\line( 1, 0){ 70}}
\put(510,770){\line( 4, 3){ 60}}
\put(511,752){\line( 4,-3){ 60}}
\put( 85,620){\line( 1, 0){ 60}}
\put(165,620){\line( 1, 0){ 60}}
\put(245,620){\line( 1, 0){ 45}}
\put(955,760){\line( 1, 0){ 45}}
\put(1020,770){\line( 1, 0){110}}
\put(1020,750){\line( 1, 0){110}}
\put(1115,730){\line( 6, 5){ 28}}
\put(1115,790){\line( 6,-5){ 28}}
\put(807,620){\line( 1, 0){ 53}}
\put(580,805){\vector( 0,-1){ 85}} \put(580,805){\vector( 0,1){ 0}}

\end{picture}

\noindent where the Galois action in the left-hand diagram is present only for
$n\equiv 0(4)$. The boundary components of $\cD'=G'(\fR)/K'$ are:
\begin{itemize}\item $\{pt\}\inn \{ \hbox{1-disc}\}^*$,  ($s=2$)
\item $\{pt\}$, ($s=1$).
\end{itemize}
\item[{\bf O.2}] The index in this case is $D_{{n\over 2},s}^{(2)}$ ($n$
  even) or ${^2D}_{{n-1\over2},s}^{(2)}$ ($n$ odd), where $s$ is the Witt
  index of $h$. The corresponding diagrams are (with non-trivial Galois
  action identifying the two right most vertices for $n$ odd):

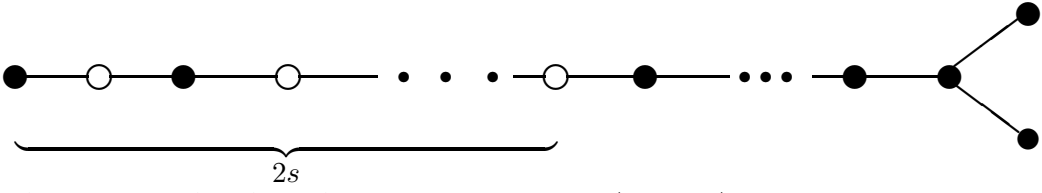
\setlength{\unitlength}{0.005500in}%
\begin{picture}(987,141)(14,670)
\thicklines
\put(435,740){\circle*{10}}
\put(480,740){\circle*{10}}
\put( 25,740){\circle*{22}}
\put(185,740){\circle*{22}}
\put(105,740){\circle{22}}
\put(285,740){\circle{22}}
\put(540,740){\circle{22}}
\put(625,740){\circle*{22}}
\put(740,740){\circle*{10}}
\put(760,740){\circle*{10}}
\put(720,740){\circle*{10}}
\put(825,740){\circle*{22}}
\put(395,740){\circle*{10}}
\put(915,740){\circle*{22}}
\put(921,732){\line( 4,-3){ 60}}
\put(990,680){\circle*{20}}
\put(990,800){\circle*{22}}
\put( 35,740){\line( 1, 0){ 60}}
\put(115,740){\line( 1, 0){ 60}}
\put(195,740){\line( 1, 0){ 80}}
\put(295,740){\line( 1, 0){ 75}}
\put(500,740){\line( 1, 0){ 30}}
\put(550,740){\line( 1, 0){ 65}}
\put(635,740){\line( 1, 0){ 70}}
\put(705,740){\line(-1, 0){  5}}
\put(785,740){\line( 1, 0){ 40}}
\put(835,740){\line( 1, 0){ 70}}
\put(920,750){\line( 4, 3){ 60}}
\put(25,680){$\underbrace{\hspace*{7.2cm}}_{\hbox{$2s$}}$}

\end{picture}

\vspace*{.2cm}
\noindent The corresponding boundary components are $\bf II_{\hbf{n-2}}^*\nni
II_{\hbf{n-4}}^*\nni \cdots \nni II_{\hbf{n-2s}}$.

\item[{\bf S.1}] The index is $C_{n,n}$, with the usual diagram and the
  following boundary components: $\{pt\}\inn \bf III_{\hbf{1}}^*\inn \cdots
  \inn III_{\hbf{n-1}}^*.$
\item[{\bf S.2}] The index is $C_{n,s}^{(2)}$, with diagram
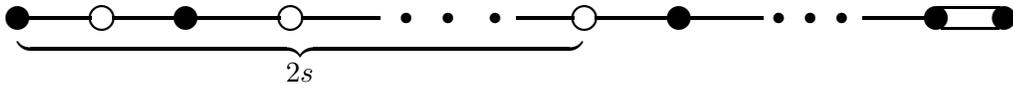
\setlength{\unitlength}{0.005500in}%
$$\begin{picture}(962,22)(14,729) \thicklines \put(435,740){\circle*{10}}
  \put(480,740){\circle*{10}} \put( 25,740){\circle*{22}} \put(
  25,720){$\underbrace{\hspace*{7.5cm}}_{\hbox{$2s$}}$}
  \put(185,740){\circle*{22}} \put(105,740){\circle{22}}
  \put(285,740){\circle{22}} \put(565,740){\circle{22}}
  \put(655,740){\circle*{22}} \put(750,740){\circle*{10}}
  \put(780,740){\circle*{10}} \put(810,740){\circle*{10}}
  \put(395,740){\circle*{10}} \put(900,740){\circle*{22}}
  \put(901,750){\line( 1, 0){ 60}} \put(965,740){\circle*{22}} \put(
  35,740){\line( 1, 0){ 60}} \put(115,740){\line( 1, 0){ 60}}
  \put(195,740){\line( 1, 0){ 80}} \put(295,740){\line( 1, 0){ 75}}
  \put(500,740){\line( 1, 0){ 55}} \put(575,740){\line( 1, 0){ 70}}
  \put(665,740){\line( 1, 0){ 70}} \put(735,740){\line(-1, 0){ 5}}
  \put(830,740){\line( 1, 0){ 60}} \put(905,730){\line( 1, 0){ 55}}
\end{picture}$$

\vspace*{.2cm}\noindent
The boundary components are then the following: $\bf III_{\hbf{n-2}}^*\nni
\cdots \nni III_{\hbf{n-2s}}$.
\item[{\bf U.1}] The index is ${^2A}_{n-1,s}$, with the diagram

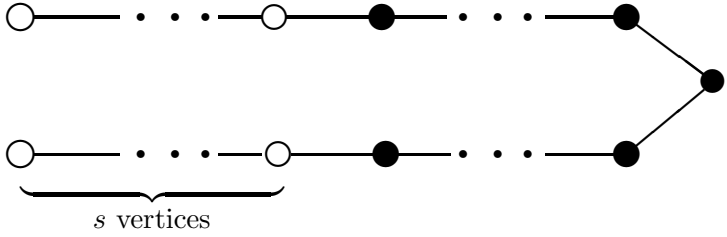
\setlength{\unitlength}{0.004500in}%
\begin{picture}(734,190)(-30,610)
\thicklines
\put(260,785){\circle*{10}}
\put(295,785){\circle*{10}}
\put(220,625){\circle*{10}}
\put(260,625){\circle*{10}}
\put(295,625){\circle*{10}}
\put(595,785){\circle*{10}}
\put(635,785){\circle*{10}}
\put(670,785){\circle*{10}}
\put(595,625){\circle*{10}}
\put(635,625){\circle*{10}}
\put(670,625){\circle*{10}}
\put(885,710){\circle*{28}}
\put(785,785){\circle*{30}}
\put(785,625){\circle*{30}}
\put(785,785){\line( 4,-3){100}}
\put(220,785){\circle*{10}}
\put(785,625){\line( 6, 5){ 90}}
\put( 95,625){\line( 1, 0){100}}
\put(690,785){\line( 1, 0){100}}
\put(690,625){\line( 1, 0){100}}
\put(375,785){\circle{28}}
\put(500,785){\circle*{30}}
\put(310,785){\line( 1, 0){ 50}}
\put(390,785){\line( 1, 0){ 95}}
\put(515,785){\line( 1, 0){ 60}}
\put(380,625){\circle{28}}
\put(505,625){\circle*{30}}
\put(310,625){\line( 1, 0){ 50}}
\put(395,625){\line( 1, 0){ 95}}
\put(520,625){\line( 1, 0){ 60}}
\put( 80,785){\circle{30}}
\put( 80,625){\circle{30}}
\put( 95,785){\line( 1, 0){100}}
\put(80,590){$\underbrace{\hspace*{3.5cm}}_{{\displaystyle s}
\hbox{ vertices}}$}
\end{picture}

\vspace*{1cm}\noindent
As above, let $(p_{\nu},q_{\nu})$ denote the signature of $h_{\nu}$, then
in the factor $\cD_{\nu}$ of $\cD$ we have the boundary components of the
type $\bf I_{\hbf{p$_{\nu}$-b,q$_{\nu}$-b}}$ for $1\leq b\leq s$. Hence a flag
of boundary components will be
\[ \prod \bf I_{\hbf{p$_{\nu}$-1,q$_{\nu}$-1}}^*\nni \prod
I_{\hbf{p$_{\nu}$-2,q$_{\nu}$-2}}^* \nni \cdots \nni \prod
I_{\hbf{p$_{\nu}$-s,q$_{\nu}$-s}}.\]
\item[{\bf U.2}] The index is in this case ${^2A}_{nd-1,s}^{(d)}$, with
  diagram
\vspace*{.5cm}
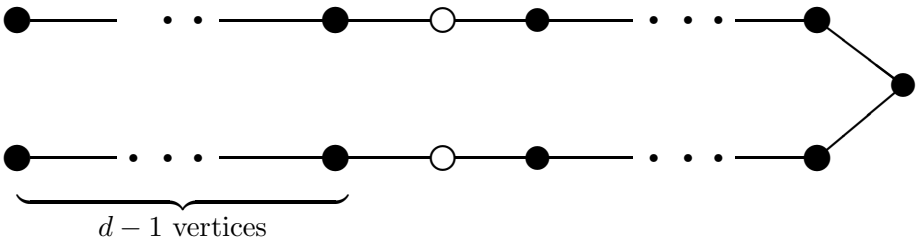
$$\setlength{\unitlength}{0.004500in}%
\begin{picture}(600,560)(400,235)
  \thicklines \put(280,780){\circle*{10}} \put(315,780){\circle*{10}}
  \put(105,780){\circle*{30}} \put(475,780){\circle*{30}}
  \put(120,780){\line( 1, 0){100}} \put(340,780){\line( 1, 0){120}}
  \put(845,780){\circle*{10}} \put(885,780){\circle*{10}}
  \put(920,780){\circle*{10}} \put(600,780){\circle{28}}
  \put(710,780){\circle*{28}} \put(1035,780){\circle*{30}}
  \put(490,780){\line( 1, 0){ 95}} \put(615,780){\line( 1, 0){ 80}}
  \put(720,780){\line( 1, 0){100}} \put(940,780){\line( 1, 0){100}}
  \put(240,620){\circle*{10}} \put(280,620){\circle*{10}}
  \put(315,620){\circle*{10}} \put(105,620){\circle*{30}}
  \put(475,620){\circle*{30}} \put(120,620){\line( 1, 0){100}}
  \put(340,620){\line( 1, 0){120}} \put(845,620){\circle*{10}}
  \put(885,620){\circle*{10}} \put(920,620){\circle*{10}}
  \put(600,620){\circle{28}} \put(710,620){\circle*{28}}
  \put(1035,620){\circle*{30}} \put(490,620){\line( 1, 0){ 95}}
  \put(615,620){\line( 1, 0){ 80}} \put(720,620){\line( 1, 0){100}}
  \put(940,620){\line( 1, 0){100}} \put(1135,705){\circle*{28}}
  \put(1035,780){\line( 4,-3){100}} \put(1035,620){\line( 6, 5){ 90}}

  \put(105,580){$\underbrace{\hspace*{4.4cm}}_{ \hbox{$d-1$ vertices}}$}
\end{picture}$$

\vspace*{-3.5cm}\noindent where there are $2s$ white vertices
altogether. Letting the notations be as for the case $\bf U.1$, we have the
following boundary components:
\[ \prod \bf I_{\hbf{p$_{\nu}$-d,q$_{\nu}$-d}}^*\nni \prod
I_{\hbf{p$_{\nu}$-2d,q$_{\nu}$-2d}}^* \nni \cdots \nni \prod
I_{\hbf{p$_{\nu}$-sd,q$_{\nu}$-sd}}.\]
\end{itemize}

We now describe briefly the parabolics in terms of the geometry of $(V,h)$
for all the cases above. Fixing a maximal $k$-split torus and an order on
it amounts to fixing a maximal totally isotropic ($s$-dimensional) subspace
$H_1\inn V$ and a basis $v_1,\ldots, v_s$ of $H_1$. There are then
$k$-vectors $v_1',\ldots, v_s'$ spanning a complementary totally isotropic
subspace $H_2$ such that $h(v_i,v_j')=\gd_{ij}$. Then each pair
$(v_i,v_i')$ spans a hyperbolic plane $V_i$ (over $D$), and $V$ decomposes:
\begin{equation}\label{E8.a} V=V_1\oplus \cdots \oplus V_s\oplus V',\quad
  V'\hbox{ anisotropic for $h$}.
\end{equation}
Furthermore, $V_1\oplus \cdots \oplus V_s=H_1\oplus H_2$. With these
notations, for $1\leq b\leq s$ the standard $k$-parabolic $P_b'\inn G'$
is given as follows:
\begin{equation}\label{E8.b} P_b'=\cN_{G'}(<v_1,\ldots, v_b>),
\end{equation}
where $<v_1,\ldots,v_b>$ denotes the span, a $b$-dimensional totally
isotropic subspace. The hermitian Levi factor of $P_b'$ is
\begin{equation}\label{E8.c}
L_b'=\cN_{G'}(V_{b+1}\oplus \cdots \oplus V_s\oplus V')/\cZ_{G'}(V_{b+1}\oplus
\cdots \oplus V_s\oplus V').
\end{equation}
It reduces to the $k$-anisotropic kernel for $b=s$.

For the exceptional cases we have the following possibilities:
\begin{itemize}\item {\bf $ E_{\hbf{6}}$:} Index: ${^2E}_{6,2}^{16'}$,
  boundary components: $\{pt\}\inn \fB_5^*$.

Index: ${^2E}_{6,1}^{35}$, boundary components: $\fB_5$.
\item {\bf $E_{\hbf{7}}$:} Index: $E_{7,3}^{28}$, boundary components
  $\{pt\}\inn \bf IV_{\hbf{1}}^*\inn IV_{\hbf{10}}^*$.

Index: $E_{7,2}^{31}$, boundary components $\bf IV_{\hbf{1}}\inn
IV_{\hbf{10}}^*$.
\end{itemize}

\section{Rational symmetric subgroups and incidence}
\subsection{Holomorphic symmetric embeddings}\label{section2.1}
Recall that an injection $i_{\cD}:\cD\hra \cD'$ of symmetric spaces is said
to be {\it strongly equivariant} if $i_{\cD}$ is induced by an injection
$i:\Gg\hra \Gg'$ of the Lie algebras $\Gg$ (resp. $\Gg'$) of the real Lie
group $G=\Aut(\cD)$ (resp. $G'=\Aut(\cD')$).  This is equivalent to the
condition that $i_{\cD}(\cD)$ is totally geodesic in $\cD'$ with respect to
the $G'$-invariant metric on $\cD'$. Assuming both $\cD$ and $\cD'$ are
hermitian symmetric, there exist elements $\xi$ (resp. $\xi'$) in the
center of the maximal compact subgroup $K$ (resp. $K'$) such that
$J=\ad(\xi)$ (resp. $J'=\ad(\xi')$) gives the complex structure, and the
condition that $i_{\cD}$ be holomorphic is

\vspace*{.2cm}
$\hbox{(H$_1$)}\hspace*{5.8cm} i\circ \ad(\xi) = \ad(\xi')\circ i.$

\vspace*{.2cm}
For any given hermitian symmetric space $\cD'$, the possible
hermitian symmetric subdomains $i_{\cD}(\cD)$ have been classified by Satake
and Ihara (see \cite{I} and \cite{S1}). Note in particular that the above
applies to $\cD_N$, where $N$ is a reductive subgroup of hermitian type and
$\cD_N$ is the associated hermitian symmetric space. We will refer to
subgroups $N\inn G$, where $G$ is the connected component of the
automorphism group of $\cD$, for which $\cD_N\inn \cD$ is a hermitian
symmetric subdomain, as {\it symmetric subgroups} $N\inn G$\footnote{The
  term ``symmetric'' arises from the fact that in most cases, $N$ can be
  defined in terms of closed symmetric sets of roots.}. For this notion it
is irrelevant whether $N$ is reductive, semisimple or even centerless.

\subsection{Incidence over $\fR$}
In this section let $G$ be a reductive Lie group of hermitian type such
that the symmetric space $\cD$ is irreducible, and let $A\inn G$ be the
maximal $\fR$-split torus (with order) defined by the maximal set of
strongly orthogonal roots of $G$ as in 1.1.1. Then we can speak
of the standard parabolics $P_b,\ b=1,\ldots, t$, $t=\rank_{\fR}G$.
We introduce the set of domains $(\cE\cD)$ as follows.

$(\cE\cD)\hspace*{5cm} \bf I_{\hbox{\scsi\bf q,q}},\ II_{\hbox{\scsi\bf
    n}},\ n \hbox{ even},\ III_{\hbox{\scsi\bf n}}.$

With
respect to a fixed $P_b$ we consider the following conditions on a symmetric
subgroup $N\inn G$ as in section \ref{section2.1}.
\begin{itemize}\item[1)] $N$ has maximal $\fR$-rank, that is,
  $\rank_{\fR}N=\rank_{\fR}G$.
\item[2)] $N$ is a maximal symmetric subgroup.
\item[2')] $N$ is a maximal subgroup of tube type, i.e.,
  such that $\cD_N$ is a tube domain.
\item[2'')] $N$ is {\it minimal}, subject to 1).
\item[3)] $N=N_1\times N_2$, where $N_1\inn P_b$ is a hermitian Levi factor of
  $P_b$ for some Levi decomposition.
\item[3')] $\cD_N^*$ contains $F$ as a boundary component.
\end{itemize}
\begin{definition}\label{d9.1} Let $G$ be a simple real Lie group of
  hermitian type, $A$ a fixed maximal
  $\fR$-split torus defined as above by a maximal set of strongly
  orthogonal roots, $\eta_i,\ i=1,\ldots, t$ the simple $\fR$-roots,
  $F_{b}$ a standard boundary component and
  $P_{b}$ the corresponding standard maximal $\fR$-parabolic, $b\in
  \{1,\ldots, t\}$. A
  symmetric subgroup $N\inn G$ (respectively the subdomain $\cD_{N}\inn
  \cD$) will be called {\it
    incident} to $P_{b}$ (respectively to $F_{b}$), if $N$ fulfills:
  \begin{itemize}\item $b<t$, then $N$ satisfies 1), 2), 3).
   \item $b=t,\ \cD\not\in (\cE\cD)$, then $N$ satisfies 1), 2) or 2'), 3').
   \item $b=t,\ \cD\in (\cE\cD)$, then $N$ satisfies 1), 2''), 3').
  \end{itemize}
  For reducible $\cD=\cD_1\times \cdots \times \cD_d$, we have the product
  subgroups $N_{b_1,1}\times \cdots \times N_{b_d,d}$, where
  $\cD_{N_{b_i,i}}$ is
  incident to the standard boundary component $F_{{b_i}}$ of $\cD_i$
  (and $N_{0,i}=G_i$).
\end{definition}
This defines the notion of symmetric subgroups incident with a
standard parabolic. Any maximal $\fR$-parabolic is conjugate to one and
only one standard maximal parabolic, $P=gP_bg^{-1}$ for some $b$. Let $N_b$
be any symmetric subgroup incident with $P_b$. Then just as one has the
pair $(P_b,N_b)$ one has the pair $(P,N)$,
\begin{equation}\label{e10.3} P=gP_bg^{-1},\quad N=g N_b g^{-1}.
\end{equation}
\begin{definition} \label{d10.1} A pair $(P,N)$ consisting of a maximal
  $\fR$-parabolic $P$ and a symmetric subgroup $N$ is called {\it
    incident}, if the groups are conjugate by a common element $g$ as in
  (\ref{e10.3}) to a pair
  $(P_b,N_b)$ which is incident as in Definition \ref{d9.1}.
\end{definition}
The existence of the symmetric subgroups $N_b$ was proved in the above
mentioned work of Ihara and Satake.
Let $P_{b}$, $1\leq b< t$ (this means $\dim(F_b)>0$)
be a standard parabolic and let $L_b$ be the ``standard'' hermitian Levi
factor, i.e., such that $Lie(L_b)=\ll_b$; then
\begin{equation}\label{e10.1} N_b:= L_b\times \cZ_G(L_b)
\end{equation}
is a subgroup having the properties given in the definition, unique since
$L_b$ is unique. We shall refer to this unique subgroup as the {\it
  standard} incident subgroup.
As to uniqueness, the following was
shown in \cite{sym}, Prop.~2.4.
\begin{proposition} If $(N,P_b)$ are incident, there is $g\in V_b$ such
  that $N$ is conjugate by $g$ to the standard $N_b$ of (\ref{e10.1}),
  where $V_b$ is the factor of $P_b$ of Theorem \ref{t4.1}.
\end{proposition}
The situation for zero-dimensional boundary components was not considered
in \cite{sym} in detail, so we take this up now.

Consider first the case where $\cD\not\in(\cE\cD)$, so $\cD$ is a
product of factors of types:
\[\bf I_{\hbf{p,q}}\ \hbox{ ($p>q$)},
II_{\hbf{n}}\ \hbox{ ($n$ even)},\ IV_{\hbf{n}},\ V\hbox{ or }\ \bf VI.\]
The
corresponding subgroups $N_t$ are: $\bf I_{\hbf{p-1,q}}\inn I_{\hbf{p,q}},\
II_{\hbf{n-1}}\inn II_{\hbf{n}},\ IV_{\hbf{n-1}}\inn IV_{\hbf{n}},\
I_{\hbf{2,4}},\ II_{\hbf{5}}$ or $\bf IV_{\hbf{8}}\inn V,\ I_{\hbf{3,3}}$
or $\bf II_{\hbf{6}}\inn VI$. Next note that if $N$ is incident to $P_t$,
so $\cD_N$ is incident to $F_t$ (=pt.), then any other domain $\cD_{N'}$
incident to $F_t$ (isomorphic to the given $\cD_N$)
will be the conjugate by some element of $G$ fixing
$F_t$, that is by $g\in P_t$. If furthermore $g\in N$, then $g$
leaves $\cD_N$ invariant. It follows that $N$ is unique (in its isomorphism
class for type $\bf V$ and $\bf VI$) up to elements in $P_t$ modulo those
in $N_t$.
Hence we must find the intersection $N_t\cap P_t$. This can be done
in the Lie algebras, i.e., we must find $\nn_t\cap \pp_t$.

Ihara has shown that all the subalgebras $(\nn_t)_{\fC}$ (with
the exception of $\bf IV_{\hbf{n-1}}\inn IV_{\hbf{n}}$, $n$ even) are {\it
  regular} subalgebras, i.e., are generated by the Cartan subalgebra $\tt$
and the root spaces $\Gg^{\ga}$ for $\ga\in \Psi_{sym}$, where $\Psi_{sym}$
is a closed, {\it symmetric} set of roots. Similarly, $(\pp_t)_{\fC}$ is
the subalgebra generated by $\tt$ and the root spaces $\Gg^{\ga}$ for
$\ga\in \Psi_{par}$, where $\Psi_{par}$ is a closed, {\it parabolic} set of
roots, $\Psi_{par}=\Phi^+\cup [\gD-\gt]$, where $\gt\inn \gD$ is some
subset of the set of simple roots, and for any subset $\Xi\inn \gD$,
$[\Xi]$ denotes the set of roots which are integral linear combinations of
the elements of $\Xi$. Then the intersection of $(\pp_t)_{\fC}$
and $(\nn_t)_{\fC}$ is given by
\[(\pp_t)_{\fC}\cap (\nn_t)_{\fC}=\tt+\left(\sum_{\ga\in \Psi_{sym}}
    \Gg^{\ga} \cap \sum_{\ga\in \Psi_{par}}\Gg^{\ga}\right)= \tt +
    \sum_{\ga\in \Psi_{sym}\cap \Psi_{par}}\Gg^{\ga}.\]
{}From this it follows that the complement of $(\pp_t)_{\fC}\cap
(\nn_t)_{\fC}$ in $(\pp_t)_{\fC}$ is given by
\[\cc=\sum_{\ga\in \Psi_{par}-(\Psi_{sym}\cap \Psi_{par})} \Gg^{\ga}.\]
This is of course not a subalgebra, but we can determine the dimension of
the parameter space of non-trivial conjugates of $N_t$ incident with $P_t$.
In other words, the homogenous space $P_t/(P_t\cap N_t)$ can be identified
with the set of symmetric subgroups $N$ incident with $P_t$; its dimension
is the cardinality of the set of roots
$\Psi_{par}-\Psi_{sym}\cap \Psi_{par}$. To
demonstrate this consider $SU(4,1)$. Let $\ga_1=\ge_1-\ge_2,\ldots,
\ga_{4}=\ge_4-\ge_5$ denote
the simple roots for $\Gg_{\fC}$, we have :
\[\Psi_{sym}=\pm(\ge_2-\ge_3),\ \pm(\ge_2-\ge_4),\ \pm(\ge_2-\ge_5),\
\pm(\ge_3-\ge_4),\ \pm(\ge_3-\ge_5),\ \pm(\ge_4-\ge_5),\]
\[\Psi_{par}=+(\ge_i-\ge_j),\ \hbox{ (10 of these) }, \pm(\ge_2-\ge_3),\
\pm(\ge_2-\ge_4),\ \pm(\ge_3-\ge_4),\]
so that $\Psi_{par}-(\Psi_{sym}\cap \Psi_{par})=+(\ge_1-\ge_j),\
j=2,\ldots, 5$. Hence, taking the relation $\sum\ge_i=0$ into account,
 there are three effective parameters. Geometrically this
can be seen as follows. The bounded symmetric domain is a four-dimensional
ball, the boundary component is a point, and the symmetric subdomain
$\cD_N$ is an embedded three-ball passing through the point. Now think of
the four-ball as embedded in $\fP^4$ via the Borel embedding; the three-ball is
the intersection of $\fB_4\inn \fP^4$ with a hyperplane $\fP^3$ passing
through the given point. There is an infinitesimal $\fP^3$ of hyperplanes
through the point, so we see three effective parameters.

Now we turn to the embedding $\bf IV_{\hbf{n-1}}\inn IV_{\hbf{n}}$, $n$
even. If $G=SO(V,h)$, $h$ symmetric of signature $(n,2)$,
let $v\in V$ be an anisotropic vector. Then $v^{\perp}$
is of codimension one, $h_{|v^{\perp}}$ has signature $(n-1,2)$ and
$N_t=N_G(v^{\perp})$. On the other hand the parabolic $P_t$ is the
stabilizer of a (maximal) two-dimensional totally isotropic subspace $I\inn
V$. Then $V$ splits off two hyperbolic planes $H_1,\ H_2$, and $v$ is in
the orthogonal complement of $H:=H_1\oplus H_2$. So the intersection
$N_t\cap P_t$ is just the stabilizer of $v$ in $P_t$, i.e.,
\[ N_t\cap P_t =\{g\in G | g(I)\subseteq I,\ g(v)\in <v>\}.\]

Finally we mention the case $\cD\in (\cE\cD)$. Then $\cD_N$ is a polydisc
and it is easy to see that the intersection $N_t\cap P_t$ is just the
parabolic in $N_t$ corresponding to the given point. Since $N_t\cong
(SL_2)^t,\ t=\rank_{\fR}\cD$, the parabolic is $(P_1)^t$, where $P_1\inn
SL_2$ is the standard one-dimensional parabolic. So the number of
parameters in this case is the dimension of $P_t$ minus $t$.

We now list the sets $\Psi_{sym}$, following Ihara, but we will use the
notations of the root systems as in \cite{bour}.
\begin{itemize}\item ${\bf I_{\hbf{p,q}}}$: $\gD=\{\ga_1,\ldots,
  \ga_{p+q-1}\},\ \Psi_{sym}=[\ga_2,\ldots, \ga_{p+q-1}].$
\item ${\bf II_{\hbf{n}}}$: ($n$ even). $\gD=\{\ga_1,\ldots, \ga_{[{n\over
    2}]},\}
  \Psi_{sym}=[\ga_2,\ldots, \ga_{[{n\over2}]}]$.
\item ${\bf IV_{\hbf{n}}}$: ($n=2\ell+1$), $\gD=\{\ga_1,\ldots,\ga_{\ell}\},\
  \gb:=\ga_{\ell-1}+2\ga_{\ell}$. The following set of roots forms a
  diagram of type $D_{\ell}$ as indicated:
$$
  \setlength{\unitlength}{1cm}\begin{picture}(3,1.5)
\put(0,1){$\ga_1$}
\put(.5,1){$\ga_2$}
\put(1,1){$\cdots$}
\put(1.5,1){$\ga_{\ell-2}$}
\put(2.5,1){$\ga_{\ell}$}
\put(1.5,.5){$\gb$}
\end{picture}$$

\vspace*{-1cm}
\item ${\bf V}$: $\gD=\{\ga_1,\ldots, \ga_6\},\
  \gb_1:=\ga_2+\ga_3+2\ga_4+\ga_5+\ga_6,\ \gb_2:=
  \ga_2+\ga_4+\ga_5+\ga_6$. Then the subalgebras are determined by the
  following sets of roots:
\[  \setlength{\unitlength}{1cm}\begin{picture}(10.5,1.5)
\put(.75,1.4){$\bf I_{\hbf{2,4}}\times SU(2)$}
\put(5.75,1.4){$\bf II_{\hbf{5}}$}
\put(9.75,1.4){$\bf IV_{\hbf{8}}$}
\put(0,1){$\gb_1$}
\put(.5,1){$\ga_1$}
\put(1,1){$\ga_3$}
\put(1.5,1){$\ga_4$}
\put(2,1){$\ga_2$}
\put(2.5,1){$\cup$}
\put(3,1){$\ga_6$}
\put(5,1){$\ga_1$}
\put(5.5,1){$\ga_3$}
\put(6,1){$\ga_4$}
\put(6.5,1){$\ga_5$}
\put(5.5,.5){$\gb_2$}
\put(9,1){$\ga_1$}
\put(9.5,1){$\ga_3$}
\put(10,1){$\ga_4$}
\put(10.5,1){$\ga_5$}
\put(10,.5){$\ga_2$}
\end{picture}\]

\vspace*{-1cm}
\item $\bf VI$: $\gD=\{\ga_1,\ldots, \ga_7\},\
  \gb_1:=\ga_6+2\ga_5+3\ga_4+2\ga_3+\ga_1+2\ga_2,\
  \gb_2:=\ga_5+2\ga_4+2\ga_3+\ga_1+\ga_2$. Then the subalgebras are
  determined by the following sets of roots:
\[  \setlength{\unitlength}{1cm}\begin{picture}(8,1.5)
\put(1,1.4){$\bf I_{\hbf{3,3}}$} \put(7,1.4){$\bf II_{\hbf{6}}$}
\put(-.2,1){$(-\ga_2)$}
\put(1,1){$\gb_1$}
\put(1.5,1){$\ga_7$}
\put(2,1){$\ga_6$}
\put(2.5,1){$\ga_5$}
\put(6,1){$\ga_7$}
\put(6.5,1){$\ga_6$}
\put(7,1){$\ga_5$}
\put(7.5,1){$\ga_4$}
\put(8,1){$\ga_2$}
\put(6.5,.5){$\gb_2$}
\end{picture}
\]
\end{itemize}

\vspace*{-.75cm}
We can also consider the converse question, i.e., given a symmetric
subgroup, what is the set of parabolics to which it is indicent? The answer
to this is easier: if $\dim(F)>0$, then for any other boundary component
$F'$ of
$\cD_N$, conjugate to $F$, the parabolic $P_{F'}=N(F')$ is also incident to
$N$. These boundary components are in 1-1 correspondence with the
zero-dimensional boundary components of the second factors $\cD_2$ of
$\cD_N=\cD_1\times \cD_2$. If $\dim(F)=0$, then, assuming $\cD_N$ is
irreducible (i.e., $\cD\not\in(\cE\cD)$), then for any other
zero-dimensional boundary component $F'$, the corresponding stabilizer
$P_{F'}$ is incident with $N$. If $\cD\in (\cE\cD)$, then we have the set
of zero-dimensional boundary components of the polydisc.

\subsection{Incidence over $\fQ$}

We now return to the notations of section 1.1.2;
$G$ is a simple $\fQ$-group of hermitian
type. The following definition gives a $\fQ$-form of Definition
\ref{d10.1}.
\begin{definition}\label{d12.1} Let $P\inn G$ be a maximal $\fQ$-parabolic,
  $N\inn G$ a reductive $\fQ$-subgroup. Then we shall say that $(P,N)$ are
  {\it incident} (over $\fQ$), if $(P(\fR),N(\fR))$ are incident in the
  sense of Definition \ref{d10.1}.
\end{definition}
Note that in particular $N$ must itself be of hermitian type, and such that
the Cartan involution of $G(\fR)$ restricts to the Cartan involution of
$N(\fR)$.
The main result of \cite{sym} is the following existence result.

\begin{theorem}\label{t12.1} Let $G$ be $\fQ$-simple of hermitian type
  subject to the restrictions above ($G$ is isotropic and $G(\fR)$ is not a
  product of $SL_2(\fR)$'s),
  $P\inn G$ a $\fQ$-parabolic. Then there exists a reductive $\fQ$-subgroup
  $N\inn G$ such that $(P,N)$ are incident over $\fQ$, with the exception
  of the indices $C^{(2)}_{2n,n}$ for the zero-dimensional boundary
  components.
\end{theorem}

We now describe the standard symmetric subgroups $N_b'$ incident to $P_b'$
for the classical cases. For details, see \cite{sym}. We consider the
vector space $V$ with the $\pm$symmetric/hermitian form $h$. In the
notation of (\ref{E8.c}), if the standard hermitian Levi factor $L_b'$ is
$L_b'=\cN_G(W)/\cZ_G(W)$, $W=V_{b+1}\oplus \cdots \oplus V_s\oplus V'$ in
the notations used there, then for $b<s$ or $c(s,\gs_i)<t_i$ for some
$i=1,\ldots, f$,
\begin{equation}\label{E12.1} N_b'=\cN_G(W).
\end{equation}
If $b=s$ and $c(s,\gs_i)=t_i$ for all $i=1,\ldots, f$, the boundary
component is a point, and $L_s'$ is the anisotropic kernel, and $N_s'$ as
in (\ref{E12.1}) is not the standard symmetric subgroup incident to $P_s'$
as we have defined it. Rather, these subgroups correspond to the following
constructions. We consider first the case where $\cD\not\in (\cE\cD)$. Pick
an anisotropic vector $v\in V$ which is defined over $k$, and consider the
subspace $W=v^{\perp}$, the space of vectors orthogonal to $v$. We describe
the subgroup $N_s'=\cN_{G'}(W)$, which depends on the choice of $v$.
\begin{itemize}\item[\ ] {\bf O.1:} $V$ is a $k$-vector space; the subgroup
   $N_s'$ gives rise to a
  subdomain $\cD_{N_s'}$ of type $\bf IV_{\hbf{n-1}}\inn IV_{\hbf{n}}$.
\item[\ ] {\bf O.2:} Here $V$ is an $n$-dimensional $D$-vector space, and
  we have $n$ odd; $W\inn V$ is of codimension one over $D$, giving rise to
  a subdomain of type $\bf II_{\hbf{n-1}}\inn II_{\hbf{n}}$.
\item[\ ] {\bf U.1:} In this case we get subdomains $\bf
  I_{\hbf{p-1,q}}\inn I_{\hbf{p,q}}$.
\item[\ ] {\bf U.2:} $V$ is $n$-dimensional over $D$, where $D$ has degree
  $d$ over $K$; the subspace $W$ gives rise to a subdomain of type $\bf
  I_{\hbf{p-d,q}}\inn I_{\hbf{p,q}}$. Iteration of this gives subdomains of
  types $\bf I_{\hbf{p-jd,q}}\inn I_{\hbf{p,q}}$, and for $j=s$ the boundary
  component will be a point $\iff$ $sd=q$, in which case $\bf
  I_{\hbf{q,q}}\inn I_{\hbf{p,q}}$ is a maximal tube domain and fulfills
  2').
\end{itemize}
Finally we consider $\cD\in (\cE\cD)$. In these cases, if the
zero-dimensional boundary component is rational, then $V$ splits into a
direct sum of hyperbolic planes (no anisotropic kernel). We can define a
unique polydisc by the prescription: letting $V=\oplus_{i=1}^s V_i$ be the
decomposition into hyperbolic planes as above, set:
\[ N_s':=\{g\in G' | g(V_i)\subseteq V_i,\ i=1,\ldots,s\}. \]
For the individual cases this gives rise to the following subdomains:
\begin{itemize}\item[\ ] {$\bf I_{\hbf{q,q}}$:} $\cD_{N_s'}\cong \bf
  I_{\hbf{d,d}} \times \cdots \times I_{\hbf{d,d}}$. In each of the factors
  $\bf I_{\hbf{d,d}}$ we can apply the results of \cite{hyp} to get a uniquely
  determined polydisc.
\item[\ ] {$\bf II_{\hbf{n}}$, $n$ even:} In this case we get a subdomain
  $\bf II_{\hbf{2}}\times \cdots \times II_{\hbf{2}}$, which is a
  polydisc, as $\bf II_{\hbf{2}}$ is a disc.
\item[\ ] {$\bf III_{\hbf{n}}$:} In case $\bf S.1$, the result is well known,
  giving just a polydisc. In case $\bf S.2$, we get as a subdomain $\bf
  III_{\hbf{2}}\times \cdots \times III_{\hbf{2}}$, and this case
  represents the exception in Theorem \ref{t12.1}; in general no
  polydisc (defined over $k$) can be found in each factor.
\end{itemize}

\section{Arithmetic groups}\label{s83.1}
By definition, an arithmetic subgroup $\gG\inn G(\fQ)$ is one which is
commensurable with $\grr^{-1}(GL(V_{\fZ}))\cap G(\fQ)$, for some ($\iff$
for any) faithful rational representation $\grr:G\lra GL(V)$, where $V$ is
a finite-dimensional $\fQ$-vector space, and $V_{\fZ}$ is a
$\fZ$-structure, i.e., a $\fZ$-lattice such that
$V_{\fZ}\otimes_{\fZ}\fQ=V_{\fQ}$. In the classical cases, it is natural to
take the fundamental representation as $\grr$ (more precisely the
fundamental representation $\grr':G'\lra GL_D(V)$ determines $\grr:G\lra
Res_{k|\fQ}GL_D(V)$), and for the exceptional structures, one has either
representations in exceptional Jordan algebras and related algebras, or
simply the adjoint representation.

Consider first the classical groups. For these, $D$ is a central simple
division algebra over $K$, where $K$ is either the totally real number
field $k$ or an imaginary quadratic extension of $k$, and $D$ has a
$K|k$-involution. The rational vector space $V$ is an $n$-dimensional right
$D$-vector space, $A=M_n(D)$ is a central simple algebra over $K$ with a
$K|k$-involution extending the involution on $D$ by (\ref{e49.2}). We have
a $\pm$symmetric/hermitian form $h:V\times V\lra D$ such that
\begin{equation}\label{E212} G'=\{g\in GL_D(V) \big| \forall_{x,y\in V},\
  h(x,y)=h(gx,gy)\}
\end{equation}
is the unitary group of the situation. We take the natural inclusion given
by (\ref{E212}), $\grr':G'\lra GL_D(V)$ and let the representation
$\grr:G\lra Res_{k|\fQ}GL_D(V)$ determined by $\grr'$ be our rational
representation. We now consider $\fZ$-structures on $V$, for which we
require an {\it order} $\gD\inn D$, i.e., a lattice that is a subring of
$D$, and consider $\gD$-lattices $\scL\inn V$. The analog of (\ref{E212}),
after fixing the $\fZ$-structure on $V$, is
\begin{equation}\label{E212.1}
\gG_{\hbox{\sscrpt L}}=\{g\in G \Big| g\scL\subseteq
  \scL \}.
\end{equation}
Then $\gG_{\hbox{\sscrpt L}}\inn G(\fQ)$ is an arithmetic subgroup, as it
is the set of elements which preserve the $\gD$-structure on $V$ defined by
$\scL$, which itself is a $\fZ$-lattice in the rational vector space $V$
(viewing $V$ as a $\fQ$-vector space). If, for example, $\gG_{\hbox{\sscrpt
    L}'} \inn \gG_{\hbox{\sscrpt L}}$ is a normal subgroup of finite index,
we get an induced representation of $\gG_{\hbox{\sscrpt
    L}}/\gG_{\hbox{\sscrpt L}'}$ in $\scL/\scL'$, where $\scL'$ is the
sublattice of $\scL$ preserved by $\gG_{\hbox{\sscrpt L}'}$. This is the
general formulation of an occurance which is well-known in specific cases.
For example, if $\gG_{\hbox{\sscrpt L}'}=\gG(N)\inn
Sp(2n,\fZ)=\gG_{\hbox{\sscrpt L}}$ is the principal congruence subgroup of
level $N$, there is a representation of $\gG/\gG(N)\cong Sp(2n,\fZ/N\fZ)$
in $(\fZ/N\fZ)^{2n}\ (=\scL/\scL')$.

Now consider the exceptional groups. In the case of $E_6$ we have the
27-dimensional representation in the exceptional Jordan algebra $\JJ$,
while in the case of $E_7$ we have the 56-dimensional representation in the
exceptional algebra of $2\times 2$ matrices over $\JJ$\footnote{This is
  what W.~Baily utilized in his beautiful paper \cite{Ba}.}. In both cases
we can also use the adjoint representation, so we require a $\fZ$-structure
on the Lie algebra itself. Such can be readily constructed, utilizing the
Tits algebras, from lattices in the constituents, composition algebras and
(exceptional) Jordan algebras.

After these introductory remarks we proceed to give a few details, which in
particular allow us to give some relevant references in each case. We start
by discussing orders, then describe the arithmetic groups these give rise
to.

\subsection{Orders in associative algebras}
A general reference for this section is \cite{reiner}. We first fix
some notations. $k$ is a totally real Galois extension of degree $f$ over
$\fQ$, and $\cO_k$ will denote the ring of integers in $k$. $D$ will denote a
division algebra (skew field), central simple of degree $d$ over $K$,
with a $K|k$ involution
($K=k$ for involutions of the first kind, and $K$ is an imaginary quadratic
extension of $k$ for involutions of the second kind).
$V$ denotes an $n$-dimensional
right $D$-vector space, so that $Hom_D(V,V)\isom M_n(D)$. $A=M_n(D)$ is a
central simple algebra over $K$ with involution extending the involution on
$D$ by
\begin{equation}\label{e49.2} M \mapsto M^*,\hbox{ where
    }(M^*)_{ij}=\overline{m}_{ji},\hbox{ for $M=(m_{ij})$},
\end{equation}
where ``--'' denotes the involution in $D$. $(V,h)$ is a $\pm$-hermitian space
with $\pm$-hermitian form $h$
(with respect to the involution on $D$). Hence $[D:K]=d^2,\ [A:K]=(nd)^2=t^2,
t=nd$.

Let $F$ be a number field, for example $F=K,k$ as above, and
let $W$ be an $F$-vector space. A full $\cO_F$-{\em lattice} $\cL$ in $W$ is an
$\cO_F$-module, finitely generated, such that $F\cdot \cL=W$. Usually we work
with full lattices and delete the word full. If $W$ is an
$F$-algebra, then an $\cO_F$-lattice $\cL$ is an $\cO_F$-{\em order}, if
$\cL$ is a subring of $W$. In particular in $W=D$, an $\cO_F$-order is a
(full) lattice which is a subring. Let $\gD\inn D$ denote an order in $D$,
and let $V$ be an $n$-dimensional vector space over $D$. Then a (full)
$\gD$-{\em lattice} in $V$ is a $\gD$-module $\cM$ with $\cM\cdot D=V$; if
again $A$ is the algebra $M_n(D)$, then a $\gD$-lattice in $A$ is a
$\gD$-{\em order}, if it is a subring of $A$.

Let an $\cO_F$-lattice $\cL\inn A$ be given. $\cL$ determines a right
(respectively left) $\cO_F$-order:
\begin{equation}\label{e73.1}\cO_r(\cL)=\{x\in A \Big| \cL\cdot x\inn \cL\},
(\hbox{respectively } \cO_l(\cL)=\{x\in A\Big| x\cdot \cL\inn \cL\}).
\end{equation}
If $\cL$ is a $\gD$-lattice, then $\cO_r(\cL)$ and $\cO_l(\cL)$ are
$\gD$-orders. If an $\cO_F$-order $\cO\inn A$ is given, and $\cL\inn A$ is a
lattice with $\cO=\cO_r(\cL)$ (respectively $\cO_l(\cL)$), then one also
calls $\cL$ an {\it $\cO$-lattice}, and says that $\cL$ and $\cO$ are {\em
associated}. An element $a\in A$ is called {\em integral}, if its
characteristic polynomial has integer coefficients,
$\chi_a\in\cO_F[X]$. It is a basic result that
every element $a\in \cO$ is integral for any $\cO_F$-order $\cO$ in $A$. An
order $\cO$ is {\em maximal}, if it is not properly contained in any other
order. It is a basic fact that maximal orders exist in $D$ and in $A$, and
that any order is contained in a maximal one (\cite{reiner}, 10.4).
One has the following description of maximal orders in $A$:
\begin{theorem}[\cite{reiner}, 21.6]\label{t74.1}
Notations as obove, let $\gD\inn D$ be a fixed maximal $\cO_F$-order in $D$,
and let $\cM$ be any (full) right $\gD$-lattice in $V$. Then
$Hom_{\gD}(\cM,\cM)$ is a maximal $\cO_F$-order in $A$, and for any maximal
$\cO_F$-order $\cO$ in $A$, there exists a (full) right $\gD$-lattice
$\cN\inn V$ with $\cO=Hom_{\gD}(\cN,\cN)$.
\end{theorem}
The following result of Chevally describes maximal orders
in associative algebras.
\begin{theorem}[\cite{reiner}, 27.6]\label{t78.1}
Let $\gD\inn D$ be a maximal $\cO_K$-order in $D$; for each right ideal
$J\inn \gD$, set $\gD'=\cO_l(J)$. Then every maximal order of $A=M_n(D)$ is
of the form
$$\cO_J=\left( \begin{array}{cccc} \gD & \cdots & \gD & J^{-1} \\
\vdots & \ddots & \vdots & \vdots \\
\gD & \cdots & \gD & J^{-1} \\
J & \cdots & J & \gD' \end{array}\right),$$
for some right ideal $J$, and for each $J$, the lattice $\cO_J$ above is a
maximal order.
\end{theorem}
In other words, to give a maximal order in $A$ is the same as giving a
maximal order $\gD\inn D$, together with a right ideal $J\inn \gD$, i.e., the
same as giving a pair $(\gD,J)$. In particular if the class number $h(\gD)=1$
(note that $h(\gD)$, which is defined as the number of left $\gD$-ideal
classes, is also equal to the number of right $\gD$-ideal classes,
see \cite{reiner}, Ex. 7 iii), p. 232), then up to $\gD$-isomorphism there is
a 1-1 correspondence between isomorphism classes of
maximal orders in $D$ and $A$.
\subsection{Orders in Jordan algebras}
First recall the result on orders in the (definite) Cayley algebra from
\cite{BS}. Let $e_0,\ldots,e_7$ be the base of $\CC_{\fQ}$ given as follows
\[\begin{minipage}{12cm}
$\CC_{\fQ}= e_0\fQ+ e_1\fQ+\ldots+ e_7\fQ,$ with center $e_0\fQ$ and
relations:\nz $e_i\cdot e_{i+1}=e_{i+3}, e_{i+1}\cdot e_{i+3}=e_i,
e_{i+3}\cdot e_i=e_{i+1}, e_i^2=-e_0,\  i\in \fZ/7\fZ$.
\end{minipage}\]
Define
\begin{equation}\label{e80a.1} \scM:=\{x=\sum\xi_ie_i\Big| 2\xi_i\in \fZ,
\xi_i-\xi_j\in \fZ, \sum \xi_i\in 2\fZ\}.
\end{equation}
Then
\begin{lemma}[\cite{BS}, 4.6]\label{l80a.1} $\scM$ is a maximal order in
$\CC_{\fQ}$, and any other maximal order is isomorphic to $\scM$.
\end{lemma}
(In \cite{BS} the authors call $\scM$ an octave-ring: a subring of
$\CC_{\fQ}$ containing 1, on which the norm form is integral, and maximal
with these properties; we just call $\scM$ a maximal order.)

A general reference for the remainder of
this section is \cite{racine}. Let $R$ be a
commutative ring. A {\em Jordan algebra} over $R$ is an $R$-module which is
commutative and satisfies the relation
\[(x^2\cdot y)\cdot x = x^2\cdot (y\cdot x),\
\forall_{x,y}.\]

\vspace*{-.3cm}
\begin{definition}\label{d74a.1} Let $\JJ$ be a Jordan algebra over a number
field $K$, and let $\cO_K$ denote the ring of integers in $K$. A full
$\cO_K$-lattice $\cL\inn \JJ$ is an {\em order}, if $\cL$ is a Jordan
algebra over $\cO_K$.
\end{definition}
An element $x\in \JJ$ is {\em integral}, if the characteristic polynomial
is integral,
i.e., if $N(x), Q(x)$ and $T(x)$ are integral (see \cite{Jac}, pp.~91, also
\cite{Jacob}, Chapter VI, for
details).
Let $\cL$ be an order in
$\JJ$, and $x\in \cL$; then $\cO_K[x]\inn \cL$ is an associative subalgebra,
hence finitely generated,
so $x$ is integral (\cite{racine} Prop.~1, p.~19). Conversely, any
integral element of $\JJ$ is contained in an order ({\it loc.~cit.}~Prop.~2).

Once again it is a basic fact that maximal orders exist
({\it loc.~cit.}~Thm.~2) and that an order is maximal if and
only if it is maximal locally
everywhere ({\it loc.~cit.}~Lemma 1).
A maximal order $\cL\inn \JJ$ is said to be {\em distinguished}, if $\cL$ is
a maximal lattice of integral elements. For example, if $\cO\inn \fO$ is a
maximal order in an octonion algebra, then $J(\cO,\gg)\inn J(\fO,\gg)$
(notations as in \ref{d46b.1}) is a distinguished maximal
order. Conversely,
for $\gg=1$,
\begin{proposition}[\cite{racine}, Prop.~5, p.~115]\label{p74a.1}
If $\JJ=J(\fO_K,1)$ is the exceptional Jordan
algebra over the totally indefinite octonion algebra $\fO_K$,
then any distinguished order $\cP\inn
\JJ$ is of the form $J(\cO,1)\inn \JJ$, with $\cO$ a maximal order in
$\fO_K$.
\end{proposition}
This may be considered in some sense as an analogue of Theorem \ref{t74.1}
for orders in exceptional Jordan algebras.

\subsection{Lattices in Tits algebras}
Let $\AA$ be a composition algebra over $K$, and
$\JJ=J(\AA',1)$ a Jordan algebra as in \ref{d46b.1} over a second composition
algebra $\AA'$. For a totally indefinite octonion
algebra over $K$, $\AA'$, and a maximal order $\gD'\inn \AA'$,
then, as we have seen
(Proposition \ref{p74a.1}), $\cL=J(\gD',1)$ is a distinguished order in
$\JJ$ and conversely. More generally it is easy to see:
\begin{lemma}\label{l74b.1} Let $\gD\inn \AA$ be a maximal order in the
composition algebra $\AA$. Then $J(\gD,\gg)$ is a maximal order in the Jordan
algebra $J(\AA,\gg)$ of Definition \ref{d46b.1}.
\end{lemma}
{\bf Proof:} $\cL:=J(\gD,\gg)$ is clearly a
Jordan algebra over $\cO_K$, hence it
is an order in $\JJ$. To see it is maximal, the method of \cite{racine} can
be used. Let $L_1=\{a\in \fO_K \big| a[j,k]\in \cL\}$ (notations as in
(\ref{e46d.1})); this is a lattice in
$\fO_K$, and, as can be checked,
is the lattice $\gD$ which we started with. If
$\cL$ is not maximal, then $\cL\subsetneqq \cL'$,
and the corresponding $L_1'$ will
be an $\cO_K$-lattice in $\AA'$ with $L_1\subsetneqq L_1'$,
contradicting the maximality of $\gD$. \ende
Let $\gD\inn\AA$ be a maximal order and $\cL\inn \JJ$ a maximal order.
Consider the Tits algebra $\scL(\AA,\JJ)$ of Definition \ref{d46f.1}.
Recall that the construction of Tits algebras requires, in addition to the
algebras $\AA$ and $\JJ$, also the Lie algebras $Der(\AA)$ and
$Der(\JJ)$. If we have maximal orders $\gD\inn \AA,\ \cL \inn \JJ$, then we
{\it define}:
\[Der(\gD):= \{ D\in Der(\AA) \Big| D(\gD)\subseteq \gD\},\quad Der(\cL):=
\{ D\in Der(\JJ) \Big| D(\cL)\subseteq \cL\}.\]
Since we know that $Der(\AA)$ is a Lie algebra of type $G_2$ and $Der(\JJ)$
is a Lie algebra of type $F_4$, we are asking for $\fZ$-structures on these
Lie algebras. Clearly $Der(\gD)$ and $Der(\cL)$ are lattices in the
corresponding Lie algebras, which are furthermore closed under the Lie
bracket. It then is natural to consider the following lattice in the Tits
algebra:
\begin{equation}\label{E210A}
\gL_{\gD,\cL}:=Der(\gD)\oplus \gD_0\otimes \cL_0\oplus Der(\cL),
\end{equation}
and the corresponding arithmetic group it defines (for
$G=\Aut(\scL(\AA,\JJ))^0$)
\begin{equation}\label{E210B}
\gG_{\gd,\cL}:= \{ g\in G \Big| \ad(g)(\gL_{\gD,\cL})\subseteq
\gL_{\gD,\cL}\},
\end{equation}
where $G$ is acting by means of the adjoint representation on
$\scL(\AA,\JJ)$.

\subsection{Arithmetic groups -- classical cases}\label{sarithmetic}
In this subsection $G'$ will denote an absolutely (almost)
simple $k$-group ($k$ a
totally real number field) which we assume is classical,
$G=Res_{k|\fQ}G'$ the $\fQ$-simple group it
defines, which we assume is of hermitian type. We let $\grr':G'\lra GL_D(V)$
be the natural inclusion and $\grr:G\lra Res_{k|\fQ}GL_D(V)$ be the natural
representation of $G$ defined by $\grr'$. Fix a maximal order $\gD\inn D$,
and let $\cL\inn V$ be a $\gD$-lattice (which is in particular a
$\fZ$-lattice of the underlying $\fQ$-vector space). As above,
$\cO_r(\cL)$ (respectively $\cO_l(\cL)$) will denote the right
(respectively left) order of $\cL$, given by the equation
(\ref{e73.1}). First of all, we have the arithmetic subgroup
$GL_{\gD}(\cL)\inn GL_D(V)$, and we define the subgroup
\[\gG_{\cL}':= \{g\in G'(k) \Big| \grr(g)(\cL)\subseteq \cL\} =
\grr^{-1}(GL_{\gD}(\cL))\inn G'(k),\]
and similarly $\gG_{\cL}\inn G(\fQ)$. By definition these are arithmetic
subgroups of $G'(k)$ and $G(\fQ)$, respectively. Let us see how this is
related to the orders $\cO_r(\cL)$ and $\cO_l(\cL)$. By Theorem \ref{t78.1}
$\cO_r(\cL)$ is of the form $\cO_J$ for some right ideal $J\inn \cO$.
Our central simple
algebra is in this case $A=M_n(D)$, and $\cO_r(\cL)$
is a maximal order in $A$. Recall how the group $G'$ and the algebra are
related (\cite{W}, Thm.~2, p.~598). Let $U=\{z\in A | z z^*=1\}$, $U_0$ the
connected component of $U$, ($(G')^0=$)$G_0:=(\Aut(A))^0$, and let $C\inn U_0$
be the center of $U_0$. Then we have an exact sequence
\[1\lra C\lra U_0 \lra G_0 \lra 1.\]
As a lattice in $A$ we consider $\cO:=\cO_r(\cL)$ and its intersection with
$U_0$,
\[ \cO_0=\cO\cap U_0.\]
Similarly, $\cC:=C\cap \cO_0$ is the center of $\cO_0$, and we have the
sequence
\[1\lra \cC\lra \cO_0 \lra \gG'\lra 1,\]
where $\gG'\cong \cO_0/\cC$ is the arithmetic subgroup $\gG'\inn G'(k)$,
showing how the maximal orders are related to the arithmetic groups.
In our situation here, $(G')^0$ plays the role of
$G_0$, while $(U')^0=(\{z\in A \Big| z^*z=1\})^0$ plays the role of $U_0$. Let
further $C'\inn (U')^0$ be the center. We have $\cO\cong \cO_J$ for some
right ideal $J$, and $(\cO_J)_0=\cO_J\cap U_0$ plays the role of $\cO_0$.
Then $C'\cap
(\cO_J)_0=\cC$ is the center of $(\cO_J)_0$, and we have sequences:
\[\begin{array}{ccccccccc} 1 & \lra & C' & \lra & (U')^0 & \lra & (G')^0 &
  \lra & 1 \\
   & & \cup & & \cup & & \cup \\
  1 & \lra & \cC & \lra & (\cO_J)_0 & \lra & \gG_{\cL}'
   & \lra &
  1. \end{array}\]
In this sense, maximal orders give rise to arithmetic subgroups. Viewing
the $\cO_k$-lattice $\cL$ as a $\fZ$-lattice gives the corresponding
diagram for the $\fQ$-groups (with hopefully obvious notations)
\[\begin{array}{ccccccccc} 1 & \lra & C & \lra &  U^0 & \lra & G^0 & \lra &
  1 \\
   & & \cup & & \cup & & \cup \\
  1 & \lra & \cZ(\cO_0) & \lra & \cO_0 & \lra & \gG_{\cL} & \lra & 1.
\end{array}\]
We now describe this more precisely for the following special cases:
\begin{itemize}\item[a)] Siegel modular groups. \item[b)] Picard modular
  groups. \item[c)] Hyperbolic plane modular groups.
\end{itemize}
These are examples of $\fQ$-groups which are of both inner type (for a))
and outer type (for b) and c)), of split over $\fR$-type, meaning the
$\fQ$-rank is equal to the $\fR$-rank (for a) and b)) and more or less the
{\it opposite} of split over $\fR$-type ($\fQ$-rank equal to one,
$\fR$-rank unbounded) (for c)). Case a) is well-known, b) is also to a
certain extent, while c) was introduced in \cite{hyp}.
\begin{itemize}\item[a)] Siegel case:
\begin{itemize}\item $A=M_{2n}(\fQ)$ with the involution $*:X\mapsto
  JX^tJ,\quad J=\left(\begin{array}{cc}0 & \boldone  \\ -\boldone & 0
    \end{array}\right)$.
\item $\Aut(A,*)\cong PSp(2n,\fQ),\quad V=\fQ^{2n}$.
\item $D=\fQ$, a maximal order is $\gD=\fZ,\ \ V_{\fZ}=\fZ^{2n}$.
\item $\gG=PSp(2n,\fZ)$.
\end{itemize}
The sequence above becomes:
\[ \begin{array}{ccccccccc} 1 & \lra & \fZ/(2) & \lra & Sp(2n,\fQ) & \lra &
  PSp(2n,\fQ) & \lra & 1 \\ & & || & & \cup & & \cup \\
1 & \lra & \fZ/(2) & \lra & Sp(2n,\fZ) & \lra & \gG & \lra & 1.
\end{array}\]
\item[b)] Picard case:
\begin{itemize}\item $A=M_n(K)$ with involution $*:X\mapsto HX^tH$, $H$
  hermitian, where $K|\fQ$ is imaginary quadratic.
  \item $\Aut(A,*)\cong PSU(K^n,h),\quad V=K^n$, $h$ is a hermitian form
    represented by $H$.
 \item $D=K$, a maximal order is $\gD=\cO_K,\quad V_{\fZ}=\cO_K^n$.
 \item $\gG=PSU(\cO_K^n, h)$ (or $PU(\cO_K^n,h)$, which is not simple, but
   is often considered anyway).
 \end{itemize}
The sequence above becomes:
\[\begin{array}{ccccccccc} 1 & \lra & \cC & \lra & SU(K^n,h) & \lra &
  PSU(K^n,h) & \lra & 1 \\
& & \cup & & \cup & & \cup \\
1 & \lra & \scC & \lra & SU(\cO_K^n,h) & \lra & PSU(\cO_K^n,h) & \lra & 1.
\end{array}\]
Note that $\cC$ is given essentially by $\cO_K\cap U(1)$, which is $\pm1$
except for the two fields $K=\fQ(\sqrt{-1}),\ K=\fQ(\sqrt{-3})$ which
contain fourth (respectively third) roots of unity.
\item[c)] Hyperbolic plane case:
\begin{itemize}\item $D$ is a division algebra, central simple of
  degree $d\geq2$ over $K$,
  with a $K|\fQ$-involution, $\gD\inn D$ is a maximal order.
\item $A=M_2(D)$ with involution $*:X\mapsto
  {^t\overline{X}}$, where ${^t\overline{X}}=(\overline{x}_{ji})$, if
  $X=(x_{ij})$, and $\overline{x}$ denotes the involution in $D$.
\item $\Aut(A,*)$ is a $\fQ$-form of $PSU(d,d)$, and $V=D^2,$ with a
  hermitian form $h:V\times V\lra D$ which is isotropic, $V_{\fZ}=\gD^2$.
\item $\gG=PSU(\gD^2,h)$.
\end{itemize}
The above sequence becomes in this case
\[\begin{array}{ccccccccc} 1 & \lra & \cC & \lra & SU(D^2,h) & \lra &
  PSU(D^2,h) & \lra & 1 \\ & & \cup & & \cup & & \cup \\
1 & \lra & \cC\cap \gD & \lra & SU(\gD^2,h) & \lra & PSU(\gD^2,h) & \lra & 1.
\end{array}\]
As $D$ is central simple over $K$, the center is as in the last case,
$\cC\cong \cO_K\cap U(1)$, hence it is $\pm 1$ except
for the case $K=\fQ(\sqrt{-1})$ and $K=\fQ(\sqrt{-3})$ as above.
\end{itemize}

\subsection{Arithmetic groups -- exceptional cases}
We mentioned above that for the exceptional cases, there are (at least) two
natural types of representations we can consider: representations in
algebras derived from exceptional Jordan algebras (Tits algebras), and the
adjoint representation. These representations correspond to the following
fundamental weights:

\vspace*{.5cm}
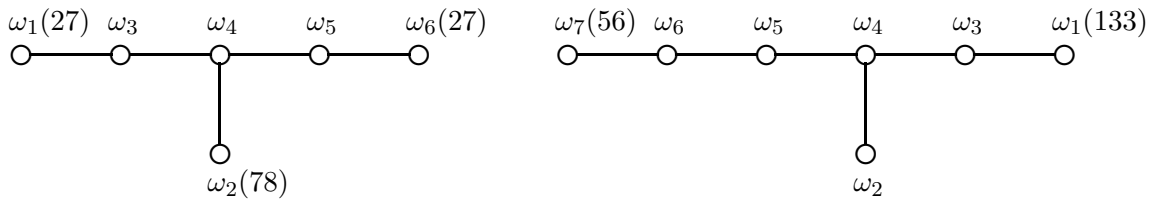
\[
\setlength{\unitlength}{0.006500in}%
\begin{picture}(854,94)(33,693)
\thicklines
\put( 40,780){\circle{14}}\put( 30,800){$\go_1 (27)$}
\put(120,780){\circle{14}}\put(110,800){$\go_3$}
\put(200,780){\circle{14}}\put(190,800){$\go_4$}
\put(280,780){\circle{14}}\put(270,800){$\go_5$}
\put(360,780){\circle{14}}\put(350,800){$\go_6 (27)$}
\put(200,700){\circle{14}}\put(190,670){$\go_2 (78)$}

\put(480,780){\circle{14}}\put(470,800){$\go_7 (56)$}
\put(560,780){\circle{14}}\put(550,800){$\go_6$}
\put(640,780){\circle{14}}\put(630,800){$\go_5$}
\put(720,780){\circle{14}}\put(710,800){$\go_4$}
\put(800,780){\circle{14}}\put(790,800){$\go_3$}
\put(880,780){\circle{14}}\put(870,800){$\go_1 (133)$}
\put(720,700){\circle{14}}\put(710,670){$\go_2$}

\put(200,773){\line( 0,-1){ 66}}
\put(487,780){\line( 1, 0){ 66}}
\put(567,780){\line( 1, 0){ 66}}
\put(807,780){\line( 1, 0){ 66}}
\put(720,773){\line( 0,-1){ 66}}
\put(727,780){\line( 1, 0){ 66}}
\put(647,780){\line( 1, 0){ 66}}
\put( 47,780){\line( 1, 0){ 66}}
\put(127,780){\line( 1, 0){ 66}}
\put(207,780){\line( 1, 0){ 66}}
\put(287,780){\line( 1, 0){ 66}}
\end{picture}\]
In the case of $E_6$, the 27-dimensional (respectively the adjoint,
78-dimensional) representation corresponds to the weights $\go_1$ and
$\go_6$ (respectively to $\go_2$), while in the case of $E_7$, the
56-dimensional (respectively the adjoint, 133-dimensional) representation
corresponds to the weight $\go_7$ (respectively to $\go_1$). We briefly
discuss the arithmetic groups arising in this way.

We first consider the 27-dimensional representation. For this we assume
$G'$ has index ${^2E}_{6,2}^{16'}$ and we use the model
\[\Gg'=\cL(\JJ)_{\gl}=\sqrt{\gl}R_{\JJ_0}\oplus Der(\JJ)\]
(Albert's twisted $\cL(\JJ)$),
where $\gl<0, \gl\in k$. We then choose a maximal order
$\cM\inn \JJ_k$ and set
\[Der(\cM)=\{a\in Der(\JJ) \Big| a(\cM)\subseteq \cM\}.\]
Then we may consider the lattice
\[\cL(\cM)_{\gl}:=\sqrt{\gl}R_{\cM_0}\oplus Der(\cM).\]
This defines an arithmetic group:
\[\gG_{\cM}:=\{g\in G' \Big| \grr(g)(\cL(\cM)_{\gl})\subseteq
\cL(\cM)_{\gl}\},\]
where $\grr$ is the 27-dimensional representation in $\cL(\JJ)_{\gl}$.

Next we consider the adjoint representation. For this we utilize the
lattice in the Tits algebra constructed in (\ref{E210A}), and the
corresponding arithmetic group (\ref{E210B}). That lattice depends on the
choice of a maximal order $\gD$ in the Cayley algebra, as well as on one in
the algebra $\BB$. More explicitly,
\begin{theorem}\label{t80b.1}
Let $\Gg'$ be a $k$-form of $\ee_{6(-14)}$ as in Corollary
\ref{c55.1}, $G'$ as in Corollary \ref{c55.2}, i.e.,
$$\Gg'\isom \scL(\CC_k, (J_1)_k^b),\quad (G')^0\isom (\Aut(\Gg'))^0.$$
Let $\gD\inn \CC_k$ be a maximal order in the Cayley algebra $\CC_k$ as
above, let $\cL\inn (J_1)_k^b$ be a maximal order in the Jordan algebra
(Definition \ref{d74a.1}), and set
$$\Gg_{(\gD,\cL)}':=\scL(\gD,\cL)=\left\{ X\in \scL(\CC_k,(J_1)_k^b) \Big|
\parbox{6cm}{$X=X_1+x\otimes y + Y_1: X_1\in Der(\gD),$

$Y_1\in Der(\cL), x\in
\gD_0, y\in \cL_0$ } \right\}.$$
Then $\Gg_{(\gD,\cL)}'$ is an $\cO_k$-lattice in the $k$-vector space $\Gg'$.
Set
$$\gG_{(\gD,\cL)}':=\{g\in G'(k) \Big| \ad(g)(\Gg_{(\gD,\cL)}')\inn
\Gg_{(\gD,\cL)}'\}.$$
Then $\gG_{(\gD,\cL)}'\inn G'(k)$ is an arithmetic subgroup.
\end{theorem}

Now consider type $E_7$. We first consider the 56-dimensional
representaion. This is the situation considered by Baily in \cite{Ba}.
In this
example $k=\fQ$, and $\JJ$ is the exceptional Jordan algebra over $\fQ$,
$\JJ_{\fQ}=J^b=J(\CC_{\fQ},(1,-1,1))$ in the notation of Definition
\ref{d46b.1}, and $\AA_{\fQ}=M_2(\fQ)$.
Let $\scM\inn \CC_{\fQ}$ be the maximal order
(\ref{e80a.1}). This determines, as in \ref{l74b.1}, a maximal order $\cL$ in
$\JJ_{\fQ}$. Also $\fZ\inn \fQ$ defines the maximal order $\gD=M_2(\fZ)\inn
M_2(\fQ)$. This then gives rise to an arithmetic group $G_{(\gD,\cL)}$, which
Baily shows is maximal and has only one cusp.

Again in this case we can also consider the adjoint representation. For
this we again utilize the lattice (\ref{E210A}), and as above, this
determines an arithmetic group as in (\ref{E210B}). This time, we need a
lattice in the totally indefinite quaternion algebra $\AA$ as well as one in
the Jordan algebra $(J_1)_k^b$. More explicitly,
\begin{theorem}\label{t80b.2}
Let $\Gg'$ be a $k$-form of $\ee_{7(-25)}$ as in Theorem \ref{t55a.2}, i.e.,
$$\Gg'\isom \scL(\AA_k,\JJ_k),\quad (G')^0\isom (\Aut(\Gg'))^0.$$
Let $\gD\inn \AA_k$ be a maximal order in the indefinite quaternion algebra
$\AA_k$ as in section 3.1, and let $\cL\inn \JJ_k$ be a maximal
order in the Jordan algebra $\JJ_k$ as in \ref{d74a.1}, and set:
$$\Gg_{(\gD,\cL)}':=\scL(\gD,\cL) \hbox{ as above }. $$
Then $\Gg_{(\gD,\cL)}'\inn \Gg'$ is an $\cO_k$-lattice, and
$$G_{(\gD,\cL)}':=\{g\in \Aut(\Gg')\Big| \ad(g)(\Gg_{(\gD,\cL)}')\inn
\Gg_{(\gD,\cL)}'\}\cap (G')^0$$
is an arithmetic subgroup in $G'(k)$.
\end{theorem}
A more detailed discussion of these matters cn be found in \cite{new}.

\section{Integral symmetric subgroups}
Let $G$ be a $\fQ$-simple algebraic group of hermitian type, and let $A\inn
G$ be a maximal $\fR$-split torus defined by the set of strongly orthogonal
roots as in section 1.1, given the canonical order.  Let $S\inn A$ be a
maximal $\fQ$-split torus with the canonical order, compatible with the
given order on $A$. Let further $\gD_{\fQ}=\{\eta_1,\ldots, \eta_s\}$ be
the set of simple $\fQ$-roots, and let $F_{\bfb},\ P_{\bfb}$ be the
standard boundary components and parabolics as explained above, ${\bf
  b}=(c(b,\gs_1),\ldots,c(b,\gs_f)), b=1,\ldots, s$. Finally, let
$N_{\bfb}$ be the standard incident symmetric subgroup (i.e., given by
(\ref{e10.1}) if $\dim(F_{\bfb})>0$, and in terms of root systems as
explained in section 2.2 for $\dim(F_{\bfb})=0$). Since $N_{\bfb}$ is a
reductive subgroup, it is not true that any $G$-conjugate $N'$ of
$N_{\bfb}$ is already $G_{\fQ}$-conjugate. Therefore we make the following
definition, yielding a proper subset of the set of $G$-conjugates of the
given $N_{\bfb}$.
\begin{definition}\label{D4.1} Let $G,\ S,\ P_{\bfb},\ N_{\bfb}$ be given as
  above. A symmetric subgroup $N'\inn G$ which is conjugate to $N_{\bfb}$
  by an element of $G(\fQ)$ is called a {\it rational symmetric} subgroup
  of $G$.
\end{definition}
The following well-known example illustrates the difference between
rational and more general symmetric $\fQ$-subgroups.
\begin{example}\label{example} Let $G'$ be the symplectic group $G'=Sp(V,h),\
  G=Res_{k|\fQ}G'$, where $V$ is a $k$-vector space of dimension $2n$ and
  $h$ is skew-symmetric. If $n=2$, the corresponding domain is a product of
  copies of the Siegel space of degree 2 (type $\bf III_{\hbf{2}}$). The
  boundary components corresponding to $P_{\hbf{1}}$ (respectively
  $P_{\hbf{2}}$) are products of one-dimensional (respectively
  zero-dimensional) boundary components. Then $N_{\hbf{1}}$ is also a
  product of two factors, $N_{\hbf{1}}=N_{\hbf{1},1}\times N_{\hbf{1},2}$,
  and each $N_{\hbf{1},i}$ is a polydisc $({\bf H})^f$. If we consider the
  universal family of abelian varieties parameterized by the domain $\cD$,
  say $\cA\lra \cD$, we may consider the following conditions on the fibres
  $A_t\in \cA$ ($t\in \cD$):
\begin{itemize}\item[1)] $A_t$ is isogenous to a product.
\item[2)] $A_t$ is simple with real multiplication by some real quadratic
  extension $k'|k$.
\end{itemize}
We claim that the locus 1) is the locus of subdomains $\cD_{N'}$, where
$N'$ is rational symmetric, while the locus 2) is the union of $\cD_{N'}$,
where $N'$ is a $\fQ$-subgroup conjugate to $N_{\hbf{1}}$, but not in
$G(\fQ)$. To see this, let us suppose $k=\fQ$; we have the familiar
description for the domains $\cD_{N'}$ of 2): in this case the standard
symmetric subdomain is $\fS_1\times \fS_1 \inn \fS_2$ (given by the
diagonal $2\times 2$ matrices), and it is conjugated (in $G_{\fR}$) by the
matrix
$$\scS=\left(\begin{array}{cc} S^{-1} & 0 \\ 0 & {^tS}
\end{array}\right),\quad S=\left(\begin{array}{cc} 1 & w \\ 1 & \overline{w}
\end{array}\right),\ w\in \cO_{k'}$$ for a real quadratic extension
$k'|\fQ$. More precisely, the subdomains $\cD_{N'}$ are given by the
equations
\begin{equation}\label{humbert} \HH_{(a,b,c,d,e)}:=\left\{
    \tau=\left(\begin{array}{cc}\tau_{11} & \tau_{12} \\ \tau_{12} &
        \tau_{22}
      \end{array}
    \right) \Big| a\tau_{11}+b\tau_{12}+c\tau_{22}+
    d(\tau_{12}^2-\tau_{11}\tau_{22})+e=0\right\},
  \end{equation}
  for some integral tuple $(a,b,c,d,e)$ with $c,d\equiv0(p)$ for some prime
  $p$. Then the {\it discriminant} is $\gD=b^2-4ac-4de$, and the field
  $k'=\fQ(\sqrt{\gD})$ is the field mentioned above; the element $w\in
  \cO_{k'}$ can be taken here, for example, as $w={b+\sqrt{\gD} \over 2}$,
  yielding a Humbert surface with $a=1,\ d=e=0$.  The standard symmetric
  subgroup $\isom SL(2,\fQ) \times SL(2,\fQ)$ gets conjugated onto groups
  $\cong SL(2,k')$ by the elements $\scS\in Sp(4,k')\inn Sp(4,\fR)$.  The
  rational boundary components of $\cD_{N'}$ which are $SL_2(k')$-rational,
  are zero-dimensional, and are also rational boundary components of $\cD$.
  Note that the domain $\cD_{N'}$ defined by the subgroup $SL(2,k')$ also
  contains one-dimensional cusps of the domain $\cD$, the normalizers of
  which are defined over $k'$, but not over $\fQ$ and these boundary
  components are consequently not rational (for either $G'$ or $N'$).

  It is clear that $\cD_{\gD}$, the union of the subdomains of given
  discriminant $\gD$, is the union of conjugates of the standard one by
  elements of $G(\fQ)$ if and only if $\gD$ is a square, giving 1). If
  $\gD$ is not a square, then $k'=\fQ(\sqrt{\gD})$, $N'$ is conjugate to
  $N_{\hbf{1}}$ by an element in $G(k')$, and these are the cases occuring
  in 2).\hfill $\Box$
\end{example}
Next we note that the set of subgroups defined in Definition \ref{D4.1} is
independent of the maximal $\fQ$-split torus used to define $N_{\hbf{1}}$;
if $S'$ is another it is conjugate in $G(\fQ)$ to $S$, and $N'$ will be
rational with respect to $S$ exactly when it is so with respect to $S'$.
For a fixed $N$ the set of rational symmetric subgroups conjugate to $N$ is
naturally identified with $\cH=G(\fQ)/N(\fQ)$. Since $G(\fQ)$ acts on $\cH$
so does any arithmetic subgroup $\gG\inn G(\fQ)$ and one can consider the
double coset space $\gG\bs \cH$. By definition, $N_{\bfb}^g$ and
$N_{\bfb}^{g'}$ will be in the same $\gG$-orbit if $g(g')^{-1}\in \gG$, so
the orbits are determined by the denominators occuring in $g$ and in
$(g')^{-1}$, respectively.  In the example above, for each prime $p$, the
group $N_{\gD}$, $\gD=p^2$ lies in a seperate $\gG$-orbit. In particular
there are in general infinitely many orbits. It turns out that the
following definition gives a convenient notion. For $b<t$ (by which we mean
$\dim(F_{\bfb})>0$) let $N_{\bfb}=N_{\hbf{b,1}}\times N_{\hbf{b,2}}$ be the
decomposition above, and for $N=N_{\bfb}^g$, let $N=N_1\times N_2$ denote
the corresponding decomposition. If $b=t$ we set $N_{\bfb}=N_{\hbf{b,1}},\
N=N_1$.
\begin{definition}\label{D5.1} Let $G,\ S,\ P_{\bfb},\ N_{\bfb}$ be fixed
  as above, $\gG\inn G(\fQ)$ arithmetic. A rational symmetric subgroup
  $N\inn G$, conjugate to $N_{\bfb}$ by $g\in G(\fQ)$,
  $N=N_{\bfb}^g:=gN_{\bfb}g^{-1}$ will be called $\gG$-{\it integral}
  (respectively {\it strongly $\gG$-integral}), if \[ N_1\cap
  \gG=g(N_{\hbf{b,1}}\cap \gG)g^{-1}\quad (\hbox{respectively } N \cap \gG
  = g(N_{\bfb}\cap \gG)g^{-1})\] for the element $g$ above.
\end{definition}
For $b=t$, both notions coincide, otherwise strongly $\gG$-integral implies
$\gG$-integral. For our purposes, the weaker notion will be most important.
Note that since $N=gN_{\bfb}g^{-1}$ the conditions are equivalent to
\begin{equation}\label{E13.1} N_{\hbf{b,1}}\cap g^{-1}\gG g
  =N_{\hbf{b,1}}\cap \gG\quad \hbox{(respectively } N_{\bfb}\cap g^{-1}\gG
  g = N_{\bfb} \cap \gG).
  \end{equation}
  This in turn means that $g^{-1}\gG g$ is integral on $N_{\hbf{b,1}}$
  (respectively integral on $N_{\bfb}$), in other words, that for some
  rational representation $\gr:G\lra GL(V)$ we have
  $\gr_{|N_{\hbfs{b,1}}}(N_{\hbf{b,1}}\cap g^{-1}\gG g)\inn GL(V_{\fZ})$
  (respectively $\gr_{|N_{\hbfs{b}}}(N_{\bfb}\cap g^{-1}\gG g)\inn
  GL(V_{\fZ})$).  Note that this definition depends on the choosen maximal
  torus, as well as on $\gG$. If $S'=xSx^{-1}$ is another maximal
  $\fQ$-split torus, then $N_{\hbf{b}}'=xN_{\bfb}x^{-1}$ is the standard
  symmetric subgroup with respect to $S'$. If $N\inn G$ is $\gG$-integral
  with respect to $N_{\bfb}$ (i.e., there is $g\in G(\fQ)$ such that $g
  N_{\bfb}g^{-1}=N$ and $N_{\hbf{b,1}}\cap g^{-1}\gG g = N_{\hbf{b,1}}\cap
  \gG$), then $N_{\hbf{b,1}}'\cap
  (gx^{-1})^{-1}\gG(gx^{-1})=N_{\hbf{b,1}}'\cap x\gG x^{-1}$; in other
  words when $N$ is $\gG$-integral with respect to $N_{\bfb}$, then $N$ is
  $x\gG x^{-1}$-integral with respect to $N_{\bfb}'$ (with similar
  statements for strongly $\gG$-integral).  At least in the classical
  cases, when $G$ is a matrix group, there is a very canonical choice for
  $N_{\bfb}$, namely as a subgroup consisting of block matrices, so this
  dependence is not unreasonable.

  Let us now suppose $G$ is classical, $\gr':G'\lra GL_D(V)$ the
  fundamental representation, $\gr:G\lra Res_{k|\fQ}(GL_D(V))$ the
  corresponding representation of $G$. We have $P_{\bfb}=Res_{k|\fQ}P_b',\
  N_{\bfb}=Res_{k|\fQ}N_b'$ and $P_b'$ (resp.~$N_b'$) are given in terms of
  $(V,h)$ by (\ref{E8.b}) (resp.~by (\ref{E12.1})). Clearly if
  $N=N_{\bfb}^g$ and $N_{\bfb}=\cN_G(W)$, then $N=\cN_G(g(W))$.
\begin{lemma}\label{L13.1} $N=\cN_G(g(W))$ is $\gG$-integral $\iff$
  $W_{\fZ}=W\cap V_{\fZ}=W\cap \gr(g^{-1})(V_{\fZ})$.
\end{lemma}
{\bf Proof:} By definition $N_1\cap \gG=gN_{\hbf{b,1}}g^{-1}\cap \gG=
g(N_{\hbf{b,1}}\cap \gG)g^{-1}$, and this holds if and only if
$N_{\hbf{b,1}}\cap g^{-1}\gG g = N_{\hbf{b,1}}\cap \gG$, i.e., $g^{-1}\gG
g$ meets $N_{\hbf{b,1}}$ in the arithmetic group $\gG$. But this holds
precisely when $g^{-1}\gG g$ maps $W_{\fZ}=W\cap V_{\fZ}$ into itself, and
this is equivalent to $W_{\fZ}=W\cap \gr(g^{-1})(V_{\fZ})$, as
$\gr(g^{-1}\gG g)$ maps $\gr(g^{-1})(V_{\fZ})$ into itself, and this is the
statement of the lemma. \ende Recall also
\begin{definition}\label{D13.2} A lattice $V_{\fZ}\inn V_{\fQ}$ being given,
  a submodule $W_{\fZ}$ is {\it pure}, if $n\cdot x\in W_{\fZ},\ n\in \fZ\
  \Ra x\in W_{\fZ}$.
\end{definition}
\begin{lemma}\label{L13.2} There is a 1-1 correspondence between rational
  subspaces $W_{\fQ}\inn V_{\fQ}$ and pure $\fZ$-submodules $W_{\fZ}\inn
  V_{\fZ}$, given by
  \[ W_{\fQ}\mapsto W_{\fQ}\cap V_{\fZ},\quad
  W_{\fZ}\mapsto W_{\fZ}\otimes_{\fZ}\fQ.\]
\end{lemma}
{\bf Proof:} Clear. \ende Note that the statement of Lemma \ref{L13.1} is
also equivalent to $g(W_{\fZ})=g(W)\cap V_{\fZ}$, and the latter is by
Lemma \ref{L13.2}, a pure submodule. This then yields:
\begin{corollary}\label{C14.1} There is a 1-1 correspondence between the
  set of $\gG$-integral symmetric subgroups and pure submodules of the form
  $g(W_{\fZ})$, with $g\in G(\fQ)$ and $W_{\fZ}\inn V_{\fZ}$ the submodule
  above (cf.~Lemma \ref{L13.1}).
\end{corollary}
{\bf Proof:} As we just remarked, $N$ is $\gG$-integral $\iff$ $g(W_{\fZ})$
is pure, and for any $g(W_{\fZ}),\ g\in G(\fQ)$, which is pure,
$N_{\bfb}^g$ is clearly $\gG$-integral.  \ende Next we consider, for
$d=\dim(W)$, the representation
\[ R=\bigwedge^d\gr:G\lra GL(\cV),\quad \cV=\bigwedge^dV.\]
Since $\gr(N_{\bfb})=\cN_G(W)$, it follows that
$R(N_{\bfb})=\bigwedge^d\gr(N_{\bfb}) =
\cN_G(\bigwedge^dW)=\cN_G(W_{\bfb})$, where $W_{\bfb}$ is one-dimensional
in $\cV$, defined over $\fQ$, and we have slightly abused notation by
denoting by $\cN_G$ the inverse image under $R$ of the corresponding
normalizer in $GL(\cV)$. Our lattice $V_{\fZ}$ produces a lattice in $\cV,\
\cV_{\fZ}=\bigwedge^dV_{\fZ}$, and $\gG$ is commensurable with
$R^{-1}(GL(\cV_{\fZ}))$.

We now return to the general situation; $G$ is $\fQ$-simple of hermitian
type, $P_{\bfb}$ is a standard parabolic and $N_{\bfb}$ is an incident
symmetric subgroup, which we take to be the standard one. Since $N_{\bfb}$
is reductive, by \cite{BHC}, Theorem 3.8, there exists a rational
representation $\pi:G\lra GL(\cV)$, defined over $\fQ$, and an element
$v\in \cV_{\fQ}$, such that $v\cdot \pi(G)$ is a closed orbit and
$N_{\bfb}=\pi^{-1}(\cN_{GL(\cV)}(v)$. For example, in the classical cases,
the representation $R$ above is such a $\pi$. We now assume that $\cV$ is
given a $\fZ$-structure $\cV_{\fZ}$ such that $\gG$ is given by
\begin{equation} \label{E5.1} \gG=\pi^{-1}(GL(\cV_{\fZ})).
\end{equation}
Let $W_{\bfb}=\fQ\langle v\rangle$ be the one-dimensional vector subspace
spanned by $v$; then we may choose a primitive integral vector $w\in
W_{\bfb}$ such that $N_{\bfb}=\pi^{-1}(\cN_{GL(\cV)}(w))$. We consider the
orbit $w\cdot \pi(G)$; as is well-known there is a natural isomorphism
$w\cdot \pi(G)\cong G/N_{\bfb}$ given by $w\cdot \pi(g) \mapsto gN_{\bfb}$.
We may consider the lattice $\cV_{\fZ}$, defining {\it integral points}
$\cV_{\fZ}\cap w\cdot \pi(G) \stackrel{\sim}{\lra} \cV_{\fZ}\cap
G/N_{\bfb}$.
\begin{lemma}\label{L14.1} Assume $\gG$ fulfills (\ref{E5.1}).
  A subgroup $N$ given by $N=N_{\bfb}^g$ is $\gG$-integral $\iff$ under the
  isomorphism $w\cdot \pi(G) \stackrel{\sim}{\lra} G/N_{\bfb}$, $N$ is
  given by an integral point $gN_{\bfb}$, i.e., $w\cdot \pi(g) \in
  \cV_{\fZ}$.
\end{lemma}
{\bf Proof:} We have the following equivalences:
\begin{eqnarray*}  & & N_1\cap \gG = g(N_{\hbf{b,1}}\cap \gG)g^{-1} \\
  & \stackrel{(\ref{E5.1})}{\iff} & N_1 \cap \pi^{-1}(GL(\cV_{\fZ})) =
  g(N_{\hbf{b,1}}\cap \pi^{-1}(GL(\cV_{\fZ})))g^{-1} \\ &
  \stackrel{\hbox{\footnotesize apply $\pi$}}{\iff} & \pi(N_1)\cap
  GL(\cV_{\fZ}) = \pi(g)(\pi(N_{\hbf{b,1}})\cap GL(\cV_{\fZ}))\pi(g^{-1})
  \\ & \iff & \pi(N_{\hbf{b,1}})\cap \pi(g^{-1})GL(\cV_{\fZ})\pi(g) =
  \pi(N_{\hbf{b,1}})\cap GL(\cV_{\fZ}) \\ &
  \stackrel{\parbox{2cm}{\footnotesize$\pi(N_{\hbfs{b,1}})=$\\
      $\cN_{\pi(G)}(W_{\hbfs{b}})/ \cZ_{\pi(G)}(W_{\hbfs{b}})$ }} {\iff} &
  \quad\quad\quad\quad\cN_{\pi(G)}(W_{\bfb})/\cZ_{\pi(G)}(W_{\bfb})\cap
  \pi(g^{-1})GL(\cV_{\fZ})\pi(g) \\ & & \quad \quad \quad \quad \quad \quad
  \quad = \cN_{\pi(G)}(W_{\bfb})/\cZ_{\pi(G)}(W_{\bfb}) \cap GL(\cV_{\fZ})
  \\ & \iff & \pi(g)(W_{\bfb}\cap \cV_{\fZ})\inn \cV_{\fZ} \\ & \iff &
  w\cdot \pi(g)\in \cV_{\fZ}
\end{eqnarray*}
where the last equivalence follows from the fact that $w$ is primitive. The
Lemma follows. \ende
\begin{corollary}\label{C5.1} The set of $\gG$-integral symmetric subgroups is
  the set of subgroups corresponding to the integral points,
\[ \left\{\parbox{5.3cm}{$\gG$-integral symmetric subgroups \\ conjugate to
    $N_{\bfb}$ }\right\} \cong G/N_{\bfb}\cap \cV_{\fZ}.\]
\end{corollary}
{\bf Proof:} This follows immediately from the proceeding Lemma. \ende
Utilizing Corollary \ref{C5.1}, we can prove finiteness of the set of
$\gG'$-equivalence classes of $\gG$-integral symmetric subgroups, for any
arithmetic subgroup $\gG'\inn G(\fQ)$. Recall the basic finiteness result
of \cite{BHC}.
\begin{theorem}[\cite{BHC}, 6.9]\label{t15a.1}
  Let $G$ be a reductive algebraic group defined over $\fQ$, $\pi:G\lra
  GL(V)$ a rational representation defined over $\fQ$, $\cL\inn V$ a
  lattice in $V_{\fQ}$ invariant under $G_{\fZ}$, and $X$ a closed orbit of
  $G$. Then $X\cap \cL$ consists of a finite number of orbits of $G_{\fZ}$.
\end{theorem}
Here $G\inn GL(n,\fC)$ and $G_{\fZ}=G\cap M_n(\fZ)$.
\begin{corollary}\label{C6.1} Given $G,\ S,\ P_{\bfb},\ N_{\bfb}$ and $\gG$
  as above ($\gG$ as in (\ref{E5.1})), there are finitely many
  $\gG'$-equivalence classes of $\gG$-integral symmetric subgroups, for any
  arithmetic subgroup $\gG'\inn G(\fQ)$.
\end{corollary}
{\bf Proof:} Since $\gG$ satisfies (\ref{E5.1}), Corollary \ref{C5.1} holds
and Theorem \ref{t15a.1} may be applied to $\gG$, hence finiteness holds
for any $\gG'$. \ende

Note that under the action of $\gG$ on $\gr(G)\cdot v\cap V_{\fZ}$, all
orbits are bijective to $\gG/(\gG\cap N_{\bfb})$. Let the orbit
decomposition with respect to $\gG'$ be
\[ \gG'\bs \gr(G)\cdot v \cap V_{\fZ}=\cO_1\cup \cdots \cup \cO_q.\]
Choose, in each orbit $\cO_i$, a representative $x_i$, and let $N_{x_i}$ be
the corresponding integral symmetric subgroup. The set $\{N_{x_i}\}$ serves
as a finite set of $\gG$-integral symmetric subgroups representing all
$\gG'$-equivalence classes of such. The following is then well defined.
\begin{definition}\label{D6.1} Given $G,\ N_{\bfb},\ \gG$ as above,
  $\gG'\inn G(\fQ)$ arithmetic, the {\it class number} of
  $\gG'$-equivalence classes of $\gG$-integral symmetric subgroups is the
  cardinality
\[ \mu(G,N_{\bfb},\gG,\gG'):= |\gG'\bs(G/N_{\bfb}\cap \cV_{\fZ})|.\]
If in a discussion $G$ and a maximal $\fQ$-split torus $S\inn G$ are fixed,
then $N_{\bfb}$ depends only on the integer $b\in \{1,\ldots, s\}$
($s=\rank_{\fQ}G$)\footnote{here again with the two exceptions for $b=s=t$
  and the two exceptional domains where there are three, resp.~two
  isomorphism classes of $N_{\bfb}$}, and we will denote this class number
by $\mu_b(\gG,\gG')$.
\end{definition}

\section{Arithmetic quotients}
In this section we keep the above notations. $G$ is $\fQ$-simple of
hermitian type, $\cD=G(\fR)/K=G(\fR)^0/K^0$ the hermitian symmetric space,
$\gG\inn G(\fQ)$ an arithmetic subgroup. The group $\gG$ acts on $\cD$ by
means of holomorphic isometries, preserving the natural Bergmann metric.
\begin{definition}\label{d83.1} The quotient $X_{\gG}:=\gG\bs \cD$, where
  $\gG\inn G(\fR)$ is arithmetic, is called an {\em arithmetic quotient}.
\end{definition}
If $\gG$ acts without fix points, then the quotient $X_{\gG}$ is a smooth
complex manifold, not compact in general. If $\gG$ has fix points, then
$X_{\gG}$ has certain quotient singularities, which can be described as
follows. Let $\gG_1\inn \gG$ be a normal subgroup of finite order without
elements of finite order, so that $\gG_1$ acts freely and hence $X_{\gG_1}$
is smooth. We have a Galois cover,
\begin{equation}\label{e83.1} X_{\gG_1}\lra X_{\gG},
\end{equation}
with $X_{\gG_1}$ smooth, and Galois group acting, yielding the
singularities of $X_{\gG}$. It is clear that the Galois group, $\gG/\gG_1$,
creates the singularities, so they are controlled by certain properties of
$\gG/\gG_1$, such as the orders of the elements, etc. In particular,
$X_{\gG}$ is {\em still} smooth if $\gG/\gG_1$ is generated by reflections,
as the quotient is then smooth by Chevally's theorem.  It is well-known
that $X_{\gG}$ is compact $\iff$ $G$ is anisotropic. Suppose this is the
case, and that in addition $\gG$ has no elements of finite order.  Then, as
Kodaira showed in 1954 as one application of his embedding theorem,
$X_{\gG}$ is a smooth projective variety, the canonical bundle
$K_{X_{\gG}}$ being ample. In this case one has Hirzebruch proportionality,
which states that the ratios of the Chern numbers of $X_{\gG}$ are equal to
the ratios of the corresponding Chern numbers of the compact hermitian
symmetric spaces $\check{\cD}$, and the overall factor of proportionality
is just the volume of $X_{\gG}$, which is the same as the volume in $\cD$
of a fundamental domain of $\gG$, where volume is taken with respect to the
Bergmann metric.

\subsection{Satake compactification and Baily-Borel embedding}
In case $G$ is not anisotropic, $X_{\gG}$ is not compact. It has a
topological compactification $X_{\gG}^*$, the so-called {\em Satake
compactification}. This is constructed by putting an appropriate topology,
the Satake topology, on $\cD^*:=\cD\cup \{\hbox{rational boundary
components}\}$ (\cite{BB}, 4.8). With the Satake topology, the action of
$\gG$ on $\cD$ extends to one on $\cD^*$ (\cite{BB}, 4.9), and the quotient
$\gG\bs \cD^*=\xgs$ is the sought for compactification. It has the following
property:
\begin{proposition}[\cite{BB}, 4.11]\label{p84.1}
$\xgs$ is a compact, Hausdorff space, and the complement $\xgs \bs \xg$ is a
finite disjoint union $\xgs\bs\xg = V_1\cup\cdots \cup V_N$, with each $V_i$
an arithmetic quotient of dimension and $\fQ$-rank less than that of $\xg$.
The length $k$ of a maximal chain $V_{i_1}\subsetneq V_{i_2}^*\subsetneq
\cdots \subsetneq
V_{i_k}^*$ is the $\fQ$-rank of $G$.
\end{proposition}
In our discussion of the $\fQ$-hermitian symmetric subgroups in section
\ref{classification} we determined the rational
boundary components for each $G$
giving rise to quotients $V_i$. One case of particular interest are the
hyperbolic planes, discussed in detail in \cite{hyp}.
We know that in this case all rational boundary
components are zero-dimensional, i.e., points. Hence the finite union of
\ref{p84.1} is a disjoint union of points; the number of such is just the
number of cusps, defined as follows. Suppose again we have the fixed
$\fQ$-split torus $S$ and the standard subgroups $P_{\bfb}$ and $N_{\bfb}$
with respect to $S$.
\begin{definition}\label{D7.1} \begin{itemize}\item[(i)] For $b\in
    \{1,\ldots, s\}$, the number of boundary varieties, conjugate to the
    $b^{th}$ standard one, is the cardinality
\[ \nu_b(\gG)=|\gG\bs G(\fQ)/P_{\bfb}(\fQ)|.\]
\item[(ii)] The {\it number of $\gG$-cusps} is the cardinality (where $B$
  is a Borel subgroup)
\[ h(\gG)=|\gG\bs G(\fQ)/B(\fQ)|.\]
\end{itemize}
\end{definition}
Note that $h(\gG)$ is also the number of {\it maximal flags} of boundary
varieties, and it is often given by a class number, hence the notation.
Since in the case of
hyperbolic planes ($s=1$), $\nu_1(\gG)=h(\gG)$, both of these are given by
the results of \cite{hyp} in terms of class numbers of certain fields.
More generally, the number
of components $N$ occurring in Proposition
\ref{p84.1} is a sum $N=r_1+\ldots
+ r_s$, where $s=\fQ$-rank of $G$, $r_b$ = \# equivalence classes of
boundary components conjugate to $F_{\hbf{b}}$. Then
\begin{proposition}\label{p84.3} For any $\xg$, the number $N$ of Proposition
\ref{p84.1} can be expressed: $N=r_1+\ldots + r_s,$ and $r_b=\nu_b(\gG)$ as
in Definition \ref{D7.1}.
\end{proposition}

The term Baily-Borel embedding of $\xgs$ refers to the following result.
\begin{theorem}[\cite{BB}, 10.11, 10.12]\label{l84.1}
$\xgs$ can be embedded in projective space as a normal algebraic variety $V$.
If $G$ has no normal $\fQ$-subgroups of dimension three, then the field of
rational functions $K(V)$ is canonically isomorphic with the field of
automorphic functions for $\gG$.
\end{theorem}
It follows in particular that $\xg$ is a normal, quasi-projective variety,
which is even smooth if $\gG$ is torsion free.

\subsection{Toroidal embeddings}
Recall the decomposition (of algebraic groups over $\fR$)
$P_{\bfb}=M_{\bfb}L_{\bfb}\cR_{\bfb}\sdprod
\cU_{\bfb}$ for the real parabolic, with the exact sequence
\[1\lra M_{\bfb}L_{\bfb}\cR_{\bfb} \lra P_{\bfb} \lra \cU_{\bfb} \lra 1;\]
this gives rise to a similar sequence for $\gG_{\bfb}=\gG\cap P_{\bfb}(\fR)$,
\[ 1\lra \gG_{\bfb}^{\ell} \lra \gG_{\bfb} \lra \gG_{\bfb}^r \lra 1,\]
with $\gG_{\bfb}^{\ell}$ being the intersection with the Levi factor and
$\gG_{\bfb}^r$ the intersection with the radical of $P_{\bfb}$. Recall
further that $\cU_{\bfb}=\cZ_{\bfb}V_{\bfb}$, where $\cZ_{\bfb}$ is the
center, and $M_{\bfb}L_{\bfb}$ acts trivially on $\cZ_{\bfb}$ and by means
of a symplectic representation on $V_{\bfb}$, while $\cR_{\bfb}$ acts
transitively on $\cZ_{\bfb}$ defining a homogenous self dual cone
$C_{\bfb}\inn \cZ_{\bfb}$, and $\cR_{\bfb}$ acts on $V_{\bfb}$ by means of
complex linear transformations, see Theorem 1.1.
This then gives us the following results
about the factors of $\gG_{\bfb}$ (and similar results hold for any
$\gG_F=N(F)\cap \gG$):
\begin{itemize}\item[i)] $M_{\bfb}(\fR)$ is compact, hence $\gG\cap
  M_{\bfb}(\fR)$ is
  finite. In particular, if $\gG$ has no torsion, $\gG\cap M_{\bfb}(\fR)=e$.
\item[ii)] $\gG\cap L_{\bfb}(\fR)$ is an arithmetic subgroup of
  $L_{\bfb}(\fR)$, and
  $\gG\cap L_{\bfb}(\fR)^0=:\gG_L^0$ acts on the boundary component $F_{\bfb}$,
  with the boundary variety $W_{\bfb}=\gG_L^0\bs F_{\bfb}\inn X_{\gG}^*$.
\item[iii)] Let $V_{\fZ}=\gG\cap V_{\bfb}(\fR)$.
  Then the group $\gG_L^0\sdprod
  V_{\fZ}$ acts on $F_{\bfb}\times V_{\bfb}(\fR)$
  (recall that $V_{\bfb}(\fR)$ has the
  structure of complex vector space), and if $\gG$ is torsion free,
  the quotient is an analytic family of abelian
  varieties over the arithmetic quotient $W_{\bfb}=\gG_{L}^0 \bs F_{\bfb}$.
\item[iv)] {\bf ([SC], p.~248)} There is an exact sequence
$$1 \lra \gG^{\p} \lra \gG_{\bfb} \lra \gG^{\p\p} \lra 1, $$
\begin{itemize}\item{} $\gG^{\p}$= subgroup of elements in $\gG_{\bfb}$
acting
trivially by conjugation on $Lie(\cZ_{\bfb})$,
\item{} $\gG^{\p\p}$= group of automorphisms of $Lie(\cZ_{\bfb})$
induced by $\gG_{\bfb}$; these map $C_{\bfb}$ into itself.
\end{itemize}
\end{itemize}
The fourth point is important for the compactification theory, as one lets
first $\gG''$ act, then $\gG'$.
A sketch of the construction is as follows: fix a boundary component $F$,
rational with respect to \gG\ (i.e., $\gG\cap N(F)$ is a lattice).
Let ${\bf E}_F \inn \cZ(F)_{\fC}\times V(F) \times F$ be the realisation of
$\cD$ as a Siegel domain as in \cite{SC}, and let $1\lra \gG'\lra \gG_F\lra
\gG'' \lra 1$ be the sequence above for $\gG_F=N(F)\cap \gG$. Furthermore,
the objects denoted above by a subscript $?_{\bfb}$ will be denoted here by
$?(F)$, for example $\cZ(F)$ instead of $\cZ_{\bfb}$, $C(F)$ instead of
$C_{\bfb}$, etc.
\begin{proposition}[\cite{SC}, p.249]\label{toroidal}
A partial compactification along $F$ can be constructed as follows:
\begin{itemize}
\item[1)] Let $\cZ(F)_{\fZ}$ act on $\cZ(F)_{\fC}$ defining the
algebraic torus $T_F$; do this in the fibration
$${\bf E}_{\{1\}}={\bf E}_F/\cZ(F)_{\fZ} \inn
\cZ(F)_{\fC}/\cZ(F)_{\fZ}\times V(F)\times F\lra F.$$
More precisely, the map $\cZ(F)_{\fC}\lra \cZ(F)_{\fC}/\cZ(F)_{\fZ}$
is given by $\exp(2\pi i\gl_1),\ldots,\exp(2\pi i\gl_k)$, where one
chooses a $\fZ$-base $\xi_1,\ldots,\xi_k$ of $\cZ(F)_{\fZ}$,
and $\gl_i:\cZ(F)_{\fC}\lra \fC$ is the dual base.
\item[2)] Now compactify the algebraic torus $T_F$ by $T_F\inn
T_{F{\{\gs_{\ga}\}}},\  \{\gs_{\ga}\}$ a $\gG^{\p\p}$-admissible polyhedral
decomposition of $C(F)\inn \cZ(F)$ ($\gG''$ as in iv) above).
$T_{F{\{\gs_{\ga}\}}}$ is locally of
finite type, but will have infinitely many components corresponding to
integral vectors $v\in \cZ(F)_{\fZ}\cap C(F)$. The cones $\gs_{\ga}$
themselves correspond to orbits of highest codimension, i.e., to points.
If $\gs_{\ga}\cap \cZ(F)_{\fZ}$ is spanned by $v_1,\ldots, v_k$, then
$\gs_{\ga}$ corresponds to $\gD_1\cap \cdots \cap \gD_k$, where $\gD_j$
is the divisor corresponding to $v_j$.
\item[3)] Glue these into ${\bf E}_{\{1\}}$ by forming the fibre product
$$({\bf E}_{\{1\}})\times^{T_F}(T_{F{\{\gs_{\ga}\}}})$$
and setting $({\bf E}_{\{1\}})_{\{\gs_{\ga}\}}$= interior of the closure of
${\bf E}_{\{1\}}$ in $({\bf E}_{\{1\}})\times^{T_F}(T_{F{\{\gs_{\ga}\}}})$.
Hence  $({\bf E}_{\{1\}})_{\{\gs_{\ga}\}}$ has a fibre structure over
$F \times V(F)$ with fibres $T_{F{\{\gs_{\ga}\}}}$. If $\gD_i$ is the divisor
corresponding to $\xi_i$ as in 1), then $\gD_i=\{z_i=0\},$ where
$(z_1,\ldots,z_k)$ are local coordinates on $T_F$.
\item[4)] $\gG_F$ still acts on  $({\bf E}_{\{1\}})_{\{\gs_{\ga}\}}$, as
follows. $\gG''$ now acts freely on  $({\bf E}_{\{1\}})_{\{\gs_{\ga}\}}$,
giving a fibre space over $F\times V(F)$ with fibre a finite $T_F$
compactification, i.e., modulo $\gG''$ there are only finitely many
integral vectors, hence components, in the fibre. Now $\gG'$ acts on
$F\times V(F)$; as $\cZ_{\bfb}$ acts trivially this amounts to an action of
$\gG_L^0\sdprod V_{\fZ}$ as in iii) above, and this action
extends to $\gG''\bs ({\bf E}_{\{1\}})_{\{\gs_{\ga}\}}$, hence an
open neighborhood of $F$
will give an open neighborhood of the boundary variety $W(F)=\gG_L^0\bs F$
in $\ove{\xg}$; this is the sought for partial compactification.
\end{itemize}
\end{proposition}
Next one glues these partial compactifications together by means of
$\{ \gs_{\ga,F}\}$, a \gG-admissible collection of polyhedral cones, one such
collection for each boundary component. The main result is:
\begin{theorem}[\cite{SC}, Main Theorem 1, p.~252]
With \gG, \cD\ as above, for every
\gG-admissible collection of polyhedral cones $\{ \gs_{\ga,F}\}$, there is
a unique compactification $\ove{\xg}\ (=(\ove{\xg})_{\{\gs_{\ga}\}})$
which is locally given at each $F$ (more precisely at $W(F)$) by the partial
compactification above, corresponding to the given collection of cones.
$\ove{\xg}$ is a compact Hausdorf, analytic variety, which is an algebraic
space. Furthermore, for properly chosen \gG-admissible collections of
polyhedral cones, the compactification is 1) a projective
resolution of the Satake
compactification: $\ove{\xg}\lra \xgs$, hence a projective variety, and 2)
smooth with $\gD_{\gG}:=\xgc-\xg$ a normal crossings divisor.
\end{theorem}

\section{Modular subvarieties}
In this paragraph, the data $G,\ S,\ \cD$ will be fixed as above, so that
for each $b=1,\ldots, s=\rank_{\fQ}G$ we have the standard boundary
components $F_{\bfb}$, the standard parabolic $P_{\bfb}$ and the standard
incident symmetric subgroup $N_{\bfb}$.

\subsection{Baily-Borel compactification}
Let $N\inn G$ be a reductive subgroup of hermitian type (this implies in
particular that $N$ is defined over $\fR$, and we assume the inclusion
$N\inn G$ is also), $\cD_N\inn \cD$ the subdomain (holomorphic symmetric
embedding) determined by $N$.
\begin{definition}\label{D9.1} The subdomain $\cD_N\inn \cD$ will be said
  to be {\it defined over $\fQ$}, if $N$ is a $\fQ$-subgroup of $G$.
\end{definition}
Suppose a subdomain $\cD_N\inn \cD$ is defined over $\fQ$, and consider an
arithmetic subgroup $\gG\inn G(\fQ)$ and the arithmetic quotient
$X_{\gG}=\gG\bs \cD$. Note that for a reductive subgroup of hermitian type
$N\inn G$, the intersection $\gG_N:= \gG\cap N$ will be an arithmetic
subgroup if and only if $N$ is defined over $\fQ$, and this is the case if
and only if $\cD_N$ is defined over $\fQ$. Hence the arithmetic quotient
$X_{\gG_N}:=\gG_N\bs \cD_N$ is defined, and clearly fits into a commutative
square
\begin{equation}\label{E9.1}
  \begin{array}{ccc} \cD_N & \hra & \cD \\ \downarrow & &
    \downarrow \\ X_{\gG_N} & \hra & X_{\gG}.
  \end{array}
\end{equation}
\begin{definition}\label{D9.2a} A {\it modular subvariety} on $X_{\gG}$ is a
  sub-arithmetic quotient $X_{\gG_N}$ as in (\ref{E9.1}), where $\cD_N$ is
  defined over $\fQ$. A modular subvariety $X_{\gG_N}\inn X_{\gG}$ will be
  called {\it rational} (resp.~{\it $\gG$-integral}), if $N$ is a rational
  (resp.~$\gG$-integral) symmetric subgroup as in Definition \ref{D4.1}
  (resp.~\ref{D5.1}).
\end{definition}
The embedding of (\ref{E9.1}) turns out to extend to one of the Baily-Borel
embeddings, legitimizing the terminology sub{\it variety}. This is given by
the following result of Satake.
\begin{theorem}\label{t17.1} Let $X_{\gG_N}^*\inn {\Bbb P}^N,\
  X_{\gG}^*\inn
  {\Bbb P}^{N'}$ be Baily-Borel embeddings. Then there is a linear
  injection ${\Bbb P}^N\hra {\Bbb P}^{N'}$ making the diagram
  $$\begin{array}{ccl} X_{\gG_N}^* & \hra & {\Bbb P}^N \\ \cap & & \cap \\
    X_{\gG}^* & \hra & {\Bbb P}^{N'}\end{array}$$ commute and
  making $X_{\gG_N}^*\inn X_{\gG}^*$ an algebraic subvariety.
\end{theorem}
{\bf Proof:} We have an injective holomorphic embedding $\cD_N\hra \cD$
which comes from a ${\Bbb Q}$-morphism $\rho:(N)_{{\Bbb C}}\hra
(G)_{{\Bbb C}}$ such that $\rho(\Gamma_{N})
    \inn \Gamma$. Hence we map apply \cite{S2}, Theorem 3, and
  the theorem follows from this. \ende

\begin{definition}\label{D9.2} We say that $X_{\gG_N}$ and a boundary
  variety $W_i$ are {\it incident}, if in $\cD^*$ there is rational
  boundary component $F$ with parabolic $P=N(F)$ covering $W_i$, such that
  $N$ and the corresponding parabolic $P$
  are incident.
\end{definition}
Note the following
\begin{lemma}\label{L9.2} $X_{\gG_N}$ and $W_i$ are incident, if and only
  if $W_i^*\inn X_{\gG_N}^*$ is a maximal-dimensional boundary component of
  $X_{\gG_N}^*$ (if $\dim(W_i)>0$), resp.~if and only if $W_i\inn
  X_{\gG_N}^*$ (if $\dim(W_i)=0$).
\end{lemma}
{\bf Proof:} If $\dim(W_i)>0$, then the groups $P$ and $N$ are incident if
$F\inn \cD_N^*$ and $F$ is maximal with this property, and if $F\in
\cD_N^*$ is rational and maximal with this property, then $P$ and $N$ are
incident. If $\dim(W_i)=0$ and $N$ is an incident symmetric subgroup, then
$W_i\inn \cD_N^*$ is a (point) rational boundary component, and
conversely. \ende

We now consider $\gG$-integral symmetric subgroups $N$ and arbitrary
arithmetic subgroups $\gG'\inn G(\fQ)$, let $\gG'_N=N\cap
\gG'$ and consider the corresponding integral modular subvarieties they define,
$X_{\gG'_N}\inn  X_{\gG'}$. As described above, the inclusion
extends to the Baily-Borel embeddings $X_{\gG'_N}^*\inn X_{\gG'}^*$. We now
take $\gG$ to be $G_{\fZ}$ for some rational representation $\gr:G\lra
GL(V_{\fZ})$, that is $\gG=\gr^{-1}(GL(V_{\fZ}))$ for some $\fZ$-structure
$V_{\fZ}$ on $V$.
Recall the notations $\nu_b(\gG'),\
b=1,\ldots, s$ and $\mu_b(\gG,\gG'),\ b=1,\ldots, s$ of Definition \ref{D6.1}
and \ref{D7.1}, respectively, for the number of $b^{th}$ boundary varieties
and the number of $b^{th}$ integral modular subvarieties, respectively. We
let $W_{b,i},\ b=1,\ldots, s,\ i=1,\ldots,\nu_b(\gG')$ be the corresponding
boundary varieties on the Satake compactification, $Y_{b,j},\ b=1,\ldots,
s,\ j=1,\ldots, \mu_b(\gG,\gG')$ the corresponding $\gG$-integral modular
varieties, everything on the arithmetic quotient $X_{\gG'}$. Then
the main result of the paper is the following.
\begin{theorem} Let $\gG$ be as above, $\gG'\inn G(\fQ)$ arithmetic, and
  $X_{\gG'}\inn X_{\gG'}^*$ the Satake compactification,
  $X_{\gG'}^*-X_{\gG'}=\sum_{b,i}W_{b,i}$. Then $\Xi:=\sum_{b,j}Y_{b,j}$ is
  a complete (finite, non-empty) set of $\gG'$-equivalence classes of
  $\gG$-integral modular subvarieties, such that for each $W_{b,i}$, there
  is at least one $Y_{b,j}$ incident to $W_{b,i}$.
\end{theorem}
{\bf Proof:} There is for each $W_{b,i}$ an incident $\gG$-integral modular
subvariety because for any representative parabolic there is a
$\gG$-integral symmetric subgroup which is incident.
The finiteness result Corollary \ref{C6.1} implies
that for each $N_{\bfb}$ (of which there are finitely many) there are
finitely many $\gG'$-equivalence classes of $\gG$-integral symmetric
subgroups of $G$ conjugate to $N_{\bfb}$, so a complete set of
$\gG'$-equivalence classes is
finite. \ende

This gives us a {\it well-defined, non-empty, finite} set of subvarieties
of the Baily-Borel embedding $X_{\gG'}^*\inn \fP^N$ for any subgroup
$\gG'\inn \gG$ of finite index. Furthermore these have a prescribed
behavior near the cusps. For example, if $f:\cD\lra \fC$ is a modular form
whose zero divisor $D_f$ on $X_{\gG'}^*$ contains the union of the integral
modular subvarieties, then $f$ is a cusp form for $\gG'$.

\subsection{Incidence}
Consider a toroidal compactification $\xgc$ which is smooth and projective;
consider what incidence means here.
Let $W_i$ be a rational boundary variety, and let $P$ be a parabolic
$P=\cN_{G(\fR)}(F)$ for some rational boundary component $F$ which covers
$W_i$. We have the decomposition $P=(ML\cR)\sdprod \cZ V$ of the parabolic.
Recall from the construction
\ref{toroidal} that the inverse image $\pi^{-1}(W_i)$ in $\xgc$ of $W_i\inn
\xgs$ is a divisor which is a torus embedding bundle over the family of
abelian varieties ``$V/V_{\fZ}$ over $W_i$''. On the other hand, if $X_{\gG_N}$
is an integral modular subvariety incident to $W_i$, and if $\dim(W_i)>0$,
then the proper transform of $X_{\gG_N}^*$ on $\xgc$ will meet
$\pi^{-1}(W_i)$ in a {\it section} of the family of abelian varieties over
$W_i$. In a sense, the standard one will meet in the zero-section, the
others meet in certain sections associated with level structures
(e.g. sections of torsion points). For $\dim(W_i)=0$ the situation is
slightly different. We now discuss this in more detail.

We consider first the case that $\dim(W_i)>0$. Then by (\ref{e10.1}),
$N_{\bfb}$ (hence any $G(\fQ)$-conjugate) has the form $N_1\times N_2$,
where $N_1\inn P_{\bfb}$ is a hermitian Levi factor. Considering the
arithmetic group $\gG_{N_{\bfb}}$ acting on $\cD_{N_1}\times \cD_{N_2}$,
since the product is defined over $\fQ$, the quotient $\gG_{N_{\bfb}}\bs
\cD_N$ is at most a finite quotient of a product itself. We assume that in
fact $X_{\gG_{N}}$ is a product (we will show below in Lemma \ref{L6.4.2}
that for $N$ $\gG$-integral this always holds);
then $X_{\gG_N}=X_1\times X_2$, where
$X_1$ is the arithmetic quotient $\gG_1\bs \cD_{N_1}$ and this is
isomorphic to the boundary variety $W_i$. It follows that $X_2=\gG_1\bs
\cD_{N_2}$ has rational boundary components which are zero-dimensional, say
$w\in X_2^*$, such that with respect to the natural inclusion
$i:X_{\gG_N}^*\inn X_{\gG}^*$ we have
\begin{equation}\label{E35} i(X_1\times \{w\})=W_i.
\end{equation}
Recall further that any two hermitian Levi factors are conjugate by an
element $g\in V$, and that, modulo $\gG$, this means a point of the abelian
variety ``$V/V_{\fZ}$''. This is of course true for any point $t\in W_i$, so we
get
\begin{lemma}\label{L11.1} Given $W_i$, $X_{\gG_N}^*$ any integral modular
  subvariety incident to it. Then the proper transform of $X_{\gG_N}^*$ in
  $\overline{X}_{\gG}$ determines a {\it section} of the
  family of abelian varieties of $\pi^{-1}(W_i)$ over $W_i$.
\end{lemma}
{\bf Proof:}
Since $X_{\gG_N}^*$ is integral, by Lemma \ref{L6.4.2} below, $X_{\gG_N}^*$
is in fact a product $X_{\gG_N}^*=X_{\gG_1}^*\times X_{\gG_2}^*$. The
boundary component $W_i$ is by (\ref{E35}) given by a zero-dimensional
boundary component $w$ of $X_{\gG_2}^*$, which gets modified under $\pi$,
$\pi^{-1}(X_{\gG_2}^*)=\overline{X}_{\gG_2}$. We have fibre spaces (at
least locally over $W_i$)
\[ \pi^{-1}(W_i)\stackrel{\eta}{\lra} A_i \stackrel{\zeta}{\lra} W_i,\]
where $A_i=W_i\times V_{\bfb}(\fR)/\gG_{\bfb}^0\sdprod V_{\bfb}(\fZ)$ is
the natural family of abelian varieties parameterized by $W_i$. Note that
the zero of $V_{\bfb}(\fR)$ determines a zero section $\gs_0:W_i\lra A_i,\
t\mapsto \hbox{ the image of } 0\in V_{\bfb}(\fR)$ in
$(A_i)_t=V_{\bfb}(\fR)/\gL_{\bfb}(t)$, where $\gL_{\bfb}(t)$ denotes the
lattice at the point $t$, and any element $x\in
V_{\bfb}(\fQ)$ determines locally a section $\gs_x=\gs_0+x$. Recall that
$N$ (=$N_{\bfb}^g$) is determined by an element $g\in V_{\bfb}(\fQ)$, so
the proper transform of $X_{\gG_N}^*$ in the abelian variety part of the
exceptional locus is $\left[
  X_{\gG_1}^*\right] = (\gs_0+g)(W_i)\inn A_i$, which is a global
section. \ende

If $\dim(W_i)=0$, two different situations occur, depending on whether
$\cD$ is a tube domain or not. They are (assume for the moment that $\cD$
is irreducible)
\begin{itemize} \item $\cD$ is a tube domain, then $V$ is trivial and there
  is no abelian variety; $\pi^{-1}(W_i)$ is a torus embedding.
\item $\cD$ is not a tube domain, $V$ is not trivial, and $\pi^{-1}(W_i)$
  has an abelian variety factor and a torus embedding factor.
\end{itemize}
In the first case there is not much more to say than that each irreducible
component $W_{ij}$ of $\pi^{-1}(W_i)$ meets the proper transform of
$X_{\gG_N}^*$ in a divisor on $X_{\gG_N}^*$ (which gets itself blown up at
the point). In the second case, the dimension of the abelian variety
factors and of the corresponding integral symmetric subvarieties are given
as follows:
\begin{equation}\label{E11.1} \begin{array}{l|c|c|c|}
 & \bf I_{\hbf{p,q}} & \bf II_{\hbf{n}} & \bf V \\ \hline
\dim(V) & q(p-q) & n-1 & 16 \\ \hline
\dim(\cD_N) & q(p-1) & {n-2 \choose 2} & 8,\ 10,\ 8
\end{array}
\end{equation}
At any rate, we have the following result:
\begin{theorem}\label{T11.1} The proper transform of $X_{\gG_N}^*$ on
  $\xgc$ is $\overline{X}_{\gG_N}$, a partial compactification for some
  $\gG_N$-admissible collection of polyhedral cones.
\end{theorem}
{\bf Proof:} Let $P_N$ be the parabolic in $N$, $P$ the corresponding
parabolic in $G$. Consider the decompositions (we omit the subscript
${\bfb}={\hbf{s}}$)
\[ P_N=(M_NL_N\cR_N)\sdprod \cZ_NV_N,\quad P=(ML\cR)\sdprod \cZ V.\]
Then $L_N=L$ is trivial (as the boundary component is a point), and there
is a natural inclusion $\cZ_N\inn \cZ$. Letting $C_N,\ C$ denote the
corresponding homogenous self dual cones, we have $C_N\inn C$, and both
inclusions are defined over $\fQ$. Finally we have $\gG_N=\gG\cap N$ which
implies $\gG_N\cap C_N=C_N\cap (C\cap \gG)$. We know by assumption that we
have a $\gG$-admissible cone decomposition of $C$, and since $C_N\inn C$ is
defined over $\fQ$, this gives one also for $\gG_N$, as follows from
\cite{oda}, Theorem 1.13. If $\{\gs\}$ is the cone decomposition of $C$,
then $\{\gs_N\},\ \gs_N:=\gs\cap C_N$ gives a corresponding cone
decomposition of $C_N$, and the theorem just mentioned applies. This
argument applies to each boundary component of $X_{\gG_N}^*$, and it is
clear that a $\gG$-admissible collection restricts to a $\gG_N$-admissible
collection. \ende

\subsection{Intersection}
First note the following:
\begin{lemma} Given $X_{\gG}$ and two modular subvarieties $X_1,X_2\inn
  X_{\gG}$, the intersection, if of dimension $\geq1$, is again a modular
  subvariety.
\end{lemma}
{\bf Proof:} We are given two injections defined over $\fQ$, $i_1:N_1\hra
G,\ i_2:N_2\hra G$, and commutative squares
\[\begin{array}{ccccc}\cD_{N_1} & \lra & \cD & \lla & \cD_{N_2} \\
  \downarrow & & \downarrow & & \downarrow \\ X_1 & \lra & X_{\gG} & \lla &
  X_2;
\end{array}\]
it follows that $X_1\cap X_2$ is covered by $\cD_{N_1}\cap \cD_{N_2}$ with
a corresponding injection $i_{12}:N_1\cap N_2 \hra G$, again defined over
$\fQ$. Since $X_1$ and $X_2$ are modular subvarieties, $\cD_{N_1}$ and
$\cD_{N_2}$ are by definition defined over $\fQ$, hence so is
$\cD_{N_1}\cap \cD_{N_2}$. It is also a symmetric subspace since
$\cD_{N_1}\cap \cD_{N_2}$ is totally geodesic in $\cD$. Consequently
$X_1\cap X_2$ is a modular subvariety. \ende This can be applied in
particular to the integral modular subvarieties. Hence for any two integral
modular subvarieties $X_i$, the intersection defines a (maybe empty)
modular subvariety. As there are finitely many possible intersections, from
the finite set of Corollary \ref{C6.1} we get a finite set of modular
subvarieties. Note that if $X_1$ and $X_2$ are both integral, then also the
intersection is, in the following sense: Let
$N_1=N_{\bfb_{\hbfs{1}}}^{g_1},\ N_2=N_{\bfb_{\hbfs{2}}}^{g_2}$, then
$\gG$-integral means:
\[ N_{\hbf{b$_{\hbfs{1}}$,1}}\cap \gG
=N_{\hbf{b$_{\hbfs{1}}$,1}}\cap g_1\gG g_1^{-1},\quad
N_{\hbf{b$_{\hbfs{2}}$,2}}\cap \gG =N_{\hbf{b$_{\hbfs{2}}$,2}}\cap g_2\gG
g_2^{-1},\] and $N_1\cap N_2=(g_1N_{\bfb_{\hbfs{1}}}g_1^{-1})\cap
(g_2N_{\bfb_{\hbfs{2}}}g_2^{-1})$. Hence
\begin{eqnarray*} (N_{\hbf{b$_{\hbfs{1}}$,1}}\cap N_{\hbf{b$_{\hbfs{2}}$,2}})
  \cap \gG &
  = & (N_{\hbf{b$_{\hbfs{1}}$,1}}\cap \gG)\cap (N_{\hbf{b$_{\hbfs{2}}$,2}}\cap
  \gG) \\ & = & (N_{\hbf{b$_{\hbfs{1}}$,1}}\cap g_1\gG g_1^{-1})\cap
  (N_{\hbf{b$_{\hbfs{2}}$,2}}\cap g_2\gG g_2^{-1}) \\ & = &
  (N_{\hbf{b$_{\hbfs{1}}$,1}}\cap N_{\hbf{b$_{\hbfs{2}}$,2}})\cap (g_1\gG
  g_1^{-1}\cap g_2\gG g_2^{-1}).
\end{eqnarray*}
Note that adjoining these to the integral modular subvarieties implies that
on an arithmetic quotient $X_{\gG'}$ for
  $\gG'\inn \gG$ of finite index, there is a well-defined, finite,
  non-empty set of subvarieties, all of which are either integral modular
  subvarieties incident to rational cusps or intersections of such.

Finally consider the boundary varieties $W_1$ and $W_2$ to which $X_1$ and
$X_2$ are incident. Since $X_1\cap X_2$ is a modular subvariety, it is
itself an arithmetic quotient (in general a product), and has a boundary
variety $W_{12}^*=W_1^*\cap W_2^*$.  In this sense, we
make the
\begin{definition} Let $X_1$ and $X_2$ be integral modular subvarieties,
  incident with $W_1$ and $W_2$, respectively. Then we say
  $X_{12}^*:=X_1^*\cap X_2^*$ is {\it incident to} $W_{12}^*:=W_1^*\cap
  W_2^*$.
\end{definition}
Next suppose that we are given the two parabolics, say $P_1$ and $P_2$,
which are the stabilizers of the boundary components $F_1$ and $F_1$, of
which $W_1$ and $W_2$ are the quotients, $W_1=\gG_1\bs F_1,\ W_2=\gG_2\bs
F_2$. Assume that $F_1^*\cap F_2^*\inn F_i^*,\ i=1,2$, is a
maximal boundary component in $F_i^*$. Under this assumption, the
intersection $P_1\cap P_2$ is a parabolic, associated with $F_1^*\cap
F_2^*$. Either of the inclusions $F_1^*\cap F_2^* \inn F_i^*$ determines
the parabolic which is the (non-maximal) parabolic stabilizing a flag of
two terms. Similarly, $X_1^*\cap X_2^*$ contains $F_1^*\cap F_2^*$ as a
rational boundary component, and either of the inclusions $X_1^*\cap
X_2^*\inn X_i^*$ determines a symmetric subgroup, also the stabilizer of a
flag with two terms. This is of course just $(N_1\cap N_2)\times
\cZ_G(N_1\cap N_2)$, where $N_i$ is the group giving rise to $\cD_{N_i}$,
covering $X_i$. So we have: $N_i$ incident with $P_i,\ i=1,2$,
$P_{12}:=P_1\cap P_2$ a parabolic, then $(N_1\cap N_2)\times \cZ_G(N_1\cap
N_2)$ is incident to $P_{12}$.

\subsection{Moduli interpretation}
In this section we suppose the algebraic group $G$ comes from a moduli
problem of {\sc Pel} structures, and will discuss the moduli-theoretic
description of the modular subvarieties $X_{\gG_N}$ and the corresponding
arithmetic quotients $\xgs$. We then also briefly describe the notion of
incidence from this point of view.

\subsubsection{{\sc Pel} structures}
\paragraph{\ }
Let $V$ be an abelian variety over $\fC$, $End(V)$ the endomorphism ring and
$End_{\fQ}(V)=End(V)\otimes_{\fZ}\fQ$ the endomorphism algebra. A
polarization, i.e., a linear equivalence class of ample divisors giving a
projective embedding of $V$, gives rise to a positive involution on
$End_{\fQ}(V)$, the so-called Rosatti involution:
\begin{eqnarray}\label{e96.1} \grr:End_{\fQ}(V) & \lra &  End_{\fQ}(V) \\
\phi & \mapsto & \phi^{\grr}. \nonumber
\end{eqnarray}
If $A$ is a central simple algebra over $\fQ$, an involution * on $A$ is called
{\em positive}, if $tr_{A|\fQ}(x\cdot x^*)>0$ for all $x\in A,\ x\neq
0$, where $tr_{A|\fQ}$ denotes the reduced trace.
Assuming $(A,*)$ to be simple with positive involution, the
$\fR$-algebra $A(\fR)$ is isomorphic to one of the following (see
\cite{shimura2}, Lemma 1)
\begin{itemize}\item[(i)] $M_r(\fR)$ with involution $X^*={^tX}$;
\item[(ii)] $M_r(\fC)$ with involution $X^*={^t\-X}$, where $^{-}$ is complex
conjugation;
\item[(iii)] $M_r(\fH)$ with involution $X^*={^t\-X}$, where $^{-}$ is
quaternionic conjugation.
\end{itemize}
The algebras $A$ occuring in (i) and (iii) are central simple over $\fR$,
while those of (ii) are central simple over $\fC$. The $\fQ$-algebra
$A$ itself is a $\fQ$-form of one of these. The central simple algebras
$A$ over $\fQ$ are known to be the $M_n(D)$, where $D$ is a division
algebra over $\fQ$. If the algebra $A$ has
a positive involution, the same holds for $D$. The division algebras $D$
which can occur are also known.
\begin{proposition}\label{p96.1} Let $D$ be a division algebra over $\fQ$
with a positive involution. Then $D$ occurs in one of the following cases:
\begin{itemize}\item[I.] A totally real algebraic number field $k$;
\item[II.] $D$ a totally indefinite quaternion algebra over $k$;
\item[III.] $D$ a totally definite quaternion algebra over $k$;
\item[IV.] $D$ is central simple over $K$ with a $K|k$ involution
of the second kind, where $K$ is an imaginary quadratic extension of $k$.
\end{itemize}
\end{proposition}
In case III the canonical involution on $D$ is the unique positive
involution, while in case II the positive involutions correspond to $x\in D$
such that $x^2$ is totally negative in $k$. If the algebra $D$ has an
involution of the second kind it is easy to see
that it admits a positive one.
It follows from the fact that $End_{\fQ}(V)$ is a semisimple algebra over
$\fQ$ with a positive involution that each simple factor is a total matrix
algebra $M_n(D)$, with $D$ as in the proposition.

\paragraph{ \ }
Let $(A,*)$ be a semisimple algebra over $\fQ$ with positive involution,
and let
\begin{equation}\label{e96.2} \Phi:A\lra GL(n,\fC) \end{equation}
be a faithful representation. Shimura considers data $\cP=(V,\cC,\gt)$ and
$\{A,\Phi,*\}$ and defines the notion of {\em polarized abelian variety
of type} $\{A,\Phi,*\}$ by the conditions:
\begin{equation}\label{e96.3}\begin{minipage}{14cm}
\begin{itemize}\item[(i)] $V$ is an abelian variety over $\fC$, $\cC$ is a
polarization;
\item[(ii)] $\gt:A\stackrel{\sim}{\lra} End_{\fQ}(V)$ is an algebra
isomorphism, and for $\gt(x):\~V\lra \~V$ (the $\~{ }$ denoting the universal
cover, i.e., $\~V$ is a complex vector space) one has $\gt(x)=\Phi(x)$;
\item[(iii)] the involution $\grr$ determined by $\cC$ as in (\ref{e96.1})
coincides on $\gt(A)$ with the involution coming from $(A,*)$, i.e.
$\gt(x)^{\grr}=\gt(x^*)$.
\end{itemize}
\end{minipage}\end{equation}
The condition (ii) is to be understood as follows. Fixing an isomorphism
\begin{equation}\label{e97.-1}
\psi:V\isom \fC^n/\gL,
\end{equation}
each $a\in End_{\fQ}(V)$ is represented by a linear
transformation of $\fC^n$ preserving $\gL$; that is each $a$ can be
represented by a matrix, and $\gt(x)=a$ is the matrix corresponding to $x\in
A$ via $\gt$. Recall also that a complex torus $\fC^n/\gL$ is an abelian
variety if and only if there exists a {\em Riemann form}: each positive (1,1)
form $\go$ gives rise to a skew symmetric matrix $(q_{ij})$:
$$\go=\sum q_{ij}dx_i\wedge dx_j,$$
where the $x_i$ are canonical coordinates on $\fC^n$. Hence if we fix a
positive divisor $C\inn V$, it determines an involution as in (\ref{e96.1})
{\em and} a Riemann form $E_C(x,y)$ on $\fC^n/\gL$, and these are related by
\begin{equation}\label{e97.0} E_C(\psi(a)x,y)=E_C(x,\psi(a^{\grr})y),
\end{equation}
where for $a\in End_{\fQ}(V)$, $\psi(a)$ denotes the matrix
representation for $a$ arising from the identification $\psi:V\isom
\fC^n/\gL$ in (\ref{e97.-1}).

Let $(V,\cC,\gt)$ be an abelian variety of type $(D,\Phi,*)$ with $D$ a
division algebra, so that $(D,*)$ is one of the algebras of Proposition
\ref{p96.1}. In the notations used there, put
\begin{equation}\label{e97.6} [k:\fQ]=f,\quad [D:K]=d^2,\hbox{ if $D$ is of
type IV} \end{equation}
defining the numbers $f$ and $d$. Let $n$=dim$(V)$; then, assuming $D$ to be
a division algebra, $2n$ is a multiple of $[D:\fQ]$, i.e., $2n=[D:\fQ]m$. Note
that $[D:\fQ]=f$ for type I, $[D:\fQ]=4f$ for types II and III, while
$[D:\fQ]=2d^2f$ if $D$ is of type IV. For the existence of $(V,\cC,\gt)$ of
type $(D,\Phi,*)$, certain restrictions are placed on $\Phi$; we assume these
are fulfilled. So under the isomorphism $\gt$, each $x\in D$ is represented
by the matrix $\Phi(x)$. This makes the lattice $\gL$ with $V\isom
\fC^n/\gL$, tensored with $\fQ$, a (left) $D$-module, i.e.,
\begin{equation}\label{e97.1} Q:=\fQ\cdot \gL = \sum_1^m\Phi(D)\cdot x_i
\end{equation}
for a suitable set of vectors $x_i$. But this is the same as saying there
exists a $\fZ$-lattice $\cM\inn D$, such that
\begin{equation}\label{e97.2} \gL=\{\sum_1^m\Phi(a_i)x_i \Big|
(a_1,\ldots,a_m)\in \cM \}.
\end{equation}
If $D$ is central over $K$, then $\cM$ is clearly also an $\cO_K$-lattice in
$D$. The integrality of the Riemann form can be expressed in terms of
$tr_{D|K}$:
\begin{equation}\label{e97.3} E_C(\sum_1^m\Phi(a_i)x_i,\sum_1^m\Phi(b_j)x_j)
= tr_{D|K}(\sum_{i,j}^ma_it_{ij}b_j^*),
\end{equation}
and $T=(t_{ij})\in M_m(D)$ is a skew-hermitian matrix:
\begin{equation}\label{e97.4} T^*=-T,
\end{equation}
where $T^*$ denotes the matrix $(t_{ji}^*)$, where $*$ is the involution on
$D$. For the lattice $\cM$ one
has
\begin{equation}\label{e97.5} tr_{D|K}(\cM T \cM^*)\inn \fZ.
\end{equation}

\paragraph{ \ }
Hence to each $(V,\cC,\gt)$ of type $(D,\Phi,*)$ one gets a $*$-skew hermitian
$T\in M_m(D)$ and a lattice $\cM\inn D$. To this situation there is a
naturally associated $\fQ$-group. On the vector space $D^m$ we consider
\begin{equation}\label{e98.1} G(D,T):=\{g\in GL(D^m) \Big| gTg^*=T\},
\end{equation}
the symmetry group of the $*$-skew hermitian form determined by $T$. It is now
easy to determine the $\fR$-group:
\begin{equation}\label{e98.2} G(D,T)(\fR)=\left\{ \begin{minipage}{10cm}
\begin{tabbing}
Type I:\quad\quad \= $Sp(m,\fR)\times \cdots \times Sp(m,\fR)$ ($m$ is even) \\
Type II:\> $Sp(2m,\fR)\times \cdots \times Sp(2m,\fR)$ \\
Type III:\> $SO^*(2m)\times \cdots \times SO^*(2m)$ \\
Type IV: \> $U(p_1,q_1)\times \cdots \times U(p_g,q_g),$
\end{tabbing}
\end{minipage} \right.
\end{equation}
where the number of factors is in each case $f$, and $p_{\nu}+q_{\nu}=md$,
and $(p_{\nu},q_{\nu})$ is the signature corresponding to the $\nu^{th}$ real
prime.
For each $\nu$, there is a matrix $W_{\nu}$ which trasforms $T_{\nu}$
into the standard form, i.e.,
\begin{eqnarray}\label{e98.3}
W_{\nu}T_{\nu}^{-1}{^tW_{\nu}} & = & \left(\begin{array}{cc} 0
& 1_l \\ -1_l & 0 \end{array} \right), l={m\over 2} \hbox{ for Type I, $l=m$
for Type II}; \\
\label{e98.4} W_{\nu}T_{\nu}^{-1}W^*_{\nu} & = & -i\left(\begin{array}{cc}
-1_m & 0 \\ 0 & 1_m \end{array} \right), \hbox{ Type III}; \\
\label{e98.5} W_{\nu}(iT_{\nu}^{-1})W^*_{\nu} & = & \left(\begin{array}{cc}
1_{p_{\nu}} & 0 \\ 0 & -1_{q_{\nu}}\end{array}\right), \hbox{ Type IV}.
\end{eqnarray}
Let $\cD=\cD_{(D,T)}$ denote the domain determined by $G(D,T)(\fR)$ (actually
a particular unbounded realisation of this domain, see \cite{shimura2}, 2.6).
Then $\cD=\prod\cD_{\nu}$, and
$z_{\nu}\in \cD_{\nu}$ gives rise to a normalised
period (i.e., one of the form $(1,\gO)$) for an abelian variety, by setting
$X_{\nu}=Y_{\nu}\-{W}_{\nu}$, where
\begin{eqnarray}\label{e98.6} Y_{\nu} & = & \left(\begin{array}{cc} z_{\nu} &
1_l \\ \-z_{\nu} & 1_l \end{array}\right), l={m \over 2}, \hbox{ Type I,
$l=m$, Type II}; \\
\label{e98.7} Y_{\nu} & = & \left(\begin{array}{cc} -z_{\nu} & 1_m \\
1_m & \-z_{\nu} \end{array}\right), \hbox{ Type III}; \\
\label{e98.8} Y_{\nu} & = & \left( \begin{array}{cc} 1_{p_{\nu}} & z_{\nu}
\\ ^t\-z_{\nu} & 1_{q_{\nu}}\end{array}\right), \hbox{ Type IV.}
\end{eqnarray}
The matrix $X_{\nu}$ determine $m$ vectors $x_1,\ldots x_m$ of $\fC^n$ (in a
rather complicated fashion, see formulas (17)-(20) in \cite{shimura2}),
which determine a lattice $\gL=\gL(z,T,\cM)$ by the formula in equations
(\ref{e97.1})-(\ref{e97.2}) above.

Note that the representation $\Phi$ contains the representations
$\chi_{\nu}$= projection on the $\nu$th real factor with multiplicities. For
Type IV, $p_{\nu}+q_{\nu}=md$, and $p_{\nu}$=multiplicity of $\chi_{\nu}$
while $q_{\nu}$=multiplicity of $\-{\chi}_{\nu}$. For things to work out one
must therefore assume, in case of Type IV, that $iT^{-1}_{\nu}$ has the {\em
same} signature $(p_{\nu},q_{\nu})$ as occurs in $\Phi$. With this
restriction, the following holds:
\begin{theorem}[\cite{shimura2}, Thm.~1]\label{t99.1}
For every $z\in \cD=\cD_{(D,T)}$, and every lattice $\cM\inn D$, we get a
polarized abelian variety $V_z=\fC^n/\gL(z,T,\cM)$ of type $(D,\Phi,*)$, and
conversely, every such $V$ is of the form $V=\fC^n/\gL(z,T,\cM)$ for some
$z\in \cD_{(D,T)},\ \cM\inn D$ a lattice.
\end{theorem}
\paragraph{ \ }
The lattice $\cM\inn D$ gives rise to an arithmetic subgroup
\begin{equation}\label{e99.1} \gG=\gG_{(D,T,\cM)}=\{g\in G(D,T) \Big|
g\cM\inn \cM\}
\end{equation}
as discussed in section \ref{s83.1}.
If one defines an isomorphism $\phi:V_z\lra V_{z'}$ of
two abelian varieties of type $(D,\Phi,*)$ as an isomorphism of the
underlying varieties, such that $\phi^{-1}(\cC')=\cC$ and
$\phi\gt(a)=\gt'(a)\phi,$ for all ${a\in D}$, then one has
\begin{theorem}[\cite{shimura2}, Thm.~2]\label{t99.2} The arithmetic quotient
$\xg=\gG\bs\cD_{(D,T)}$ is the moduli space of isomorphism classes of abelian
varieties $V_z=\fC^n/\gL(z,T,\cM)$ of type $\{(D,\Phi,*),(T,\cM)\}$, where
$\gG$ is the arithmetic group of (\ref{e99.1}).
\end{theorem}
Moreover, one calls two such pairs $(T_1,\cM_1),\ (T_2,\cM_2)$ {\em
equivalent}, if $\exists_{U\in M_m(D)}$, such that $UT_2U^*=\gd T_1$ for some
positive $\gd\in \fQ$ and $\cM_1U=\cM_2$. Equivalent pairs give rise to
isomorphic families of abelian varieties (\cite{shimura2}, Prop.~4). Summing
up, $*$-skew hermitian matrices $T\in M_m(D)$ determine certain $\fQ$-groups,
lattices $\cM\inn D$ determine certain arithmetic groups, and the
corresponding arithmetic quotients are moduli spaces for certain families of
abelian varieties.
\begin{remark}\label{r99.1} The complex multiplication by $\cM$ describes the
endomorphism ring. The {\it automorphisms} determined by
$\cM$ are the invertible
elements, i.e., $\Aut(V)=\cM^*$, the group of units.
\end{remark}

One can also accomodate level structures in this settup, introduced in
\cite{shimura3}, cf.~also \cite{shimura4}. This is done by fixing $s$ points
$y_1,\ldots,y_s$ in the $D$-module $Q$, as in (\ref{e97.1}), and $s$ points
$t_1,\ldots,t_s$ of the abelian variety $V$. One requires that the map $\psi$
of (\ref{e97.-1}) maps the $y_i$ onto the $t_i$. More precisely,
\begin{definition}\label{d99z.1}
Let $Q$ be a $D$-vector space of dimension $m$, and
$\scM\inn Q$ a $\fZ$-lattice. Consider a conglomeration:
$$\scT:=\{(D,\Phi,*),(Q,T,\scM);y_1,\ldots,y_s\},$$
where $(D,\Phi,*)$ is as above, $(Q,T,\scM)$ is a $D$-vector space with
lattice $\scM$ and $*$-skew
hermitian ($D$-valued) form $T$ on $Q$, and $y_i$ points in
$Q$. This is called a {\sc Pel}{\em -type}. Consider a conglomeration:
$$\scQ:=\{(V,\cC,\gt);t_1,\ldots,t_s\},$$
where $(V,\cC,\gt)$ is a polarized abelian variety with analytic coordinate
$\gt$ as above and $t_i$ are points of {\em finite order} on $V$. This is
called a {\sc Pel}-{\em structure}. Then $\scQ$
is {\em of type} $\scT$, if there exists a commutative diagram
\begin{equation}\label{e99z.1} \begin{array}{ccccccccc} 0 & \lra & \scM & \lra
& Q(\fR) & \lra & Q(\fR)/\scM & \lra & 0 \\
 & & \downarrow & & f\downarrow & & \downarrow & & \\
0 & \lra & \gL & \lra & \fC^n & \stackrel{\psi}{\lra} & V & \lra & 0
\end{array},
\end{equation}
satisfying the conditions\begin{itemize}\item[(i)] $\psi$ gives a holomorphic
isomorphism (strictly speaking, this is the $\psi^{-1}$ of above);\item[(ii)]
$f$ is an $\fR$-linear isomorphism, and $f(\scM)=\gL$;
\item[(iii)] $f(ax)=\Phi(a)f(x)$, and $\Phi(a)$ defines $\gt(a)$ for every
$a\in D$ as (\ref{e96.3}), (ii);
\item[(iv)] $C\in\cC$ determines a Riemann form $E_C$ as in (\ref{e97.0}).
\end{itemize}
\end{definition}
Note that the finite set of points $y_i$ and $t_i$ come both equipped with a
form; on the former the form $T$, and the Riemann form $E_C$ on the latter.
These forms are preserved under the isomorphism.

There is a natural notion of isomorphism of abelian varieties with {\sc Pel}
structures. Given two {\sc Pel}-structures $\scQ$ and $\scQ'$, an isomorphism
$\phi:V\lra V'$ is an {\em isomorphism} from $\scQ$ to $\scQ'$, if
$\phi\gt(a)=\gt'(a)\phi$ for all $a\in D$, and $\phi(t_i)=t_i'$ for all $i$.
\begin{definition}\label{d99z.2} A {\sc Pel}-type $\scT$ is {\em equivalent}
to a {\sc Pel}-type $\scT'$, if $D=D'$, $*=*'$, $s=s'$, $\Phi$ and $\Phi'$
are equivalent as representations of $D$, and there is a $D$-linear
automorphism $\mu$ of $Q$ such that $T'(x\mu, y\mu)=T(x,y),\ \scM\mu=\scM',
y_i\mu\equiv y_i'$mod$\scM'$  for all $i$. If $\scQ$ is of type $\scT$, then
$\scQ$ is also of type $\scT'$ if and only if $\scT$ and $\scT'$ are
equivalent. A {\sc Pel}-type $\scT$ is called {\em admissible}, if there
exists at least one {\sc Pel}-structure of that type.
\end{definition}
One has an anolgue of Theorems \ref{t99.1} and \ref{t99.2} in this situation
also.
\begin{theorem}[\cite{shimura4}, Thm.~3]\label{t99z.1} For every
admissible {\sc Pel}-type $\scT$ there exists a bounded symmetric domain $\cD$
(this is the same domain as in Theorem \ref{t99.1}) such that the statement
of Theorem \ref{t99.1} holds in this situation, and every {\sc Pel}-structure
$\scQ$ of type $\scT$ occurs in this family.
\end{theorem}
Now define a corresponding arithmetic group as follows:
\begin{equation}\label{e99z.2} \gG=\{g\in G(D,T) \Big| \scM g=\scM,\
y_i g\equiv y_i\hbox{mod}\scM,\ i=1,\ldots,s\}
\end{equation}
Then the analogue of Theorem \ref{t99.2} is
\begin{theorem}[\cite{shimura4}, Thm.~4]\label{t99z.2}
Two members of the family of Theorem \ref{t99z.1} corresponding to points
$z_1,z_2\in \cD$ are isomorphic if and only if $z_1=\gg(z_2)$ for some
$\gg\in \gG$, $\gG$ as in (\ref{e99z.2}).
\end{theorem}
In Table \ref{table17} we list the data for each of the cases II, III and
IV of (\ref{e98.2}).

\begin{table}\caption{\label{table17}
 Rational groups for {\sc Pel}-structures.}
\medskip\begin{center}
\begin{tabular}{|c|c|c|c|}\hline & Type II & Type III & Type IV \\ \hline
$D$ & \begin{minipage}{3.5cm} A totally \\ indefinite quaternion algebra
over $\fQ$ \end{minipage} &
\begin{minipage}{3.5cm} a totally \\
definite quaternion \\
algebra over $\fQ$
\end{minipage} &
\parbox{4cm}{\footnotesize
simple division algebra, central over $K$, an imaginary
quadratic extension of $\fQ$, with an involution of the second kind. One may
assume $D$ to be a cyclic algebra} \\ \hline
$d$ & 2 & 2 & $d$ \\ \hline
dim($V$) & $2m$ & $2m$ & $d^2m$ \\ \hline
Tits index & $C_{m,s}^{(2)}$ & ${^iD}_{m,s}^{(2)}, i=1,2$ &
 $^2A^{(d)}_{dm-1,s}$ \\ \hline
\end{tabular}

\medskip The types listed are absolutely $\fQ$-simple, the
$\fQ$-rank is $s$, and this is the Witt index of the
$\pm$hermitian form.
\end{center}
\end{table}

\subsubsection{Modular subvarieties}
We continue with the notations above, $G,\ S,\ P_{\bfb}$ and $N_{\bfb}$
being fixed, $\gG$ an arithmetic group satisfying (\ref{E5.1}), and
$\gG'\inn G(\fQ)$ another arithmetic group. We consider the arithmetic
quotient $X_{\gG'}$, its Baily-Borel embedding $X_{\gG'}^*$, and a smooth
toroidal embedding $\overline{X}_{\gG'}$. Let $N\inn G$ be a rational
symmetric subgroup, $N=N_{\bfb}^g$ for some $b=1,\ldots, s$ and some $g\in
G(\fQ)$. Let us first suppose for the boundary point in question $F_{\bfb}$
that $\dim(F_{\bfb})>0$. Under this assumption we know that $N_{\bfb}$ is a
product
\[ N_{\bfb}=L_{\bfb}\times \cZ_G(L_{\bfb}),\]
and $L_{\bfb} (\fR)^0=(\Aut(F_{\bfb}))^0$. This implies immediately that
the domain $\cD_{N_{\hbfs{b}}}$ is also a product,
\[ \cD_{N_{\hbfs{b}}}\cong \cD_1\times \cD_2.\]
Let $\imath_i:\cD_i\hra \cD_1\times \cD_2$ be the natural inclusion, and
consider the inclusion $\eta:\cD_{N_{\hbfs{b}}}\hra \cD$. Then
$\eta(\imath_i(\cD_i))$ is a symmetric subdomain, which itself has an
interpretation in terms of {\sc Pel} structures, which is a sub-{\sc Pel}
structure of that attached to $\cD$. Let us now explain this for the
individual cases. For {\sc Pel} structures, only the domains of type $\bf
I_{\hbf{p,q}},\ II_{\hbf{n}},\ III_{\hbf{n}}$ occur, types {\bf U.1, U.2,
  O.2, S.1, S.2}.
\begin{itemize}\item[{\bf U.1}]: If $b<t$, then $F_{\bfb}$ is of type $\bf
  I_{\hbf{p-b,q-b}}$, and $\cD_{N_{\hbfs{b}}}$ is of the type $\bf
  I_{\hbf{p-b,q-b}}\times I_{\hbf{b,b}}$. The moduli interpretation is
  complex multiplication on abelian $(p+q)$-folds. In the locus $\bf
  I_{\hbf{p-b,q-b}}\times I_{\hbf{b,b}}$, the variety $A^{p+q}$ splits into
  $A^{p+q-2b}\times A^{2b}$, where the complex multiplication has signature
  $(p-b,q-b)$ and $(b,b)$, respectively.

  If $b=s=t$, $F_{\hbf{t}}$ is a point, $\cD_{N_{\hbfs{b}}}$ is of type $\bf
  I_{\hbf{p-1,q}}$. Here the abelian variety $A^{p+q}$ splits off an
  elliptic curve, $A^{p+q}=A^{p+q-1}\times A^1$. Since $A^1$ has no
  moduli, only the moduli of $A^{p+q-1}$ contributes.
\item[{\bf U.2}]: If $b<t$, then $F_{\bfb}$ is of type $\bf
  I_{\hbf{p-db,q-db}}$ and $\cD_{N_{\hbfs{b}}}$ is of type $\bf
  I_{\hbf{p-db,q-db}}\times I_{\hbf{db,db}}$. The moduli involved here is a
  degree $d$ division algebra $D$, central simple over $K$ with
  $K|k$-involution, as endomorphism algebra. In the locus $\bf
  I_{\hbf{p-db,q-db}}\times I_{\hbf{db,db}}$ the abelian variety $A^{p+q}$
  splits $A^{p+q}=A^{p+q-2db}\times A^{2db}$, and each factor retains the
  endomorphisms by $D$.

  If $b=s=t$, again $F_{\hbf{t}}$ is a point, $\cD_{N_{\hbfs{b}}}$ is of type
  $\bf I_{\hbf{q,q}}\inn I_{\hbf{p,q}}$. Here the abelian variety splits as
  $A^{2q}\times A^{p-q}$, where the endomorphism ring on $A^{p-q}$ is
  definite, and contributes no moduli. In this case the only moduli
  contributing is the modulus of $A^{2q}$.
\item[{\bf O.2}]: If $b<s$ or $s<[{n\over 2}]$, then $F_{\bfb}={\bf
  II_{\hbf{n-2b}}},\ \cD_{N_{\hbfs{b}}}=\bf II_{\hbf{n-2b}}\times
  II_{\hbf{2b}}$, and the splitting is evident. For $b=t$, $F_{\hbf{t}}=\bf
  II_{\hbf{0}}$ ($n$ even) or $\bf II_{\hbf{1}}$ ($n$ odd), both of which
  are points. Then for $n$ even, $N_{\hbf{t}}$ is a polydisc by definition,
  $\bf II_{\hbf{2}}\times \cdots \times II_{\hbf{2}}$, and again the
  splitting is evident, this time as a product of
  abelian surfaces with multiplication
  by the quaternion division algebra $D$. For $n$ odd, $N_{\hbf{t}}$ is of
  type $\bf II_{\hbf{n-1}}$, and the splitting of $A^{2n}$ is as
  $A^{2n}\cong A^{2n-s}\times A^2$.
\item[{\bf S.1}]: If $b<t$, then $F_{\bfb}$ is of type $\bf
  III_{\hbf{n-b}}$, the subdomain $\cD_{N_{\hbfs{b}}}$ is of type $\bf
  III_{\hbf{n-b}}\times III_{\hbf{b}}$, $A^n$ splits $A^n=A^{n-b}\times
  A^b$. If $b=s=t$, then $F_{\hbf{t}}$ is a point, $\cD_{N_{\hbfs{t}}}
  \cong ({\bf
  III_{\hbf{1}}})^n$. Here the abelian variety splits into a product of
  elliptic curves.
\item[{\bf S.2}]: If $b<t$, $F_{\bfb}={\bf III_{\hbf{n-2b}}}$,\
  $\cD_{N_{\hbfs{b}}}={\bf III_{\hbf{n-2b}}\times III_{\hbf{2b}}}$,
  where $b<
    [{n\over 2}]$. If $b=s=t$, $\cD_{N_{\hbfs{t}}}=({\bf
      III_{\hbf{2}}})^{n\over 2}$, ($n$ even follows from $s=t$). Once
    again, in both cases the moduli-theoretic meaning is evident.
  \end{itemize}
This explains the expression of sub-{\sc Pel} structures, and we have
established
\begin{proposition}\label{p6.4.1} For any subdomain $\cD_{N_{\hbfs{b}}}\inn
  \cD$, the corresponding abelian varieties split in the manner described
  above.
\end{proposition}
{\bf Proof:} We prove a typical case and leave the verification of the
other cases to the reader. We will do case {\bf U.2}. For this we consider
the matrix $Y_{\nu}$ of (\ref{e98.8}). We may assume that the realisation
of the domain $\cD$ is such that for the subdomain $\cD_N$, the
corresponding $z_{\nu}$ splits, i.e.,
\begin{equation}\label{E161.1}
   z_{\nu}\in \cD_N \Ra z_{\nu}=\left(\begin{array}{cc} z_{\nu,1} & 0 \\ 0
    & z_{\nu,2}
  \end{array}\right),
\end{equation}
where $z_{\nu,1}\in \cD_1$ and $z_{\nu,2}\in \cD_2$ for the decomposition
$\cD_N=\cD_1\times \cD_2$. In this case, $\cD_1$ is of type $\bf
I_{\hbf{p-jd,q-jd}}$, while $\cD_2$ is of type $\bf I_{\hbf{jd,jd}}$.
We need the vectors $x_1,\ldots, x_m$ determined by $X_{\nu}$, given in
this case by the formula (20) in \cite{shimura2}
\[ X_{\nu}=\left[\begin{array}{cccccccccc} u_{11}^{\nu} & \cdots &
    u_{m1}^{\nu} & u_{12}^{\nu} & \cdots & u_{m2}^{\nu} & \cdots &
    u_{1d}^{\nu} & \cdots & u_{md}^{\nu} \\
    \overline{v}_{11}^{\nu} & \cdots & \overline{v}_{m1}^{\nu} &
    \overline{v}_{12}^{\nu} & \cdots & \overline{v}_{m2}^{\nu} & \cdots &
    \overline{v}_{1d}^{\nu} & \cdots & \overline{v}_{md}^{\nu}
  \end{array}
  \right]. \]
Here the vectors $x_i$ are given by ${^tx}_i^{\nu}=({^tu}_{i1}^{\nu}\cdots
{^tu}_{1d}^{\nu} {^tv}_{i1}^{\nu} \cdots {^tv}_{id}^{\nu})$ and
$u_{ik}^{\nu}\in \fC^{p_{\nu}},\ v_{ik}^{\nu}\in \fC^{q_{\nu}}$. Now
from the particular form of our $z_{\nu}$, we can conclude that also the
vectors $x_i$ have a special form.
Indeed, comparing the above with (\ref{E161.1}), we see that for
$W_{\nu}=id$ we have
\[ X_{\nu} = \left( \begin{array}{cc|cc}\boldOne_{\hbf{p$_{\nu}$-jd}} & 0 &
    z_{\nu,1} & 0 \\
  0 & \boldOne_{\hbf{jd}} & 0 & z_{\nu,2} \\ \hline
  {^t\overline{z}}_{\nu,1} & 0 & \boldOne_{\hbf{q$_{\nu}$-jd}} & 0 \\
  0 & {^t\overline{z}}_{\nu,2} & 0 & \boldOne_{\hbf{jd}}
\end{array}\right), \]
so that the vector ${^tx}_1$, for example, has the form\footnote{for
  convenience the transpostition $t$ is placed to the right of the vector
  in this expression}
\[ {^tx}_1 = \left( \left(
    \begin{array}{c} 1 \\ 0 \\ \vdots \\ \vdots
      \\ \vdots \\ 0
    \end{array}\right)^{t}\cdots
   \left(\begin{array}{c} 0 \\ \vdots
      \\ 1 \\ 0 \\ \vdots \\ 0
    \end{array}\right)^{t}
  \left(\begin{array}{c} z_{\nu,1}^{(1)} \\ 0
    \end{array}
    \right)^{t} \cdots
  \left(\begin{array}{c} z_{\nu,1}^{(p-(j-1)d+1)} \\ 0
      \end{array}
      \right)^{t} \cdots
  \left( \begin{array}{c} 0 \\ z_{\nu,2}^{(1)}
        \end{array}
        \right)^{t} \cdots
   \left( \begin{array}{c} 0 \\ z_{\nu,2}^{((j-1)d+1)}
          \end{array}
          \right)^{t} \right),\]
where $z_{\nu,j}^{(k)}$ denotes the $k^{th}$ column of $z_{\nu,j}$, and
similarly for the other ${^tx}_i$.  From this it follows that the lattice
$\gL$ of (\ref{e97.2}) splits, $\gL=\gL_1\oplus \gL_2$, where $\gL_1$ and
$\gL_2$ are orthogonal to each other, each being itself a normalized period
matrix
\[ \gL_1=\left(\begin{array}{c} 1 \\ 0 \\ \vdots \\ 0 \\ \vdots \\ 0
  \end{array}
  \right)\fZ\oplus \cdots \oplus \left(\begin{array}{c} 0 \\ \vdots \\ 1 \\
      0 \\ \vdots
    \end{array}
    \right)\fZ \oplus z_{\nu,1}^{(1)}\fZ \oplus \cdots \oplus
    z_{\nu,1}^{(q-jd)}\fZ, \gL_2=\left(\begin{array}{c} 0 \\ \vdots \\ 1 \\
        0 \\ \vdots \\ 0
      \end{array}
      \right) \fZ\oplus \cdots \oplus \left(\begin{array}{c} 0 \\ \vdots \\
          0 \\ 1
        \end{array}
        \right) \fZ\oplus z_{\nu,2}^{(1)}\fZ\oplus \cdots \oplus
        z_{\nu,2}^{(jd)}\fZ. \]
The proposition follows from this for the case that $W_{\nu}=id$. Finally
we note that if $W_{\nu}\neq id$, this does not influence the reasoning
above, and the splitting remains (only the polarization is no longer
principal). In particular, the case $\dim(F_{\bfb})=0$, which occurs for
$b=s,\ sd=q_{\nu}$, is covered by the above,
\[ z_{\nu}=\left(\begin{array}{cc} \boldOne & 0 \\ 0 & z_{\nu,2}
  \end{array}
  \right)\]
with $z_{\nu,2}\in \bf I_{\hbf{sd,sd}}$. \ende

We now consider conjugates $N=N_{\bfb}^g$. The following two lemmas apply
to any $G$ as considered in this paper so we assume for the moment only
that $G$ is $\fQ$-simple of hermitian type, $\gG$ fulfills (\ref{E5.1}) and
$\gG'\inn G(\fQ)$ is arithmetic.
\begin{lemma}\label{L6.4.1} If $g\in G(\fQ)$, $N=N_{\bfb}^g$, then the
  modular subvariety $X_{\gG'_N}\inn X_{\gG'}$ is a finite quotient of a
  product. Consequently, for all $\cD_N$, $N$ rational symmetric,
  the arithmetic subvariety $X_{\gG'_N}$ is in the locus of isomorphism
  classes of abelian varieties which are isogenous to products, i.e., are
  not simple.
\end{lemma}
{\bf Proof:} Since $g\in G(\fQ)$, we see that $\cD_N$ is $\fQ$-equivalent
to $\cD_{N_{\hbfs{b}}}$. For $\cD_{N_{\hbfs{b}}}$
the statement follows from the
fact that $N_{\bfb}=N_1\times N_2$ is a product over $\fQ$, so for $g\in
G(\fQ)$ it is likewise true for $N=N_{\bfb}^g$. Consequently, the action of
$\gG$ on $\cD_N$ is up to a finite action a product action. The second
statement follows from this, as on the finite cover which is a product, the
splitting property follows as discussed above. \ende
Now suppose $N$ is in fact $\gG$-integral, i.e., $N_1\cap \gG
=g(N_{\hbf{b,1}}\cap \gG)g^{-1}$.
\begin{lemma}\label{L6.4.2} If $N$ is $\gG$-integral, then the discrete
  subgroup $N\cap \gG$ is in fact a product, $\gG_N=N\cap \gG=\gG_1\times
  \gG_2$, $\gG_i\inn \Aut(\cD_i),\ i=1,2$.
\end{lemma}
{\bf Proof:}
We note that there is a natural inclusion $N\cap \gG\inn
g\gG_{\hbf{b,1}}g^{-1}\times g\gG_{\hbf{b,2}}g^{-1}$, and since $N_1\cap
\gG\inn N\cap \gG$ is equal to $g\gG_{\hbf{b,1}}g^{-1}$ we get the exact
diagram
\[\begin{array}{ccccccccc} & & 1          & & 1         \\
                           & & \downarrow & & \downarrow \\
1&\lra & \gG_1 & \lra & g\gG_{\hbf{b,1}}g^{-1} & \lra &  1 \\
 &     & \downarrow & & \downarrow             &      & \downarrow \\
1&\lra & N\cap \gG & \lra & g\gG_{\hbf{b,1}}g^{-1}\times
g\gG_{\hbf{b,2}}g^{-1} & \lra & K_1 & \lra & 1 \\
& & \downarrow & & \downarrow & & \downarrow \\
1&\lra & Q & \lra & g\gG_{\hbf{b,2}}g^{-1} & \lra & K_2 & \lra & 1 \\
& & \downarrow & & \downarrow & & \downarrow \\
& & 1& \lra & 1 & \lra & 1,
\end{array}
\]
and the splitting $N\cap \gG\cong \gG_1\times Q$ follows from that of
$g\gG_{\hbf{b,1}}g^{-1}\times g \gG_{\hbf{b,2}}g^{-1}$: $Q$ is a subgroup
of finite index in $g\gG_{\hbf{b,2}}g^{-1}$, and giving the
injection $\gG_1\times
Q\hra g\gG_{\hbf{b,1}}g^{-1}\times g \gG_{\hbf{b,2}}g^{-1}$ is equivalent
to giving the injection $Q\hra g\gG_{\hbf{b,2}}g^{-1}$. \ende

It may well be that $N$ is in fact $\gG$-integral if and only if
$X_{\gG_N}$ is a product, but I have no argument for this. At any rate,
this can now be applied to derive the moduli interpretation of $X_{\gG'_N}$
for $N$ $\gG$-integral.

Applying the two lemmas above again in the situation that $G$ corresponds
to a {\sc Pel}-structure yields the following.
\begin{theorem}\label{t6.4.1} Let $G,\ S,\ P_{\bfb},\ N_{\bfb}$ and $\gG$
  be as above ($b<t$), $\gG'\in G(\fQ)$ arithmetic, and let $X_{\gG'_N}$ be a
  modular subvariety of $X_{\gG'}$ for $N$ rational symmetric, conjugate to
  $N_{\bfb}$. Then $X_{\gG'_N}$ is a finite quotient of a product, and the
  set of $\gG'$-equivalence classes of such modular subvarieties forms a
  locus in $X_{\gG'}$ where the corresponding abelian varieties are
  isogenous to products, i.e., are not simple.
  If $N$ is $\gG$-integral, then $X_{\gG'_N}$ is a product, and
  the set of $\gG'$-equivalence classes of such modular subvarieties forms
  a locus in $X_{\gG'}$ where the corresponding abelian varieties split
  while preserving the endomorphisms.
\end{theorem}
{\bf Proof:} The first statement follows immediately from Lemma
\ref{L6.4.1}. By Lemma \ref{L6.4.2}, if $N$ is $\gG$-integral,
the discrete subgroup $\gG'_N$ is a
product, hence so is the quotient $X_{\gG'_N}$, giving the second
assertion. We know by the discussion above the moduli interpretation
upstairs in $\cD$, given in Proposition \ref{p6.4.1}. Since $X_{\gG'_N}$
itself is a product, it follows that the abelian varieties $A^n$ also split
as $A^n\cong A^q\times A^{n-q}$, where $\tau_q$, the modulus of $A^q$,
defines a point in one of the factors of $X_{\gG'_N}$, while $\tau_{n-q}$,
the modulus of $A^{n-q}$, defines a point in the second factor. That the
endomorphisms are preserved was shown above in the proof of
\ref{p6.4.1}. \ende
We leave it to the reader to derive the correct result for $b=t$.

Finally we briefly mention the moduli interpretation of incidence. For
this, recall that one has the Satake compactification and the (smooth
projective) toroidal compactifications. The former relate to degenerations
of the abelian varieties as follows. A quasi-abelian variety $A'$ is an
extension of an abelian variety by an algebraic torus
\begin{equation}\label{E161.2} 1 \lra (\fC^*)^h \lra A' \lra B \lra 0.
\end{equation}
Thus $A'$ is still an abelian group. Let $c$ denote the dimension of the
abelian variety $B$, $n=h+c$ the dimension of $A'$.
We now suppose that $X_{\gG}$ is a moduli space of {\sc Pel} structures,
and assume the notations used above in this case.
Let $F_b$ be a standard boundary component of the domain $\cD$, and
$W_i$ a boundary variety which is covered by $F_b$. Let $n=\dim(A)$ for
the abelian varieties parameterized by $X_{\gG}$, $m=\dim_D(V)$ so that
$n=mg,\ g=f, 4f$ and $2d^2f$ in the respective cases. Since $F_b$ has rank
$b$, it corresponds to a vector subspace $W\inn V$ with
$\dim_D(W)=\dim_D(V)-b=m-b\ (b=1,\ldots, s=\hbox{Witt index of the form},\
s\leq [{m\over 2}])$,
and hence to abelian varieties $B$ with $\dim(B)=(m-b)g$ and $g$
as above. We
observe that the sequence (\ref{E161.2}) is relevant here, with $h=bg,\
c=(m-b)g$. An extension as in
(\ref{E161.2}) is far from being unique, and the precise degenerations have
been constructed in many cases by utilizing methods from the theory of
mixed variations of Hodge structures, and this can be brought into relation
to the toroidal compactifications mentioned above, where the parameter
spaces of the degenerations are divisors on $\xg$. For our purposes
(\ref{E161.2}) is sufficient. We now consider an integral modular subvariety
$X_{\gG_N}$ incident with a boundary variety $W_i$. As we have seen above,
the abelian varieties parameterized by $X_{\gG_N}$ split, in this case as
\begin{equation}\label{E161.3} 0\lra A^h\lra  A' \lra A^c\lra 0,
\end{equation}
and the relation to (\ref{E161.2}) is obvious; the boundary varieties are
the loci in $X_{\gG_N}$ where the $A^h$ of (\ref{E161.3}) totally
degenerate.

\end{document}